\documentclass[10pt]{article}
\usepackage[T1]{fontenc}
\usepackage{amsmath,amssymb,amsthm,mathtools,microtype}
\usepackage[margin=1in]{geometry}
\usepackage{graphicx}
\usepackage{xcolor}

\usepackage[hidelinks]{hyperref}
\numberwithin{equation}{section}
\numberwithin{figure}{section}
\newtheorem{theorem}{Theorem}[section]
\newtheorem{corollary}[theorem]{Corollary}
\newtheorem{proposition}[theorem]{Proposition}
\newtheorem{lemma}[theorem]{Lemma}
\newtheorem{definition}[theorem]{Definition}
\newcommand{\be}{\begin{equation}}
\newcommand{\ee}{\end{equation}}
\newcommand{\bse}{\begin{subequations}}
\newcommand{\ese}{\end{subequations}}
\newcommand{\del}{\partial}
\newcommand{\R}{\mathbb R}
\newcommand{\Mext}{M_{\textbf{ext}}}
\newcommand{\diver}{\operatorname{div}}
\newcommand{\Ocal}{\mathcal O}

\title{\textbf{The positive mass theorem\\ under integral decay and weak regularity}}
\author{Philippe G. LeFloch\thanks{Laboratoire Jacques-Louis Lions, Sorbonne Universit\'e \& Centre National de la Recherche Scientifique, 4~Place Jussieu, 75252 Paris, France. Email: {\tt contact@philippelefloch.org}.
}}
\date{November 2026}
\begin{document}
\maketitle
\begin{abstract}
We prove the positive mass theorem for a class of asymptotically Euclidean, three-dimensional Riemannian metrics that we call $H^1$-tame. These metrics enjoy a \emph{natural integral decay} (for defining the ADM mass, without a pointwise rate), belong locally to $H^1$, satisfy uniform spectral bounds relative to a complete smooth metric $h$ that is Euclidean on the asymptotic end, and obey the regularity condition $g-h\in\dot H^1(M,h)$ in the homogeneous Sobolev space; in addition, we require that, locally, up to a conformal factor, the metric satisfies a small Sobolev multiplier bound; for instance, local Dini mean oscillation gives a sufficient criterion under nonnegative scalar curvature. 
In this class of metrics, we prove that, if the scalar curvature density is non-negative in the sense of distributions, the ADM mass is non-negative and vanishes exactly for Euclidean metrics. Our result \emph{almost} reaches the natural $H^1_{\mathrm{loc}}$ threshold, which allows one to define curvature. Our proof is based on new local Dirac estimates and a bounded Witten-spinor construction. We identify the physical Dirac domain under bi-Lipschitz $W^{2,2}_{\mathrm{loc}}$ coordinate changes and bounded $H^1_{\mathrm{loc}}$ conformal changes. 
The proposed construction turns a solution in the completed energy space into an actual $H^1_{\mathrm{loc}}$ spinor and makes its energy \emph{well defined} against a possibly \emph{singular} scalar curvature measure. We use the limit of \emph{dyadic annular averages} to define the mass in this weak class of metrics, and we prove that it is independent of the admissible Euclidean reference at infinity, and is finite exactly when the scalar curvature measure has finite mass on asymptotic ends.
\end{abstract}

\

{

\setcounter{secnumdepth}{2}
\setcounter{tocdepth}{1}

\tableofcontents
}

\

\newpage 

\section{Introduction}


\subsection{Purpose of this paper}
The positive mass theorem of Schoen and Yau~\cite{SchoenYau1979}, followed by Witten's spinorial proof~\cite{Witten1981}, relates the local sign of scalar curvature to a global invariant of an asymptotically Euclidean manifold. For a complete three-dimensional manifold without boundary and with nonnegative scalar curvature, the theorem asserts that the total ADM mass is nonnegative. Its rigidity statement identifies Euclidean space as the zero-mass case. These pioneering results have played a seminal role in the development of a broad range of topics in differential geometry and general relativity. Bartnik~\cite{Bartnik1986} established a geometric framework for the ADM mass and proved its invariance under coordinate changes satisfying suitable decay and integrability conditions. Beig and Chru\'sciel~\cite{BeigChrusciel1996} later studied the rigidity of the positive energy theorem and its relation to asymptotic Killing fields.

We seek to extend positivity and rigidity to metrics of \emph{weak regularity} when curvature can only be defined in the sense of distributions. Such metrics can describe geometric singularities. The curvature inequality must then be understood through its action on test functions. The geometry at infinity in the presence of asymptotically flat ends also requires care: standard approaches are based on a pointwise decay condition at some prescribed rate to ensure that the ADM mass is well defined, while in this paper we assume only \emph{integrability conditions} at asymptotic ends. Our aim is to retain the geometric conclusions of Schoen and Yau~\cite{SchoenYau1979} and Witten~\cite{Witten1981} under these weaker notions of curvature and decay.

Lee and LeFloch~\cite{LeeLeFloch} provided a foundation for this program by proving a distributional positive mass theorem for continuous $W^{1,3}$ Sobolev metrics. Their argument combined an integrated form of the Lichnerowicz identity with an annular definition of mass. We build on this framework and develop the analytic tools needed to approach the natural energy level $H^1_{\mathrm{loc}}$, where the first derivatives of the metric are only locally square-integrable.

We overcome this difficulty by separating the local regularity needed for the spinorial proof from the decay needed to define mass. The local condition is geometric and allows rough conformal changes. In the finite-mass case, a boundedness argument yields Witten spinors whose energy is well defined even when the scalar curvature measure is singular. At infinity, integral decay controls an average of the mass over annuli. At zero mass, the rigidity argument reconstructs global Euclidean coordinates and controls their regularity.

The regularity assumptions are intrinsic and are kept separate from the scalar curvature inequality. One calls the resulting class of metrics $H^1$-tame. We now define this class and states the positive mass and rigidity theorem under nonnegative distributional scalar curvature.

\subsection{The regularity class and the main theorem}
\paragraph{The \texorpdfstring{$H^1$}{H1}-tame regularity class.}
For the positivity and rigidity properties, we control the $L^2$ norm of the product of a test function with the covariant derivative of a local conformal representative by the Dirichlet norm of the test function. Given a smooth Riemannian metric $h$ and an open set $U$, the Sobolev multiplier norm of a scalar function $f\in L^2_{\mathrm{loc}}(U,h)$ is the best constant $A\in[0,+\infty]$ in
\be
 \|vf\|_{L^2(U,h)}\le A\|dv\|_{L^2(U,h)}
 \qquad\bigl(v\in C_c^\infty(U)\bigr).
\ee
The same definition applies to a locally square-integrable tensor field $F$, using its pointwise norm $|F|_h$. For a tensor $q\in H^1_{\mathrm{loc}}$, we apply it to the global covariant derivative $\nabla^h q$ and define the \underline{Sobolev multiplier norm of $\nabla^h q$} on $U$ by
\be
\label{eq:derivative-concentration-on-open}
 K_h(q;U)\coloneqq
 \sup_{\substack{v\in C_c^\infty(U)\\\|dv\|_{L^2(U,h)}>0}}
 \frac{\bigl\|v\,\nabla^h q\bigr\|_{L^2(U,h)}}{\|dv\|_{L^2(U,h)}}\in[0,+\infty].
\ee
All tensor norms and the volume density are those of $h$; explicitly,
\be
 \|v\nabla^h q\|_{L^2(U,h)}^2=\int_U|v|^2|\nabla^h q|_h^2\,dv_h,
 \qquad
 \|dv\|_{L^2(U,h)}^2=\int_U h^{ij}\del_i v\del_jv\,dv_h.
\ee
In coordinates, $\nabla^h q$ is computed with the Levi-Civita connection of $h$; when $h$ is Euclidean in Cartesian coordinates, it reduces to the ordinary derivatives of the coefficients of $q$. Thus $K_h(q;U)$ is the norm of multiplication by $\nabla^h q$ from the scalar Dirichlet space to $L^2$, in the standard terminology of Sobolev multipliers~\cite[Part~I]{MazyaShaposhnikova2009}. This is a geometric quantity for the pair $(q,h)$ and is unchanged by simultaneous smooth pullback of the tensors and the domain. For a bounded $H^1$ tensor on a relatively compact coordinate domain, its finiteness is independent of the smooth reference metric: changing the reference connection adds a bounded tensor, whose contribution is controlled by Poincar\'e's inequality. The numerical value can change with $h$, and a smallness threshold must be adjusted to the spectral bounds. The local condition in the definition below uses a flat reference after a conformal change of the metric.

\begin{definition}\label{def:H1-tame}
Let $M$ be a connected smooth three-manifold without boundary, let $g$ be a measurable symmetric two-tensor on $M$, and let $h$ be a smooth Riemannian metric on $M$. 

One calls $g$ an \underline{$H^1$-tame metric with reference $h$} if the following four conditions hold.
\begin{enumerate}
\item\label{item:tame-reference} \textbf{Euclidean reference at infinity.} The reference metric $h$ is complete, and there is a compact set $K\subset M$ such that $\Mext\coloneqq M\setminus K$ is a single end identified with an exterior Euclidean region, on which $h=g_{\text{Eucl}}$.
\item\label{item:tame-regularity} \textbf{Uniform spectral bounds.} The tensor $g$ belongs to $H^1_{\mathrm{loc}}(M;\operatorname{Sym}^2T^*M)$. The spectra of the positive $h$-self-adjoint endomorphisms $h^{-1}g$ lie almost everywhere in one fixed compact interval $J\subset(0,+\infty)$:
\be
 \operatorname{spec}\bigl((h^{-1}g)_x\bigr)\subset J
 \quad\text{for almost every }x\in M.
\ee
\item\label{item:tame-decay} \textbf{Integral decay.} The metric difference satisfies $g-h\in\dot H^1(M,h)$ (and hence belongs to $L^6(M,h)$).
\item\label{item:tame-local} \textbf{Local conformal regularity.} Every point has a neighborhood $U$ carrying a flat reference metric $f$ with bi-Lipschitz $W^{2,2}_{\mathrm{loc}}$ Euclidean coordinates onto a ball, and a positive function $w\in H^1_{\mathrm{loc}}(U)$ with $w,w^{-1}\in L^\infty(U)$, such that the conformal representative $\widehat g\coloneqq w^{-4}g$ satisfies, for some $\lambda\ge1$ and $0\le a<a_*(\lambda)$,
\be
\label{eq:tame-geometric-comparison}
 \begin{gathered}
 \operatorname{spec}(f^{-1}\widehat g)\subset[\lambda^{-1},\lambda]
 \quad\text{almost everywhere on }U,\\
 \|v\nabla^f\widehat g\|_{L^2(U,f)}\le a\|dv\|_{L^2(U,f)}
 \qquad\bigl(v\in C_c^\infty(U)\bigr).
 \end{gathered}
\ee
The positive threshold $a_*(\lambda)$ depends only on $\lambda$.\footnote{The local Dirac estimate allows the choice $a_*(\lambda)=c\lambda^{-3}$ for a sufficiently small universal constant $c>0$.}
\end{enumerate}
\end{definition}

The following observations about Definition~\ref{def:H1-tame} are in order.
\begin{enumerate}
\item \textbf{Intrinsic formulation.} For quantitative estimates, one chooses $C_{\mathrm{ell}}\ge1$ such that $J\subset[C_{\mathrm{ell}}^{-1},C_{\mathrm{ell}}]$. The uniform spectral bounds are equivalently expressed by $C_{\mathrm{ell}}^{-1}h\le g\le C_{\mathrm{ell}}h$, comparing lengths in every tangent direction independently of coordinates. All quantities in condition~\ref{item:tame-local} are formed from the tensor pair $(\widehat g,f)$: the spectral bounds compare the two metrics, and the covariant derivative, norms and volume density in the multiplier estimate are those of $f$. The flat reference has $H^1_{\mathrm{loc}}$ coefficients; its Levi-Civita connection is defined almost everywhere and belongs to $L^2_{\mathrm{loc}}$. The stated regularity of its Euclidean coordinates is part of the condition.

\item \textbf{Coordinate expression.} We use the coordinates $X:U\to B$, so $f=X^*g_{\text{Eucl}}$, and write $k=(X^{-1})^*\widehat g$. Then $k\in H^1_{\mathrm{loc}}(B)$ and
\be
\label{eq:tame-conformal-representation}
 g=w^4X^*k.
\ee
The weak tensor chain rule gives $\nabla^f\widehat g=X^*(\nabla^{g_{\text{Eucl}}}k)$. Change of variables identifies the multiplier norms, so condition~\ref{item:tame-local} reads equivalently
\be
\label{eq:tame-concentration-threshold}
 \lambda^{-1}g_{\text{Eucl}}\le k\le\lambda g_{\text{Eucl}},\qquad
 K_{g_{\text{Eucl}}}(k;B)<a_*(\lambda).
\ee
Compactly supported Lipschitz tests give the same norm by bounded smooth approximation. This justifies the equivalence for the Sobolev coordinate map $X$ and gives the coefficient formulation used in the Dirac estimates below.

\item \textbf{Role in the Dirac energy.} The multiplier condition addresses a specific difficulty in the Dirac energy. For a metric with only $H^1_{\mathrm{loc}}$ regularity, the spin connection belongs locally to $L^2$. An $H^1$ spinor belongs locally to $L^6$, and these two facts give a connection product only in $L^{3/2}$, whereas the Dirac energy requires it in $L^2$. In the canonical spin frame of $\widehat g$ built from the flat coordinates, the connection coefficients are bounded by $C(\lambda)|\nabla^f\widehat g|_f$. Condition~\ref{item:tame-local} therefore bounds multiplication by this connection from the Dirichlet space to $L^2$. Its small coefficient allows the connection contribution to be absorbed in the estimate for the first derivatives of the spinor. The conformal law handles the factor $w$ separately, and the coordinate and conformal transport identifies the physical Dirac domain in Section~\ref{sec:isotropic} below. The bounded Witten-spinor construction in Section~\ref{sec:witten-end-estimates} below then recovers actual $H^1_{\mathrm{loc}}$ spinors and their scalar-curvature energy.
\end{enumerate}

\paragraph{Euclidean references and scalar curvature.}
A neighborhood of infinity is an open set containing the complement of some compact set in $M$. One calls a Euclidean reference at infinity \underline{admissible} if it is given by a $W^{1,1}_{\mathrm{loc}}$ homeomorphism $X$ from an exterior region in $\R^3$ onto a neighborhood of infinity, taking infinity to infinity, and if the metric transferred to this exterior region, $g_X=X^*g$, is locally $H^1$ and satisfies uniform spectral bounds relative to $g_{\text{Eucl}}$. One also requires
\be
\label{eq:intro-end-integrability}
 \nabla^{g_{\text{Eucl}}}g_X\in L^2(g_{\text{Eucl}}),\qquad
 (g_X-g_{\text{Eucl}})/r\in L^2(g_{\text{Eucl}}),
\ee
where $r$ is the Euclidean radius. Both norms are taken over the entire exterior region of the reference. For a Sobolev map, the pullback is evaluated almost everywhere using the weak differential and a representative of $g$ satisfying the spectral bounds everywhere. The bounds in the two descriptions imply that both transition maps are locally Lipschitz, so the resulting tensor is independent of this choice of representative.

In the prescribed smooth identification at infinity, these integral conditions are equivalent to condition~\ref{item:tame-decay} of Definition~\ref{def:H1-tame}, by Proposition~\ref{prop:automatic-end-reference} below. We use the completion of compactly supported smooth tensor fields; its interpretation and decay consequences are explained after Proposition~\ref{prop:H1-tame-criteria} below.

The local $H^1$ regularity is a property of the tensor $g$: it can be expressed by $g,\nabla^h g\in L^2_{\mathrm{loc}}$ for any smooth reference metric $h$, and is independent of this auxiliary choice. Together with locally uniform spectral bounds, it suffices to define the scalar curvature density: the connection coefficients in a local frame lie in $L^2_{\mathrm{loc}}$, and in local coordinates the density is the divergence of a vector field in $L^2_{\mathrm{loc}}$ plus a term in $L^1_{\mathrm{loc}}$. We denote this distributional density by $R_g^{\mathrm w}$; for a smooth metric it equals $R_g\,dv_g$, where $dv_g$ denotes the Riemannian volume density of $(M,g)$, expressed by the volume form in oriented charts. When $R_g^{\mathrm w}$ is nonnegative, it defines a measure that can describe curvature concentrated on a surface. We give the formulation in Section~\ref{sec:weak-coframes} below.

\noindent\begin{minipage}{\textwidth}
\paragraph{The positive mass theorem.}
The curvature inequality is now imposed as a separate hypothesis on the regularity class just defined. The scalar curvature density and the dyadic annular average $m_{\text{dya}}(R)$ are defined in Sections~\ref{sec:weak-coframes} and~\ref{sec:mass-existence} below.

\begin{theorem}[Positive mass and rigidity for $H^1$-tame metrics]
\label{thm:small-concentration-mass}
Let $g$ be an $H^1$-tame metric with reference $h$, and suppose that its scalar curvature density $R_g^{\mathrm w}$ is nonnegative in the sense of distributions. Then the following conclusions hold.
\begin{enumerate}
\item The limit $m_{\rm ADM}\coloneqq\lim_{R\to\infty}m_{\text{dya}}(R)$ of the dyadic annular average lies in $[0,+\infty]$ and is finite exactly when $R_g^{\mathrm w}(\Mext)<+\infty$.
\item The mass is independent of the choice of admissible Euclidean reference at infinity.
\item The mass vanishes exactly when there is a global bi-Lipschitz homeomorphism $\Phi:(M,d_h)\to\R^3$ with $g=\Phi^*g_{\text{Eucl}}$ almost everywhere. The inverse of $\Phi$ gives an admissible Euclidean reference at infinity by restriction to an exterior region. Both $\Phi$ and $\Phi^{-1}$ belong to $W^{2,2}_{\mathrm{loc}}$ in smooth charts.
\end{enumerate}
\end{theorem}
\end{minipage}

Sections~\ref{sec:isotropic} and~\ref{sec:witten-end-estimates} below establish the physical Dirac domain and construct bounded Witten spinors with actual $H^1_{\mathrm{loc}}$ regularity. At zero mass, the resulting parallel spinor determines a closed orthonormal coframe; Section~\ref{sec:rigidity} below integrates it into global Euclidean coordinates with the stated second-derivative regularity.

\paragraph{Sufficient criteria for $H^1$-tame metrics.}
The following two criteria give ways to verify the local requirement in condition~\ref{item:tame-local} of Definition~\ref{def:H1-tame}. The first uses a fixed smooth reference metric, while the second provides coordinates in which the scalar curvature inequality controls the derivatives of the metric.

The \underline{limiting local multiplier norm at $x$} is obtained by restricting to $h$-geodesic balls $B_h(x,r)$ and letting the radius tend to zero; its supremum over $M$ gives a uniform local multiplier bound relative to $h$.
\be
\label{eq:derivative-concentration}
 \kappa_h(g;x)\coloneqq\lim_{r\downarrow0} K_h(g;B_h(x,r)),
 \qquad K_h(g)\coloneqq\sup_{x\in M}\kappa_h(g;x).
\ee
By monotonicity under restriction, the limit exists in $[0,+\infty]$ at every point. Both its finiteness and its vanishing are independent of the smooth reference metric: on shrinking balls, the bounded difference between reference connections contributes $\Ocal(r)$ to the multiplier norm. The multiplier inequality with an arbitrarily small coefficient of the Dirichlet energy in \eqref{eq:full-metric-derivative-form} below is equivalent to $K_h(g)=0$. We will prove the reference comparisons and the relation with Sobolev multipliers in Section~\ref{sec:sobolev-multiplier-concentration} below.

The \underline{volume-normalized metric} associated with $g$ is $b\coloneqq J_g^{-2/3}g$, where $J_g\coloneqq dv_g/dv_h=(\det(h^{-1}g))^{1/2}$ is the ratio of volume densities. Then $dv_b=dv_h$ and $g=w^4b$ with $w\coloneqq J_g^{1/6}$. For each tangent vector $Y$, the weak matrix chain rule gives
\be
\label{eq:trace-free-shape-derivative}
 \nabla^h_Y b=J_g^{-2/3}\left(\nabla^h_Y g-\tfrac13\operatorname{tr}_g(\nabla^h_Yg)\,g\right).
\ee
The multiplier condition can therefore be imposed on the trace-free part in the two metric indices, which describes changes of shape after the volume has been normalized.

The local multiplier condition also has a formulation in which its limiting numerical value is intrinsic. We choose any smooth connection $\nabla^0$ on $TM$, and let $\operatorname{tf}_g$ remove the $g$-trace in the two metric indices. We define
\be
\label{eq:metric-concentration}
 \mathcal K(g;x)\coloneqq\lim_{U\downarrow x}\sup_{0\ne v\in C_c^\infty(U)}
 \frac{\|v\,\operatorname{tf}_g(\nabla^0g)\|_{L^2(U,g)}}{\|dv\|_{L^2(U,g)}},
 \qquad \mathcal K(g)\coloneqq\sup_{x\in M}\mathcal K(g;x).
\ee
The limit uses any neighborhood basis at $x$, and all norms and volume densities are those of $g$. A change of smooth connection adds a tensor bounded in these norms on compact sets, whose contribution vanishes by the local Poincar\'e inequality. Consequently $\mathcal K(g;x)$ is independent of the auxiliary connection and invariant under smooth pullback. Under fixed spectral bounds relative to $h$, it is comparable to $\kappa_h(b;x)$; we prove this comparison in Section~\ref{sec:sobolev-multiplier-concentration} below. This intrinsic quantity expresses the local regularity of the conformal shape of the metric, while $\kappa_h(b;x)$ is a convenient version for quantitative estimates in a fixed reference.

\begin{definition}\label{def:tame-sufficient-criteria}
Let $g$ be a measurable Riemannian metric with $H^1_{\mathrm{loc}}$ regularity and locally uniform spectral bounds relative to a smooth reference metric.
\begin{enumerate}
\item Given a smooth reference metric $h$ and $\eta>0$, one says that $g$ has \underline{a small multiplier bound} on its volume-normalized derivative relative to $h$ at threshold $\eta$ at a point $x$ if its volume-normalized metric $b$ satisfies $\kappa_h(b;x)<\eta$. One uses the same term on a set when the inequality holds at every point of that set.
\item One says that $g$ has \underline{local Dini mean oscillation} if every point has a smooth coordinate neighborhood with the following property. In the chart image $U$, one writes $g=(g_{ij})$ for the coefficient matrix and $g_B$ for its coefficientwise average on a ball $B$. For every compact subball $V\subset U$, one can choose $0<s_V<\operatorname{dist}(V,\partial U)$ such that
\be
\label{eq:dini-mean-oscillation}
 \int_0^{s_V}\omega_{g,V}(s)\,\frac{ds}{s}<+\infty,\qquad
 \omega_{g,V}(s)\coloneqq\sup_{y\in V}|B_s(y)|^{-1}\int_{B_s(y)}
                  |g-g_{B_s(y)}|\,dx .
\ee
\end{enumerate}
\end{definition}

The Dini condition is independent of the choice of smooth coordinates: on compact subsets, a smooth change of coordinates compares oscillations on balls of comparable radii and adds a term of order $s$ to the modulus. Together with the reference independence of local multiplier finiteness and vanishing, this shows that the regularity properties are intrinsic. Quantitative smallness is expressed relative to $h$, with a threshold adapted to the spectral bounds.

\begin{figure}[htbp]
\centering
\includegraphics[width=0.78\linewidth]{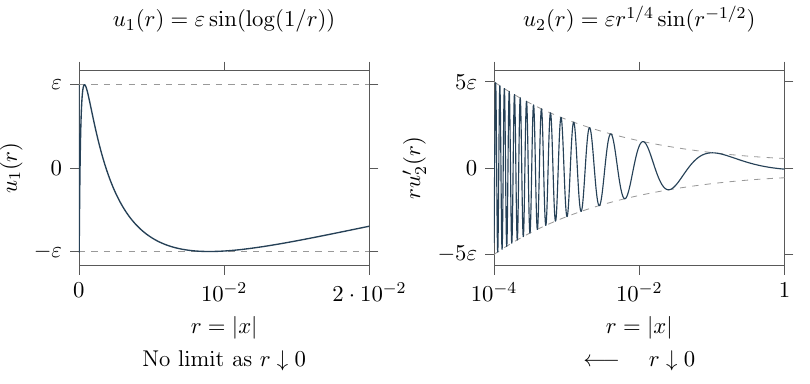}
\caption{The conditions that fail in the two examples. Left: $u_1(r)$ has oscillations of fixed amplitude accumulating at $r=0$, so it has no limit there and fails the Dini condition. Right: $r u_2^{\prime}(r)$ has growing amplitude, reflecting the failure of the derivative multiplier bound although $u_2$ has Dini mean oscillation. The radial scale is linear on the left and logarithmic on the right; in both panels $r\downarrow0$ toward the left.}
\label{fig:two-regularity-conditions}
\end{figure}

The two analytic conditions in Definition~\ref{def:tame-sufficient-criteria} are independent in a fixed smooth reference. Figure~\ref{fig:two-regularity-conditions} illustrates this distinction with radial functions on the three-dimensional Euclidean unit ball. For $r=|x|>0$ and $\varepsilon>0$, we set
\be
 u_1(x)\coloneqq\varepsilon\sin\bigl(\log(1/r)\bigr),\qquad
 u_2(x)\coloneqq\varepsilon r^{1/4}\sin(r^{-1/2}),
\ee
where we put $u_2(0)=0$ and may choose any value for $u_1(0)$; both functions are bounded and belong to $H^1(B_1)$. The three-dimensional Hardy inequality gives $K_{g_{\text{Eucl}}}(u_1;B_1)\le2\varepsilon$, while its mean oscillation on balls centered at the origin is bounded below by $c\varepsilon$; hence $u_1$ does not have Dini mean oscillation. In contrast, $u_2\in C^{0,1/6}$ and therefore has Dini mean oscillation, but $K_{g_{\text{Eucl}}}(u_2;B_\rho)=+\infty$ for every $0<\rho\le1$. Indeed, tests supported in $R<r<2R$ give a multiplier quotient bounded below by $c\varepsilon R^{-1/4}$. The metrics $\operatorname{diag}(e^{2u_i},e^{-u_i},e^{-u_i})$ have unit volume density and the same distinction between the two regularity conditions. The left panel displays $u_1$ itself to show its lack of a limit at the origin; the right panel displays the second of the following weighted derivatives:
\be
\label{eq:radial-example-derivatives}
 \begin{aligned}
 r u_1'(r)&=-\varepsilon\cos\bigl(\log(1/r)\bigr),\\
 r u_2'(r)&=\tfrac14\varepsilon r^{1/4}\sin(r^{-1/2})
             -\tfrac12\varepsilon r^{-1/4}\cos(r^{-1/2}).
 \end{aligned}
\ee
Both unweighted derivatives are square integrable against $r^2\,dr$ near zero. For $u_1$, scaling a ball about the origin gives a continuous periodic function of $\log(1/r)$ for its normalized mean oscillation; this function is strictly positive at every phase, so it has a positive minimum. For $u_2$, the amplitude bound and $|u_2'(r)|\le C\varepsilon r^{-5/4}$ give the H\"older exponent $1/6$ by separating radial increments smaller and larger than $r^{3/2}$. If $v_R(x)=\phi(|x|/R)$ with a fixed nonzero $\phi\in C_c^\infty((1,2))$, the squared Dirichlet norm is proportional to $R$, while the squared norm of $v_Rdu_2$ is bounded below by $c\varepsilon^2R^{1/2}$ for small $R$. This follows from the second formula above by integrating its rapidly oscillating squared cosine, whose weighted average tends to one half. The resulting quotient is at least $c\varepsilon R^{-1/4}$, which proves the asserted divergence of the multiplier norm. These examples concern the analytic conditions; the second local criterion in Proposition~\ref{prop:H1-tame-criteria} below also requires a nonnegative scalar curvature density.

\begin{proposition}[Sufficient conditions for $H^1$-tame metrics]
\label{prop:H1-tame-criteria}
Let $g$ be a measurable Riemannian metric with $H^1_{\mathrm{loc}}$ regularity and locally uniform spectral bounds relative to a smooth reference metric. One fixes a point $x\in M$. The representation and bounds required in condition~\ref{item:tame-local} of Definition~\ref{def:H1-tame} hold near $x$ if either of the following conditions holds at $x$.
\begin{enumerate}
\item In a neighborhood of $x$, the metric satisfies uniform spectral bounds relative to a smooth reference metric $h$, with $\operatorname{spec}(h^{-1}g)\subset[C_{\mathrm{ell}}^{-1},C_{\mathrm{ell}}]$ and $C_{\mathrm{ell}}\ge1$, and it satisfies the small multiplier condition relative to $h$ at $x$ for the threshold $\eta_*(C_{\mathrm{ell}})>0$. This threshold depends only on $C_{\mathrm{ell}}$.
\item The coefficients of $g$ have local Dini mean oscillation in a smooth chart about $x$, and $R_g^{\mathrm w}\ge0$ in that neighborhood.
\end{enumerate}
If conditions~\ref{item:tame-reference}--\ref{item:tame-decay} of Definition~\ref{def:H1-tame} hold and one of these local criteria holds at every point, then $g$ is an $H^1$-tame metric with reference $h$. The choice of criterion may vary from one neighborhood to another.
\end{proposition}
\begin{proof}
For the first criterion, we use smooth coordinates in which the metric $h$ agrees with $g_{\text{Eucl}}$ at the point, and shrink the ball so that $h$ satisfies $g_{\text{Eucl}}/2\le h\le2g_{\text{Eucl}}$. We take $X$ to be this chart, $k$ the coefficients of $b$, and $w=J_g^{1/6}$. The factor and its inverse are bounded and the Sobolev chain rule gives $w\in H^1_{\mathrm{loc}}$. In the chart, $|dk|\le C|\nabla^h b|_h+C_B$, where $C$ is uniform under the metric comparisons and $C_B$ is finite on a relatively compact subball. The bounds on $g$ and $h$ give $(2C_{\mathrm{ell}}^2)^{-1}g_{\text{Eucl}}\le k\le 2C_{\mathrm{ell}}^2 g_{\text{Eucl}}$. A coordinate ball of radius $r$ is contained in an $h$-ball of radius $Cr$. The comparison of norms, the definition of $\kappa_h$, and Poincar\'e's inequality therefore give $\||dk|\varphi\|_{L^2}\le(C_0\eta_*(C_{\mathrm{ell}})+C_1r)\|d\varphi\|_{L^2}$ for sufficiently small $r$. The constant $C_0$ depends only on the metric comparison bounds, and $C_1$ is finite on the chosen subball. We choose $\eta_*(C_{\mathrm{ell}})$ so that $C_0\eta_*(C_{\mathrm{ell}})<a_*(2C_{\mathrm{ell}}^2)/2$, then shrink $r$ to obtain the estimate in condition~\ref{item:tame-local} of Definition~\ref{def:H1-tame}.

For the second criterion, Proposition~\ref{thm:dini-harmonic-mass} below constructs a bi-Lipschitz $W^{2,2}_{\mathrm{loc}}$ harmonic chart in which the metric coefficients are continuous. We take $w=1$ and let $k$ be the metric in this chart. Its scalar density is nonnegative by covariance under the coordinate change. After a constant linear normalization at the point, Proposition~\ref{prop:harmonic-scalar-control} below bounds the multiplier coefficient by $C\|k-g_{\text{Eucl}}\|_{L^\infty}^{1/2}$, which tends to zero as the ball shrinks.
\end{proof}

\paragraph{The homogeneous Sobolev condition at infinity.}
We use the homogeneous Sobolev notation~\cite{BahouriCheminDanchin}. For symmetric two-tensors, $\dot H^1(M,h)$ is the completion of $C_c^\infty(M;\operatorname{Sym}^2T^*M)$ for the norm $\|\nabla^h q\|_{L^2(M,h)}$. The global Sobolev inequality and Kato's inequality give $\|q\|_{L^6(M,h)}\le C\|\nabla^h q\|_{L^2(M,h)}$, so this completion has a canonical image in $L^6(M,h)$. We identify its elements with these tensor fields. Thus membership implies $L^6$ integrability; for an actual tensor $q$, it is equivalent to $q\in L^6$ and $\nabla^h q\in L^2$. In particular, the condition on $g-h$ fixes the constant reference at infinity. After a cutoff at the inner boundary of the end, Hardy's inequality recovers the weighted term in~\eqref{eq:intro-end-integrability}. Proposition~\ref{prop:automatic-end-reference} below proves these characterizations and their independence of the smooth extension of $h$ across the compact core.

The integral conditions~\eqref{eq:intro-end-integrability} describe decay at the scale of $W^{1,2}_{-1/2}$ in dimension three. Under $R_g^{\mathrm w}\ge0$, Lemma~\ref{lem:dyadic-annular-mass} below gives finite negative variation of $m_{\text{dya}}$, finite total variation when its limit is finite, and agreement with the smooth annular expressions used in the proof.

Section~\ref{sec:coordinate-invariance} below proves mass invariance: the transition between admissible references and its inverse have bounded differentials and square-integrable second derivatives near infinity, and each differential differs from a constant orthogonal map by an $L^6$ field.

\subsection{Historical background}
\paragraph{Curvature and positive mass.}
We compare the curvature, decay, and regularity assumptions with earlier positive mass theorems.

The distributional connection and curvature framework used here was developed by LeFloch and Cristinel Mardare~\cite{LeFlochMardare}. Geroch and Traschen~\cite{GerochTraschen} had studied distributional curvature for locally bounded metrics and inverses with locally square-integrable weak first derivatives. Lott~\cite{LottRicciMeasure} defined a Ricci curvature measure for certain singular torsion-free connections through an integral formula. His later work on singular metrics distinguished the assumptions for distributional scalar curvature from those needed for a self-adjoint Dirac operator~\cite{Lott}. These distinctions are relevant here: the scalar curvature density is defined before the multiplier estimate used in the spinorial argument.

Arnowitt, Deser, and Misner~\cite{ADM1961} introduced the energy expression at spatial infinity. The analytic foundations of the spinorial proof were developed by Reula~\cite{ReulaWitten} and Parker and Taubes~\cite{ParkerTaubes1982}. Bartnik and Chru\'sciel~\cite{BartnikChruscielDirac} developed a boundary and regularity theory for Dirac-type operators with Sobolev coefficients, while Bartnik~\cite{Bartnik1986} studied invariance of the mass and the harmonic-coordinate and weighted elliptic framework.

Chru\'sciel~\cite{ChruscielSpatialInfinity} studied the ADM energy in a Hamiltonian formulation. His discussion of the critical rate $r^{-1/2}$ included flat metrics exhibiting coordinate dependence. His work on weaker boundary conditions also included examples of infinite mass~\cite{ChruscielPositiveEnergy}. The homogeneous Sobolev condition used here fixes the Euclidean reference through its $L^6$ representative, while the admissible coordinate class is specified by the same quadratic conditions in each reference.

Miao~\cite{MiaoCorners} proved a positive mass theorem for metrics with corners along a hypersurface. Lee~\cite{LeeLipschitz} considered Lipschitz metrics with small singular sets, and Grant and Tassotti~\cite{GrantTassotti} considered continuous metrics in $W^{2,3/2}_{\mathrm{loc}}$. In 2015, Lee and LeFloch~\cite{LeeLeFloch} proved positivity and rigidity for spin manifolds with $C^0\cap W^{1,3}_{\mathrm{loc}}$ metrics, nonnegative distributional scalar curvature and weighted $W^{1,3}_{-q}$ asymptotics, $q\ge1/2$.\footnote{The positive mass and rigidity result included the critical case $q=1/2$, whereas the conclusions on the annular limit and the finiteness criterion were proved under the stricter assumption $q>1/2$.} Their proof used an integrated Lichnerowicz identity and an annular definition of mass, and allowed a scalar curvature measure.

Hafemann~\cite{Hafemann} obtained positivity and metric-space rigidity for continuous $W^{1,3}_{\mathrm{loc}}$ metrics with a smooth asymptotically flat end. Lee, Litzinger, and Simon~\cite{LeeLitzingerSimon} proved positivity and rigidity for $L^\infty\cap W^{1,3}_{\mathrm{loc}}$ metrics with a smooth end, pointwise $C^2$ decay of order greater than $1/2$, and integrable scalar curvature there. Our local multiplier criterion applies directly to the Dirac equation at $H^1_{\mathrm{loc}}$ regularity, and the alternative Dini criterion obtains the required local estimate from harmonic coordinates and nonnegativity of the scalar curvature density.

Lundgren and Meco~\cite{LundgrenMeco} studied an annular definition of mass under weighted Sobolev asymptotics. They obtained a finite, cutoff-independent mass under weighted $W^{1,2}_{-1/2}$ asymptotics and scalar curvature represented by an $L^1$ function. They established invariance under weighted $W^{1,p}_{-\tau}$ asymptotics in both charts, with $p>3$ and $\tau>1/2$, once the annular mass is well defined in one chart. Here the scalar density may be a singular measure, the mass may be infinite, and the two quadratic integrability conditions determine the admissible coordinate class. The dyadic annular average fixes the averaging rule without an auxiliary cutoff.

\paragraph{Outline of this paper.}
Section~\ref{sec:weak-coframes} defines the scalar curvature density and proves its invariance, and Section~\ref{sec:mass-existence} constructs the dyadic annular mass. Section~\ref{sec:isotropic} establishes the local Dirac estimates and the integrated Lichnerowicz identity. Section~\ref{sec:witten-end-estimates} constructs bounded Witten spinors, identifies their energy with the ADM mass, and proves positivity. Section~\ref{sec:coordinate-invariance} compares Euclidean references at infinity and constructs the harmonic coordinates needed for the Dini criterion. Section~\ref{sec:rigidity} proves Euclidean rigidity and completes the main theorem.

\section{The distributional scalar curvature of \texorpdfstring{$H^1_{\mathrm{loc}}$}{H1loc} metrics}\label{sec:weak-coframes}


\subsection{Coframes and connection}
For $H^1_{\mathrm{loc}}$ metrics with locally uniform spectral bounds, we define the scalar curvature density using the distributional connection framework of LeFloch and Mardare~\cite{LeFlochMardare}. We compare its coframe expression, used in the spinorial proof, with the coordinate formula of Lee and LeFloch~\cite{LeeLeFloch}, which gives the current defining mass. The local construction uses a smooth reference metric $h$ and imposes no condition at infinity.

On an open set $U$, we choose a smooth $h$-orthonormal frame $(e_i)$ with dual coframe $(\eta^i)$, for example by orthonormalizing a coordinate frame. In this construction, the indices $i,j,k$ and $a,b,c$ run from $1$ to $3$ and distinguish the reference frame from the $g$-orthonormal frame, respectively. We assume that the measurable metric $g$ is positive definite almost everywhere and satisfies
\be
 g,\ g^{-1}\in L^\infty_{\mathrm{loc}}(U),\qquad g\in H^1_{\mathrm{loc}}(U).
\ee
We write $g_{ij}\coloneqq g(e_i,e_j)$ for the coefficients of $g$ in the reference frame and denote by $g^{1/2}$ the unique positive symmetric matrix whose square is $(g_{ij})$. This matrix represents the positive $h$-self-adjoint endomorphism $(h^{-1}g)^{1/2}$, and the one-forms $\theta^a\coloneqq(g^{1/2})_{ai}\eta^i$ form a $g$-orthonormal coframe with dual frame $(e_a)$. The local upper and lower bounds for the eigenvalues of $h^{-1}g$ allow the bounded Sobolev chain rule for the square root and the inverse, so the coframe and its inverse are locally bounded and locally $H^1$. Its exterior derivative is
\be
\label{eq:reference-coframe-exterior-derivative}
 d\theta^a=d(g^{1/2})_{ai}\wedge\eta^i+(g^{1/2})_{ai}\,d\eta^i,
\ee
in which both terms are locally in $L^2$. The second includes the commutators of the reference frame through $d\eta^i(e_j,e_k)=-\eta^i([e_j,e_k])$. We use a local orientation for each coframe calculation and prove agreement on overlaps below.

\begin{lemma}[Connection forms from the torsion equation]
Let $\theta=(\theta^a)$ be an oriented coframe on $U$ with locally bounded coefficients and inverse. Assume that $d\theta\in L^1_{\mathrm{loc}}(U)$. There is a unique locally integrable skew matrix of one-forms $\omega=(\omega^a{}_b)$ satisfying
\be
                    d\theta^a+\omega^a{}_b\wedge\theta^b=0.
\ee
On every relatively compact subdomain one has $C^{-1}|d\theta|\le|\omega|\le C|d\theta|$ almost everywhere, with $C$ depending only on the bounds for the coframe and its inverse on that subdomain.
\end{lemma}
\begin{proof}
At each point, we consider the linear torsion map from skew
matrix-valued one-forms to triples of two-forms, which is a map between spaces of dimension nine. We write an element of its kernel as $\sigma_{ab}=s_{abc}\theta^c$. Skew symmetry and vanishing torsion give
$s_{abc}=-s_{bac}$ and $s_{abc}=s_{acb}$. Repeatedly applying these
relations yields
\be
 s_{abc}=s_{acb}=-s_{cab}=-s_{cba}=s_{bca}=s_{bac}=-s_{abc}.
\ee
The kernel is therefore zero, and the inverse map depends continuously on the
coframe in the bounded invertible matrix class. Applying it almost
everywhere to $-d\theta$ constructs $\omega$ and gives the corresponding norm bounds.
The torsion equation then holds as an equality of integrable forms,
hence in the sense of distributions.
\end{proof}


\subsection{The scalar density and exterior products}
We write $g=\sum_a\theta^a\otimes\theta^a$ and denote its volume density by $dv_g$. In the dual frame $(e_a)$ we contract the connection forms to define a vector field $F$ and a scalar $T$, following Lott~\cite{Lott}. The product $\omega\wedge\omega$ uses matrix multiplication. Repeated frame and coordinate indices are summed: 
\be
 F=\omega^{ab}(e_b)e_a-\omega^{ab}(e_a)e_b,
 \qquad T=(\omega\wedge\omega)^{ab}(e_a,e_b).
\ee
We form the contraction $T$ before testing its integrability. The torsion equation gives $F\in L^1_{\mathrm{loc}}$ when the exterior derivative of the coframe satisfies $d\theta\in L^1_{\mathrm{loc}}$. The definition below also requires $T\in L^1_{\mathrm{loc}}$. For $H^1_{\mathrm{loc}}$ coframes the connection belongs to $L^2_{\mathrm{loc}}$, so both integrability requirements follow.
\begin{definition}
If $F,T\in L^1_{\mathrm{loc}}(U)$, the \underline{weak scalar curvature density} $R_g^{\mathrm w}$ in this coframe is defined by
\be
\label{eq:coframe-density-test}
 \langle R_g^{\mathrm w},\varphi\rangle
 =-\int_U d\varphi(F)\,dv_g-\int_U T\varphi\,dv_g,
 \qquad \varphi\in C_c^\infty(U).
\ee
One writes $R_g^{\mathrm w}\ge0$ if the pairing is nonnegative on every nonnegative test function, and denotes the resulting Radon measure by $\mu$.
\end{definition}

To express the scalar density in coordinates, we write $\theta^a=A_{ai}dx^i$ and $e_a=e_a^i\del_i$, using $g_{ij}=g(\del_i,\del_j)$ for the coordinate matrix of the metric. Its volume density is $dv_g=\sqrt{\det(g_{ij})}\,dx$. For the coordinate calculations, we abbreviate $\rho\coloneqq\sqrt{\det(g_{ij})}$ and choose the orientation so that $\rho=\det A$. In these coordinates, $F^i=(e_a^ie_b^j-e_a^je_b^i)\omega^{ab}_j$ and the formula becomes $R_g^{\mathrm w}=\operatorname{div}(\rho F)-\rho T\,dx$, where the divergence acts on the vector density $\rho F$. We prove invariance under changes of coframe below and recover $R_g\,dv_g$ for smooth metrics in Proposition~\ref{prop:critical-coordinate-scalar}.

The next identity computes a divergence from the exterior derivatives of the coframe and will be used both in changes of coframe and in the spinorial energy calculation.

\begin{lemma}[A divergence identity for exterior products]
For the oriented coframe above, let
$C_{ab}^{ij}=e_a^ie_b^j-e_a^je_b^i$ and define
\be
 B_{ab}^{i}=e_a^ie_c^j\omega^c{}_{bj}
 -e_b^ie_c^j\omega^c{}_{aj}
 +e_c^ie_b^j\omega^c{}_{aj}-e_c^ie_a^j\omega^c{}_{bj}.
\ee
Then $B\in L^1_{\mathrm{loc}}$ and
\be
 \rho C_{ab}^{ij}=\varepsilon_{abc}\varepsilon^{ijk}\theta^c_k,
 \qquad \del_i(\rho C_{ab}^{ij})=-\rho B_{ab}^j,
 \qquad \del_j(\rho B_{ab}^j)=0
\ee
in the sense of distributions. The hypotheses control the exterior derivative of the coframe.
\end{lemma}
\begin{proof}
The first formula expresses the duality between an oriented coframe and its dual bivectors, with the volume density included. On differentiation, its right side is
$\varepsilon_{abc}\varepsilon^{ijk}\del_i\theta^c_k$;
antisymmetry makes this one half of the corresponding contraction
of $(d\theta^c)_{ik}$, an integrable function. We substitute
\be
 (d\theta^c)_{ik}=-\omega^c{}_{di}\theta^d_k
                         +\omega^c{}_{dk}\theta^d_i.
\ee
This volume-duality identity and the skew symmetry of $\omega$
identify this contraction with $-\rho B_{ab}^j$, giving the second
formula with all its products integrable. Finally the contraction
$C_{ab}^{ij}$ is antisymmetric in $i,j$. Commutativity of
distributional derivatives gives its double divergence equal to
zero, proving the last identity. These calculations use only the
exterior derivative of the fixed coframe.
\end{proof}


\subsection{Invariance under changes of coframe}
\paragraph{Invariance of the scalar density.}
We first compare the densities under a smooth rotation, allowing the exterior derivative of the original coframe to be only locally integrable.

\begin{lemma}[Smooth rotations of a coframe]\label{lem:scalar-smooth-rotation}
Let $\theta$ be an oriented coframe on $U$ with locally bounded coefficients and inverse. Suppose $d\theta\in L^1_{\mathrm{loc}}(U)$ and the contraction $T$ in the scalar density belongs to $L^1_{\mathrm{loc}}(U)$. For a smooth map $R:U\to SO(3)$, set $\theta'=R\theta$. Then $T'\in L^1_{\mathrm{loc}}(U)$ and the scalar densities computed from $\theta$ and $\theta'$ agree.
\end{lemma}
\begin{proof}
We set $\Xi=R^{-1}dR$ and use the uniqueness of the torsion-free connection to obtain
$\omega'=R\omega R^{-1}-dR\,R^{-1}$. We express the contractions defining $F'$ and $T'$
in the original frame indices; the rotation factors cancel algebraically and
leave the connection $\omega-\Xi$. Thus, using the antisymmetric coefficients $C_{ab}^{ij}$ defined above,
\be
 \rho(F'^i-F^i)=\mathcal E^i\coloneqq -\rho C_{ab}^{ij}\Xi^{ab}_j,
 \qquad
 T'-T=B_{ab}^j\Xi^{ab}_j
       +C_{ab}^{ij}(\Xi_i\Xi_j)^{ab}.
\ee
The second identity follows by expanding the two cross terms and
the term quadratic in $\Xi$ in the defining contraction of $T$;
its cross coefficient is precisely the four terms defining $B$.
Each new term is locally integrable, because $B\in L^1_{\mathrm{loc}}$ and
$\Xi$ is smooth. The quadratic contraction is formed before the comparison of the integrals.

Differentiating $R^{-1}R=I$ gives the Maurer--Cartan identity for the rotation form $\Xi$.
\be
\label{eq:rotation-maurer-cartan}
 d\Xi+\Xi\wedge\Xi=0.
\ee
Combining \eqref{eq:rotation-maurer-cartan} with the divergence identity for $\rho C_{ab}^{ij}$ gives
\be
 \diver\mathcal E
 =\rho B_{ab}^j\Xi^{ab}_j
       -\rho C_{ab}^{ij}\del_i\Xi^{ab}_j
 =\rho B_{ab}^j\Xi^{ab}_j
       +\rho C_{ab}^{ij}(\Xi_i\Xi_j)^{ab}
 =\rho(T'-T).
\ee
The differentiation involves smooth coefficients of $\Xi$ and the $L^1_{\mathrm{loc}}$ divergence of $\rho C_{ab}^{ij}$ just established. Subtracting the two formulas for the scalar curvature density therefore
gives $(R_g^{\mathrm w})'-R_g^{\mathrm w}=\diver\mathcal E-\rho(T'-T)dx=0$ as a
distribution. This proves equality of the entire measure whenever
the distributions are measures.
\end{proof}

We now prove invariance under rotations of class $H^1_{\mathrm{loc}}$, which relate the coframes built from different smooth references for a metric of this regularity.

\begin{lemma}[Sobolev rotations of a coframe]\label{lem:scalar-sobolev-rotation}
Let $\theta$ be an oriented $H^1_{\mathrm{loc}}$ coframe on $U$ with locally bounded coefficients and inverse. If $R\in H^1_{\mathrm{loc}}(U,SO(3))$ and $\theta'=R\theta$, then the scalar density computed from $\theta'$ equals that computed from $\theta$.
\end{lemma}
\begin{proof}
We use the notation of Lemma~\ref{lem:scalar-smooth-rotation}. The Sobolev product rule gives $\theta'=R\theta\in H^1_{\mathrm{loc}}$ and $\Xi=R^tdR\in L^2_{\mathrm{loc}}$. Torsion uniqueness yields $\omega'=R\omega R^t-dR\,R^t\in L^2_{\mathrm{loc}}$. The algebraic transformation formulas above remain valid almost everywhere; their linear terms belong to $L^2_{\mathrm{loc}}$ and their quadratic terms to $L^1_{\mathrm{loc}}$. The identity~\eqref{eq:rotation-maurer-cartan} continues to hold in the sense of distributions. Indeed, we work on a relatively compact subpatch and convolve the entries of $R$. The uniformly bounded smooth approximations converge strongly in $H^1$ and almost everywhere. The products $R^t dR$ are obtained in $L^2$ and $dR^t\wedge dR$ in $L^1$, so the smooth identity gives $d(R^t dR)=dR^t\wedge dR$ distributionally. Differentiating $R^tR=I$ gives $dR^t=-\Xi R^t$, and thus $dR^t\wedge dR=-\Xi\wedge\Xi$.

For $A_{ab}^{ij}=\rho C_{ab}^{ij}$, the exterior product identity gives $A\in H^1_{\mathrm{loc}}\cap L^\infty_{\mathrm{loc}}$ and $\partial_iA_{ab}^{ij}=-\rho B_{ab}^j\in L^2_{\mathrm{loc}}$. For $\varphi\in C_c^\infty(U)$, we define the derivative pairing by $L^2$ duality:
\be
 \bigl\langle A_{ab}^{ij}\partial_i\Xi^{ab}_j,\varphi\bigr\rangle
 \coloneqq -\int\Xi^{ab}_j\partial_i(A_{ab}^{ij}\varphi)\,dx.
\ee
Antisymmetry in $i,j$ allows us to test \eqref{eq:rotation-maurer-cartan} against the bounded $H^1$ coefficient $A^{ij}\varphi$: we use bounded smooth test coefficients converging strongly in $H^1$ and almost everywhere, with dominated convergence for the $L^1$ quadratic term. Consequently
\be
 \operatorname{div}\mathcal E=\rho B_{ab}^j\Xi^{ab}_j+\rho C_{ab}^{ij}(\Xi_i\Xi_j)^{ab}=\rho(T'-T)
\ee
as distributions. This is the cancellation in the proof of Lemma~\ref{lem:scalar-smooth-rotation}, and it gives $(R_g^{\mathrm w})'=R_g^{\mathrm w}$. The limiting calculation approximates the test coefficients with $g$ fixed.
\end{proof}

A constant orthogonal reflection preserves the contractions defining $F$ and $T$. On each connected overlap the determinant of the Sobolev change of coframe belongs to $H^1_{\mathrm{loc}}$ and takes only the values $1$ and $-1$. Its derivative vanishes by differentiating its square. Thus its sign is constant on the overlap, and the rotation argument also applies when the local orientations differ. Proposition~\ref{prop:critical-coordinate-scalar} below identifies the resulting density with its expression in smooth coordinates.

\paragraph{The divergence identity for the frame.}
We now derive the divergence identity used in the Dirac argument.

\begin{lemma}[Exterior products and the frame flux]
If $\alpha,\beta$ are bounded differential forms with coefficients
in $H^1_{\mathrm{loc}}$, then
\be
 d(\alpha\wedge\beta)=d\alpha\wedge\beta
                       +(-1)^{\deg\alpha}\alpha\wedge d\beta
\ee
in the sense of distributions. Consequently, for the $H^1_{\mathrm{loc}}$ coframe above,
\be
\label{eq:frame-flux-new}
                  \del_i(\rho e_a^i)=-\tfrac12\rho F_a,
 \qquad F_a=\langle F,e_a\rangle_g.
\ee
\end{lemma}
\begin{proof}
The bounded Sobolev chain rule applied to $(f+g)^2$ and $(f-g)^2$
gives $d(fg)=f\,dg+g\,df$ for bounded $H^1_{\mathrm{loc}}$ scalar
coefficients. Applying this rule coefficient by coefficient and
using the antisymmetry of the exterior product gives the first
identity. We apply the product rule to the two-form
$\iota_{e_a}(\theta^1\wedge\theta^2\wedge\theta^3)
=\tfrac12\varepsilon_{abc}\theta^b\wedge\theta^c$.
Its coordinate coefficients are $\rho e_a^i$ under contraction with
$dx^1\wedge dx^2\wedge dx^3$. Substitution of the torsion equation
and skew symmetry of $\omega$ gives its exterior derivative equal
to $-\rho F_a\,dx^1\wedge dx^2\wedge dx^3/2$, proving
\eqref{eq:frame-flux-new}.
\end{proof}


\subsection{Coordinate formula and positive scalar curvature measures}
\paragraph{The coordinate expression of scalar curvature.}
We now identify the coframe expression with the coordinate scalar distribution. Under the same $H^1_{\mathrm{loc}}$ metric hypotheses, the sign of this distribution also gives a capacity property of its associated measure.

\begin{proposition}[Agreement with the coordinate scalar distribution]
\label{prop:critical-coordinate-scalar}
For $g\in H^1_{\mathrm{loc}}$ with $g,g^{-1}\in L^\infty_{\mathrm{loc}}$, the
scalar curvature density $R_g^{\mathrm w}$ computed in a canonical coframe agrees with the coordinate scalar density used by Lee and LeFloch \cite{LeeLeFloch} and by Lee, Litzinger, and Simon \cite{LeeLitzingerSimon}.
In a Euclidean background chart, one defines the Christoffel symbols
$\Gamma^k_{ij}$ by the weak first derivatives of $g$, and sets
\be
 V^k=g^{ij}\Gamma^k_{ij}-g^{ik}\Gamma^j_{ij},\qquad
 \mathcal Q=-(\del_k g^{ij})\Gamma^k_{ij}
       +(\del_jg^{ij})\Gamma^k_{ik}
       +g^{ij}(\Gamma^k_{k\ell}\Gamma^\ell_{ij}
                         -\Gamma^k_{j\ell}\Gamma^\ell_{ik}).
\ee
Then for every $\varphi\in C_c^\infty(U)$,
\be
 \langle R_g^{\mathrm w},\varphi\rangle
   =\int\{-V^k\del_k(\rho\varphi)
                         +\rho\mathcal Q\varphi\}\,dx.
\ee
The current $V$ belongs to $L^2_{\mathrm{loc}}$ and $\mathcal Q$ belongs to $L^1_{\mathrm{loc}}$. The pairings are finite, and both definitions give the same Radon measure when the scalar curvature density is nonnegative.
\end{proposition}
\begin{proof}
By coframe invariance, we may use $\theta=A\,dx$, with $A=\sqrt{(g_{ij})}$.
For a smooth metric, we write $\Gamma^i{}_j=\Gamma^i_{kj}dx^k$. The relation
$\omega=A\Gamma A^{-1}-dA\,A^{-1}$ gives
$d\omega+\omega\wedge\omega=A(d\Gamma+\Gamma\wedge\Gamma)A^{-1}$,
so the coordinate and orthonormal curvature contractions agree.
To check the sign in the coframe formula, we put $w_{ab,c}=\omega_{ab}(e_c)$ and
\be
 X=\sum_{a,b,c}w_{ac,a}w_{cb,b}
   =-\sum_b\Bigl(\sum_cw_{cb,c}\Bigr)^2,\qquad
 Y=\sum_{a,b,c}w_{ac,b}w_{cb,a}.
\ee
Then $T=X-Y$. The first two terms defining $B$ contract with $\omega$ to $-2X$, and the last two to $2Y$; hence $B_{ab}^j\omega_j^{ab}=-2T$.
The exterior product identity therefore gives
\be
 \rho e_a^ie_b^j(\partial_i\omega_j^{ab}-\partial_j\omega_i^{ab})
       =\operatorname{div}(\rho F)-2\rho T.
\ee
The quadratic curvature contraction contributes $\rho T$, leaving
$\operatorname{div}(\rho F)-\rho T$. In coordinates the same contraction is
\be
 \rho g^{ij}\bigl(\partial_k\Gamma^k_{ij}-\partial_j\Gamma^k_{ik}
       +\Gamma^k_{k\ell}\Gamma^\ell_{ij}
       -\Gamma^k_{j\ell}\Gamma^\ell_{ik}\bigr).
\ee
Integrating its derivative terms by parts, and separating the derivative of $g^{ij}$ from that of $\rho\varphi$, yields the stated formula with $V$ and $\mathcal Q$.

We now fix a compactly supported smooth test function and a slightly larger relatively compact chart. We convolve $g$ there with nonnegative smooth kernels. The resulting metrics $g_\varepsilon$ preserve uniform spectral bounds and converge to $g$ strongly in $H^1$ and, along a subsequence, almost everywhere. The matrix square root, inverse and determinant maps are smooth with bounded derivatives on this spectral range. Their compositions therefore converge strongly in $H^1$, remain uniformly bounded in $L^\infty$, and converge almost everywhere.
It follows that $\Gamma_\varepsilon,\omega_\varepsilon,F_\varepsilon,V_\varepsilon$ converge strongly in $L^2$. Every quadratic term in the two scalar formulas converges in $L^1$, including $T_\varepsilon$, $\mathcal Q_\varepsilon$ and $V_\varepsilon\,d\rho_\varepsilon$. Thus both pairings converge on the fixed test, and the smooth agreement passes to $g$. These same estimates give $V\in L^2_{\mathrm{loc}}$, $\mathcal Q\in L^1_{\mathrm{loc}}$ and finiteness of all pairings.

For a smooth background $h$, the calculation uses $\nabla^h$ and $\Gamma_g-\Gamma_h$, with the additional term $\operatorname{tr}_g\operatorname{Ric}_h$, as in~\cite{LeeLitzingerSimon}; its smooth coefficients cause no change in the convergence argument. 
Finally, pullback by a fixed smooth coordinate change preserves this local convergence. Passing to the limit in the smooth transformation law gives the absolute-Jacobian rule for a scalar density, including on overlaps with opposite orientations. The density is therefore defined on $M$, and its two expressions give the same Radon measure whenever it is nonnegative.
\end{proof}

\paragraph{The scalar measure and Sobolev capacity.}
The nonnegativity of the scalar distribution gives a measure that vanishes on sets of zero local $H^1$ capacity, a property that we will use to test the density against Sobolev functions.

\begin{lemma}[The scalar measure and Sobolev capacity]\label{lem:scalar-measure-capacity}
Let $g$ be a metric in $H^1_{\mathrm{loc}}$ with locally uniform spectral bounds and suppose $R_g^{\mathrm w}\ge0$. The corresponding Radon measure vanishes on every Borel set of zero local $H^1$ capacity. For a bounded $H^1$ function compactly supported in a chart, its quasi-continuous representative may be used to test the scalar density by the same formula as for smooth tests.
\end{lemma}
\begin{proof}
On a chart the scalar formula reads $\mu=\diver J-q\,dx$, with $J\coloneqq\rho F\in L^2_{\mathrm{loc}}$ and $q\coloneqq\rho T\in L^1_{\mathrm{loc}}$. We fix a compact set $L$ of zero $H^1$ capacity in the chart. We choose smooth compactly supported $0\le\eta_j\le1$ equal to one near $L$, with support in a fixed relatively compact subdomain and $\|\eta_j\|_{H^1}\to0$. Passing to a subsequence gives $\eta_j\to0$ almost everywhere. The nonnegativity of $\mu$ gives
\be
 0\le\mu(L)\le\int\eta_j\,d\mu
 =-\int J\cdot d\eta_j\,dx-\int q\eta_j\,dx\longrightarrow0.
\ee
The current term tends to zero by Cauchy--Schwarz and the last term by dominated convergence. Inner regularity and a countable chart cover give the conclusion for all Borel sets of zero local capacity.

We consider a bounded function $\varphi$ in $H^1$ with compact support in the chart. Smooth convolution after extension by zero gives uniformly bounded tests with support in a fixed compact subdomain that converge strongly in $H^1$. A subsequence converges outside a set of zero $H^1$ capacity to the quasi-continuous representative $\widetilde\varphi$. This follows from the capacity estimate for superlevel sets of successive differences and a subsequence with summable capacity bounds. The convergence is also $\mu$-almost everywhere by the first part. Thus the scalar pairing extends to
\be
\label{eq:scalar-density-sobolev-test}
 \int \widetilde\varphi\,d\mu=-\int J\cdot d\varphi\,dx-\int q\varphi\,dx .
\ee
The two sides converge by dominated convergence for the measure term and the $L^1$ term, and by $L^2$ duality for the current.
\end{proof}


\section{Definition and existence of the annular ADM mass}\label{sec:mass-existence}


\subsection{Annular mass and integral asymptotics}\label{sec:end-integrals}
\paragraph{The dyadic annular average.}
We define the mass at a finite scale by averaging the ADM expression on the annulus between that scale and twice that scale. The resulting dyadic annular average uses the prescribed Euclidean reference and its radius, while its limit will be independent of the admissible reference.

We identify $\Mext$ with $\{x\in\R^3:|x|>R_0\}$, write $r\coloneqq|x|$ and $A_R\coloneqq\{R<r<2R\}$, and denote the metric perturbation by $h\coloneqq g-g_{\text{Eucl}}$. All contractions and volume integrals in this paragraph use $g_{\text{Eucl}}$, with the same notation for the Euclidean tensor and its Cartesian coefficient matrix. For $g\in H^1_{\mathrm{loc}}(\Mext)$, the ADM expression is the locally square-integrable field with components
\be
\label{eq:adm-linear-expression}
 P_i(g)\coloneqq\sum_{j=1}^3(\del_jg_{ij}-\del_i g_{jj}).
\ee

\begin{definition}
For $g\in H^1_{\mathrm{loc}}(\Mext)$ and $R>R_0$, the \underline{dyadic annular average} of the ADM expression is
\be
\label{eq:dyadic-annular-mass}
 m_{\text{dya}}(R)\coloneqq\frac1{16\pi R}\int_{A_R}P(g)\cdot dr\,dx.
\ee
When its limit as $R\to+\infty$ exists in $\R\cup\{+\infty\}$, one calls it the \underline{ADM mass} in this exterior reference and denotes it by $m_{\rm ADM}$.
\end{definition}
By polar integration, the volume integral is the arithmetic mean over $R<t<2R$ of the spherical ADM integrals, which are defined for almost every $t$. The volume integral itself is defined for every $R>R_0$; at a fixed scale it depends on the chosen exterior reference. The ratio of the two radii is fixed at two, as in a dyadic decomposition, while $R$ may vary continuously.

For the distributional calculations, we also use smooth radial cutoff functions. Such a function $0\le\chi\le1$ equals one on $B_1$ and zero outside $B_2$, and its rescaling $\chi_R(x)\coloneqq\chi(x/R)$ has gradient vanishing outside $R\le r\le2R$. For locally integrable coefficients $h$, the corresponding smooth annular expression is
\be
\label{eq:smooth-annular-mass}
 16\pi m_\chi(R)\coloneqq\int_{\Mext}\bigl(h_{ij}\del_i\del_j\chi_R-(\operatorname{tr}_{g_{\text{Eucl}}}h)\Delta\chi_R\bigr)\,dx.
\ee
For $h\in W^{1,1}_{\mathrm{loc}}(\Mext)$, integration by parts gives
\be
\label{eq:annular-mass-definition}
 m_\chi(R)=-\frac1{16\pi}\int_{\Mext}P(g)\cdot d\chi_R\,dx.
\ee
Under the integral and curvature hypotheses used below, Lemma~\ref{lem:dyadic-annular-mass} below identifies the limits of these smooth expressions with the limit of $m_{\text{dya}}(R)$. Thus smooth cutoff functions serve as test functions in the proofs, while the dyadic annular average fixes the rule for computing the mass. We will compare admissible references in Section~\ref{sec:coordinate-invariance} below.

\paragraph{Consequences of integral decay.}
We now relate the two quadratic integrals on the end to a global Sobolev condition on the metric difference, and derive the estimates needed for the mass. One defines the homogeneous Sobolev space by completion of compactly supported smooth fields in the gradient norm, with its canonical Sobolev embedding into $L^6$, as explained in the introduction. 
The $L^6$ condition below is a consequence of this completion and provides a useful characterization of its elements. The two exterior integrals are
\be
\label{eq:end-energy}
 \int_{r>R_0}\left(|\del h|^2+r^{-2}|h|^2\right)\,dx<+\infty,
 \qquad h=g-g_{\text{Eucl}}.
\ee
The following equivalences concern a sufficiently far exterior region; for a field assumed only to belong to $H^1_{\mathrm{loc}}$ on an abstract exterior, they do not impose an additional condition near its inner boundary.

\begin{proposition}
\label{prop:automatic-end-reference}
Let $h\in H^1_{\mathrm{loc}}(\Mext)$ be a symmetric two-tensor in the Euclidean exterior reference.
\begin{enumerate}
\item One fixes $R_1>R_0$ and chooses a smooth function $\zeta$ that vanishes for $r\le R_1$ and equals one for $r\ge2R_1$. The following three conditions are equivalent:
\be
\label{eq:end-homogeneous-equivalences}
 \begin{aligned}
 &\int_{r>R_1}\bigl(|dh|^2+r^{-2}|h|^2\bigr)\,dx<+\infty;\\
 &\text{the zero extension of }\zeta h\text{ belongs to }\dot H^1(\R^3);\\
 &dh\in L^2(\{r>R_1\}),\qquad h\in L^6(\{r>R_1\}).
 \end{aligned}
\ee
The choice of $\zeta$ does not affect these conditions.
\item Under these equivalent conditions, the following two quadratic annular errors tend to zero:
\be
\label{eq:quadratic-annular-errors}
 \frac1R\int_{A_R}|h|\,|\del h|\,dx\longrightarrow0,
 \qquad \frac1{R^2}\int_{A_R}|h|^2\,dx\longrightarrow0.
\ee
\item Suppose that $\Mext$ is the single Euclidean end of a connected smooth three-manifold $M$ with compact core, that $g\in H^1_{\mathrm{loc}}(M)$, and that $h=g-g_{\text{Eucl}}$ on this end. For every complete smooth reference metric $g_0$ equal to $g_{\text{Eucl}}$ there, the conditions in the first item are equivalent to
\be
\label{eq:global-homogeneous-regularity}
 \begin{aligned}
 &g-g_0\in\dot H^1(M,g_0),\\
 &g-g_0\in L^6(M,g_0),\qquad\nabla^{g_0}(g-g_0)\in L^2(M,g_0).
 \end{aligned}
\ee
Here the fields are symmetric two-tensors, with tensor norms and volume density induced by $g_0$. This condition is independent of the smooth extension of the prescribed Euclidean reference across the compact core.
\end{enumerate}
\end{proposition}
The global condition concerns the difference between the physical and reference metrics. Ordinary $H^1$ regularity of this difference would also require an unweighted $L^2$ bound, which need not hold for perturbations with a nonzero $r^{-1}$ leading term; the homogeneous condition instead has the decay expressed by \eqref{eq:end-energy}.

\begin{proof}
\bse
{\bf 1.} \emph{Sobolev and Hardy inequalities on the end.} We extend $H\coloneqq\zeta h$ by zero to $\R^3$. The compact transition region contributes finite norms by local $H^1$ regularity. Under the first condition in \eqref{eq:end-homogeneous-equivalences}, the field $H$ is locally $H^1$ and has square-integrable gradient. For the radial cutoff functions above, we set $H_R\coloneqq\chi_RH$, whose gradient error satisfies
\be
 \|d(H_R-H)\|_2\le\|dH\|_{2,\{r>R\}}+
 C\left(\int_{A_R}r^{-2}|H|^2\,dx\right)^{1/2}\longrightarrow0.
\ee
Each $H_R$ admits compact smooth approximation in the gradient norm. The Euclidean Sobolev inequality therefore makes these approximations Cauchy in $L^6$, and their limit agrees with $H$ on compact sets. This proves the homogeneous Sobolev condition and the $L^6$ condition.

Conversely, if $dh\in L^2$ and $h\in L^6$ on the tail, then $dH\in L^2(\R^3)$ and the error from differentiating the cutoff is bounded by
\be
 R^{-1}\|H\|_{L^2(A_R)}\le C\|H\|_{L^6(A_R)}\longrightarrow0.
\ee
The same compact approximation proves $H\in\dot H^1(\R^3)$, and Hardy's inequality on its approximants gives
\be
\label{eq:exterior-hardy-normalization}
 \int_{\R^3}r^{-2}|H|^2\,dx\le4\int_{\R^3}|dH|^2\,dx.
\ee
This recovers the weighted term on the tail. Membership in the homogeneous Sobolev completion also implies these two gradient and integrability conditions by its Sobolev embedding, so all three conditions are equivalent. Changing $\zeta$ changes $H$ only by a compactly supported $H^1$ tensor and leaves the conditions unchanged.

\ese

\vskip.2cm \noindent
{\bf 2.} \emph{The quadratic annular errors.} Since $r/R\le2$ on $A_R$, the Cauchy--Schwarz inequality gives
\be
 \frac1R\int_{A_R}|h|\,|\del h|\,dx
 \le2\left(\int_{A_R}r^{-2}|h|^2\,dx\right)^{1/2}
       \left(\int_{A_R}|\del h|^2\,dx\right)^{1/2}\longrightarrow0.
\ee
The second error in \eqref{eq:quadratic-annular-errors} tends to zero because $R^{-2}\le4r^{-2}$ on the same annulus and the weighted integral has a vanishing tail.

\vskip.2cm \noindent
{\bf 3.} \emph{The global tensor condition.} For a complete smooth reference $g_0$ with the stated geometry, the scalar Sobolev inequality proved in Section~\ref{sec:witten-end-estimates} below and Kato's inequality give
\be
 \|q\|_{L^6(M,g_0)}\le C\|\nabla^{g_0}q\|_{L^2(M,g_0)}
 \qquad\bigl(q\in C_c^\infty(M;\operatorname{Sym}^2T^*M)\bigr).
\ee
Thus the homogeneous completion has a representative in $L^6$, with its covariant derivative in $L^2$. Conversely, if these two norms are finite for a locally $H^1$ tensor $q$, radial cutoffs on the end satisfy
\be
 \|\nabla^{g_0}(\chi_Rq-q)\|_2
 \le\|\nabla^{g_0}q\|_{2,\{r>R\}}+CR^{-1}\|q\|_{2,A_R}
 \le\|\nabla^{g_0}q\|_{2,\{r>R\}}+C\|q\|_{6,A_R}\longrightarrow0.
\ee
Compact smooth approximation proves membership in the completion. We apply this characterization to $q=g-g_0$. The compact core contributes finite norms by local $H^1$ regularity, while on the end the reference connection and volume are Euclidean. The first part of the proof and Hardy's inequality therefore give the equivalence with the two exterior integrals. Two smooth reference metrics with the same Euclidean end differ only on a compact set, where their norms are equivalent and the difference of their connections is bounded; their difference as tensors is itself smooth and compactly supported. Consequently the global condition is unchanged.
\end{proof}


\subsection{Negative variation and existence of the mass at infinity}
The extended-limit and finiteness conclusions below have a precedent
in Lee--LeFloch \cite[Proposition 2.4]{LeeLeFloch}, whose proof explicitly
used the monotonicity of the scalar curvature measure tested against increasing
cutoffs. The following estimate bounds the negative variation under the two quadratic integrability conditions on the asymptotic end. The positive part of the scalar curvature measure may have infinite mass.

\begin{proposition}[Finite negative variation of the annular mass]
\label{prop:annular-negative-variation}
On $\Mext=\{r>R_0\}$ let $g\in H^1_{\mathrm{loc}}$ with $g,g^{-1}\in L^\infty$. Suppose the distributional scalar density is a signed Radon measure $\mu$ with $\mu^-(\{r>R_1\})<+\infty$ for some $R_1>R_0$. Assume, with $h=g-g_{\text{Eucl}}$,
\be
 \int_{\Mext}\bigl(|dg|^2+r^{-2}|h|^2\bigr)dx<+\infty.
\ee
One fixes a smooth radially nonincreasing cutoff $\chi$, equal to one on
$B_1$ and zero outside $B_2$, and put
$m_\chi(R)=-(16\pi)^{-1}\int_{\Mext} P(g)\cdot d\chi_R$.
For sufficiently large $R_1$ the function $m_\chi$ is the sum of locally absolutely continuous functions $A$ and $B$ on $[R_1,\infty)$ satisfying
\be
 m_\chi=A+B,\qquad A'\ge0,\qquad
 \int_{R_1}^{\infty}|B'(R)|\,dR<+\infty.
\ee
Consequently $m_\chi$ has finite negative variation and a limit in
$\R\cup\{+\infty\}$. Its limit is finite if and only if
$\mu^+(\{r>R_1\})<+\infty$. For a smooth end these two conditions on $\mu^\pm$ amount to integrability of $(R_g)_\pm$ near infinity. In particular, finiteness of the total scalar variation at infinity follows
if $\liminf_{R\to\infty}m_\chi(R)<+\infty$. In the finite case
$m_\chi$ has finite total variation for every smooth radial cutoff function. The limit is independent of the radial cutoff function, even when that cutoff function is not monotone.
\end{proposition}
\begin{proof}
In the coordinate formula of Proposition~\ref{prop:critical-coordinate-scalar}, we expand $d(\rho\varphi)$ for a smooth test function $\varphi$. With $\mathcal V=\rho V$, this gives
\be
 \mu=\diver\mathcal V+Q\,dx,\qquad Q=\rho\mathcal Q-V\cdot d\rho,\qquad |Q|\le C|dg|^2 .
\ee
We choose a smooth cutoff $0\le\zeta\le1$ vanishing on a collar of the inner boundary and equal to one far out. We extend the current $J=\zeta\mathcal V$ by zero to $\R^3$. Its divergence satisfies
\be
 \diver J=\zeta\mu^++\nu,\qquad
 \nu=(d\zeta\cdot\mathcal V-\zeta Q)\,dx-\zeta\mu^- .
\ee
The term with $d\zeta$ is supported in a compact shell, $Q\in L^1(\Mext)$, and the cutoff may be chosen to vanish on $\{r\le R_1\}$, where finiteness of $\mu^-$ was not assumed. Thus $\nu$ has finite total variation. We enlarge the starting radius $R_1$ for the flux estimates so that the support of $d\chi_R$ lies where $\zeta=1$ for every $R\ge R_1$.
Testing the distributional identity for $\operatorname{div}J$ gives
\be
 m_\chi(R)=\frac1{16\pi}\int\chi_R\zeta\,d\mu^+
       +\frac1{16\pi}\int\chi_R\,d\nu+e(R),\qquad
 e(R)=-\frac1{16\pi}\int_{\Mext}(P(g)-\mathcal V)\cdot d\chi_R.
\ee
On each compact interval of radii, $\partial_R\chi_R$ is uniformly bounded and supported on a fixed compact shell. Both measure integrals can therefore be differentiated under the integral sign, even if the radial measure has atoms. The first term defines $A$, which is nondecreasing. The variation of
the second term is at most $(16\pi)^{-1}|\nu|(\R^3)$, since each
$\chi_R(x)$ increases with $R$. To control the error $e(R)$,
we differentiate only its cutoff factor. The bounds on the derivatives of the cutoff give
\be
 |e'(R)|\le C R^{-2}\int_{R<r<2R}|h||dg|\,dx,
 \qquad
 \int_{R_1}^{\infty}|e'(R)|\,dR
 \le C\int_{\Mext} r^{-1}|h||dg|\,dx
 \le C\|h/r\|_2\|dg\|_2.
\ee
The middle inequality follows by integrating $R^{-2}$ over
$r/2<R<r$ and applying Tonelli's theorem. These estimates also
justify local absolute continuity and differentiation under the
integral. Moreover,
\be
 |e(R)|\le C\|h/r\|_{L^2(R<r<2R)}
                  \|dg\|_{L^2(R<r<2R)}\longrightarrow0.
\ee
Thus $B=(16\pi)^{-1}\int\chi_R\,d\nu+e(R)$ converges to $(16\pi)^{-1}\nu(\R^3)$ and its total variation is at most $(16\pi)^{-1}|\nu|(\R^3)+C\|h/r\|_{L^2(\Mext)}\|dg\|_{L^2(\Mext)}$. Monotone convergence gives
$A(R)\to(16\pi)^{-1}\int\zeta\,d\mu^+$, possibly $+\infty$.
Local finiteness of the measure on compact annuli gives the stated finiteness criterion. For a smooth metric $d\mu=\rho R_g\,dx$, with $\rho$ bounded above and below. Finally
$(m_\chi')_-\le|B'|$; if $A$ is bounded, its total variation is
finite as well. For any radial cutoff function, not necessarily monotone, $0\le\chi_R\le1$ and $\chi_R\to1$ pointwise. Dominated convergence gives the same limit if the positive scalar term is integrable; otherwise Fatou's lemma shows that it tends to $+\infty$. The finite signed remainder converges by dominated convergence and $e(R)\to0$ as before. This proves independence of the cutoff function also for an infinite limit. For a finite limit both measures in the decomposition are finite. Writing $\chi(x)=\chi_0(|x|)$, the change of variable $s=|x|/R$ gives $\int_{R_1}^\infty |\partial_R\chi_R(x)|\,dR\le\int_1^2|\chi_0'(s)|\,ds$ uniformly in $x$. Tonelli and the bound on $e'$ give finite total variation without any monotonicity of $\chi_0$.
\end{proof}

When the scalar curvature measure is nonnegative, an upper bound along a sequence of annular masses thus implies that this measure has finite mass on the end. We will establish the nonnegativity of the limiting mass in Section~\ref{sec:mass-identity} below.


\subsection{One-sided traces on coordinate spheres}
We write $S_t=\{|x|=t\}$, with outward Euclidean unit normal $n$ and Euclidean area measure $dS$. We use the density $\mathcal V$ from the scalar curvature formula $\mu=\operatorname{div}\mathcal V+Q\,dx$. The total flux is defined by one-sided annular limits at every radius. The surface integral in the definition below is first taken at almost every radius by polar integration.

\begin{lemma}[Spherical traces of the scalar curvature current]\label{prop:sphere-trace}
On $\Mext$, assume only $\mathcal V\in L^1_{\mathrm{loc}}(\Mext)$, $Q\in L^1_{\mathrm{loc}}(\Mext)$ and that $\mu=\operatorname{div}\mathcal V+Q\,dx$, tested on radial functions, defines a locally finite signed measure on $(R_0,\infty)$. One denotes this radial measure by $\lambda$, so $\int\psi\,d\lambda=\langle\mu,\psi(|x|)\rangle$. One fixes $s>R_0$ and sets, for almost every $t$,
\be
 f(t)=\frac1{16\pi}\int_{S_t}\mathcal V\cdot n\,dS,\qquad
 q(t)=\int_{S_t}Q\,dS,\qquad
 H(t)=f(t)+\frac1{16\pi}\int_s^t q(u)\,du .
\ee
Then $f$ and $H$ have unique left-continuous representatives $f_-,H_-$ in $\mathrm{BV}_{\mathrm{loc}}((R_0,\infty))$. For $R_0<a<b$ they satisfy
\be
\label{eq:spherical-trace-balance}
 16\pi\big(H_-(b)-H_-(a)\big)=\lambda([a,b)),\qquad
 16\pi\big(f_+(a)-f_-(a)\big)=\lambda(\{a\}),
\ee
where $f_+$ is the right-continuous representative. In particular $H_-$ is nondecreasing when $\lambda\ge0$. For every $a>R_0$ the two traces of $f$ are the limits of its inner and outer annular averages. If $\mu$ itself is a measure, the jump of $f$ equals $\mu(S_a)/(16\pi)$.
\end{lemma}
\begin{proof}
Polar integration gives $f,q\in L^1_{\mathrm{loc}}$ on the radius interval. Testing with $\psi(r)$, for $\psi\in C_c^\infty((R_0,\infty))$, gives $16\pi Df=\lambda-q\,dt$ and $16\pi DH=\lambda$ as one-dimensional distributions. A left-continuous primitive of $\lambda$ is given, up to a constant, by $\lambda([s,t))$ for $t\ge s$ and $-\lambda([t,s))$ for $t<s$. Its derivative is $\lambda$ by Fubini, so it agrees with $16\pi H$ up to a constant almost everywhere. This gives both identities, since the primitive of $q$ is continuous. Left continuity and almost-everywhere equality determine the representative uniquely. At every radius, including an atom, polar integration gives
\be
 f_-(a)=\lim_{\epsilon\downarrow0}
 \frac1{16\pi\epsilon}\int_{a-\epsilon<r<a}\mathcal V\cdot n\,dx ,
 \qquad
 f_+(a)=\lim_{\epsilon\downarrow0}
 \frac1{16\pi\epsilon}\int_{a<r<a+\epsilon}\mathcal V\cdot n\,dx .
\ee
Indeed these expressions average a function having the indicated one-sided limits. The same conclusion holds for nonnegative averaging kernels of integral one supported on the respective shrinking intervals.
\end{proof}

If $\mu$ is a Radon measure on the end, we can also test the flux on each sphere against $\eta\in C^\infty(S^2)$. We extend $\eta$ by $\eta(x)=\eta(x/r)$ and apply Lemma~\ref{prop:sphere-trace} to the current $\eta\mathcal V$. The product rule gives
\be
 \operatorname{div}(\eta\mathcal V)
   =\eta\mu+(\mathcal V\cdot\nabla\eta-\eta Q)\,dx .
\ee
Thus the inner and outer traces $T_t^\pm(\eta)$ of the normalized flux exist at every radius and satisfy
$16\pi(T_t^+-T_t^-)(\eta)=\int_{S_t}\eta\,d\mu$.
They define distributions of order at most one on $S^2$. Indeed on any compact interval of radii with $t$ in its interior, the value of a one-sided representative is bounded by its mean absolute value plus its variation. Polar integration bounds these two terms by $C\|\eta\|_{C^1}$, using the local masses of $|\mu|$, $|Q|\,dx$ and $|\mathcal V|\,dx$. Taking $\eta=1$ recovers the total trace.

For $g\in H^1_{\mathrm{loc}}$ with $g,g^{-1}\in L^\infty_{\mathrm{loc}}$, the coordinate scalar formula gives the local integrability of $\mathcal V$ and $Q$ required above. Under the hypotheses of Proposition~\ref{prop:annular-negative-variation}, $q\in L^1$ near infinity and $\lambda^-$ has finite mass there. Thus $f_-(r)$ has a limit in $\R\cup\{+\infty\}$. The same limit is obtained from the right-continuous convention. At finite radii, nonnegative scalar curvature makes the corrected flux $H$ nondecreasing.

\begin{lemma}\label{lem:dyadic-annular-mass}
Under the hypotheses of Proposition~\ref{prop:annular-negative-variation}, the dyadic annular average $m_{\text{dya}}$ is locally absolutely continuous for $R>R_0$ and has finite negative variation near infinity. It has the same limit as the spherical current traces $f_\pm$ and the smooth expressions $m_\chi$, including the possible value $+\infty$. This common limit is finite exactly when the positive part of the scalar curvature measure has finite mass near infinity. In that case $m_{\text{dya}}$ has finite total variation near infinity.
\end{lemma}
\begin{proof}
The dyadic annular average has the same variation properties. For almost every $t$, we define
\be
 E(t)\coloneqq\frac1{16\pi}\int_{S_t}(P(g)-\mathcal V)\cdot n\,dS.
\ee
The quadratic end integrals give $E\in L^1(dt/t)$ near infinity. By polar integration, we can write
\be
 m_{\text{dya}}(R)=B(R)+A(R),\qquad
 B(R)\coloneqq\int_1^2f(Rs)\,ds,\qquad
 A(R)\coloneqq\frac1R\int_R^{2R}E(t)\,dt.
\ee
Both averages are locally absolutely continuous, and differentiation of the error average gives
\be
 A'(R)=\frac{2E(2R)-E(R)}R-\frac1{R^2}\int_R^{2R}E(t)\,dt
 \quad\text{almost everywhere}.
\ee
For every sufficiently large $S$, Tonelli's theorem and the substitutions $t=2R$ and $R<t<2R$ give
\be
 \int_S^{\infty}|A'(R)|\,dR
 \le4\int_S^{\infty}|E(t)|\,\frac{dt}{t}<+\infty,
 \qquad
 |A(R)|\le2\int_R^{2R}|E(t)|\,\frac{dt}{t}\longrightarrow0.
\ee
The negative variation of $B$ is bounded by that of $f$ on the same tail: for every increasing partition of radii, the negative part of each increment of $B$ is at most the average of the negative parts of the corresponding increments of $f(Rs)$. Integration over $1<s<2$ and then passage to the supremum over partitions give the bound. The same argument with absolute values applies to total variation. Since $16\pi Df=\lambda-q\,dt$, the finite negative part of the scalar measure and $q\in L^1$ give finite negative variation of $f$, and hence of $m_{\text{dya}}$. If the limit of $f$ is finite, $f$ has finite total variation, so the same holds for $m_{\text{dya}}$. Finally $B(R)$ converges to the limit of $f_-$ and $A(R)\to0$, proving the same ADM limit and finiteness criterion, including the value $+\infty$.
For a radially nonincreasing smooth cutoff, polar integration expresses the scalar-current part of $m_\chi$ as an average of $f$ against $-\partial_r\chi_R$, a nonnegative weight of integral one supported in $[R,2R]$. It converges to the same extended limit, while the quadratic annular error tends to zero. Proposition~\ref{prop:annular-negative-variation} gives independence of the smooth cutoff and the finiteness criterion for the scalar measure. Thus these smooth test expressions and the dyadic annular average define the same $m_{\rm ADM}$.
\end{proof}


\section{Dirac estimates under Sobolev multiplier bounds}\label{sec:isotropic}


\subsection{Interior Dirac regularity for a critical Sobolev metric}
Throughout Section~\ref{sec:isotropic}, $g\in H^1_{\mathrm{loc}}$ and its coefficients and inverse are locally bounded. The integrated identity holds for bounded, compactly supported $H^1$ spinors; extending it to arbitrary $H^1$ spinors, proving local Dirac regularity and identifying the global operator domain use the additional multiplier estimates stated in the respective results. The scalar curvature density is assumed to be nonnegative in the sense of distributions only in the statements that require this hypothesis.

We first prove local $H^1$ regularity for the adjoint Dirac equation with $L^2$ source, under $W^{1,3}_{\mathrm{loc}}$ regularity or a small multiplier bound on the metric derivative. These estimates are independent of the sign of scalar curvature. We then establish the scalar trace and operator domain needed for the global mass identity.

On each oriented chart, we use the spin frame associated with the positive symmetric coframe of Section~\ref{sec:weak-coframes}. The Clifford action $c$ is skew-Hermitian and satisfies
$c(\alpha)c(\beta)+c(\beta)c(\alpha)=-2g^{-1}(\alpha,\beta)I$.
The torsion-free connection induces the spin connection $\nabla^g$, whose matrix in a local spin frame is denoted by $\Omega$. With $c^a\coloneqq c(\theta^a)$, the Dirac operator is
\be
 D_g\psi\coloneqq\sum_{a=1}^3c(\theta^a)\nabla^g_{e_a}\psi=c^a e_a^j(\partial_j\psi+\Omega_j\psi).
\ee
We abbreviate $\nabla^g$ and $D_g$ to $\nabla$ and $D$ when the metric is fixed.

The Hermitian fiber norm is the standard norm on $\mathbb C^2$ in an orthonormal spin frame. Changes of orthonormal spin frame act unitarily, so $|\psi|$, the $L^p(dv_g)$ norms, and the $L^2$ norms of $\nabla\psi$ and $D\psi$ are intrinsic. In contrast, $d\psi$ and the ordinary $H^1$ norm of its coefficients refer to the specified frame. An arbitrary $H^1$ change of spin frame need not preserve this coefficient space; the multiplier estimates below justify the frame changes used here. If $dv_g=\rho\,dx$, the map $\psi\mapsto\rho^{1/2}\psi$ is unitary from physical $L^2(dv_g)$ to coordinate $L^2(dx)$.

Since $\Omega\in L^2_{\mathrm{loc}}$, in this spin frame we interpret $D\psi=f$ in the sense of distributions through the adjoint identity
\be
 \int_U\langle \psi,D\phi\rangle\,dv_g
   =\int_U\langle f,\phi\rangle\,dv_g
 \qquad (\phi\in C_c^\infty(U;\mathbb C^2)).
\ee
Here $D\phi\in L^2$ for every test, so the equation makes sense for $\psi,f\in L^2_{\mathrm{loc}}$.

\begin{proposition}[Interior regularity for a critical Sobolev metric]
\label{thm:critical-Sobolev-Dirac-entry}
On a three-dimensional chart let $g\in W^{1,3}_{\mathrm{loc}}$ with
$g,g^{-1}\in L^\infty_{\mathrm{loc}}$. In the canonical spin frame the spin connection belongs to $L^3_{\mathrm{loc}}$. If $\psi,f\in L^2_{\mathrm{loc}}$
satisfy $D\psi=f$ in the sense of distributions, through the adjoint identity against compactly supported smooth
coefficient spinors, then $\psi\in H^1_{\mathrm{loc}}$ and
$\Omega \psi\in L^2_{\mathrm{loc}}$. On sufficiently small nested balls,
\be
 \|\psi\|_{H^1(B')}+\|\Omega \psi\|_{L^2(B')}
       \le C_{B',B}\bigl(\|f\|_{L^2(B)}+\|\psi\|_{L^2(B)}\bigr).
\ee
\end{proposition}
\begin{proof}
{\bf 1. \sl The comparison operator.} We first estimate the Dirac operator of a comparison metric on $\R^3$. We assume that $k$ has uniform spectral bounds relative to $g_{\text{Eucl}}$ on $\R^3$, is constant outside a ball, and
$\|dk\|_3\le\epsilon$. Conjugate its physical Dirac by the unitary
volume map $\phi\mapsto\rho_k^{1/2}\phi$ into coordinate $L^2$.
The resulting operator has the form $S=a^j\del_j+b$, with
\bse
\be
 (a^i)^*a^j+(a^j)^*a^i=2k^{ij}I,\qquad
 \|da\|_3+\|b\|_3\le C\|dk\|_3.
\ee
The constants depend only on the spectral bounds for $k$. These statements follow
from the Sobolev chain rule for the square root, inverse and volume
factor, and the algebraic formula for the spin connection. We establish symmetry on spinors with compactly supported smooth coefficients. We write $D=a^j\del_j+P$ in the physical coefficient frame,
where $a^j$ is skew-Hermitian and $P=a^j\Omega_j$. Compatibility
of the spin connection with Clifford multiplication, together with
the torsion and frame-flux identities, gives in the sense of distributions
\be
 \del_j(\rho a^j)=\rho(P-P^*).
\ee
Indeed $a^j=e_a^j c^a$, $c^{a*}=-c^a$ and $\Omega_j^*=-\Omega_j$.
The Clifford relations and the connection convention give
\be
 P-P^*=[a^j,\Omega_j]=e_a^j\omega_{abj}c^b,
 \qquad F_b=-2e_a^j\omega_{abj}.
\ee
The frame-flux identity gives
$\partial_j(\rho e_b^j)=\rho e_a^j\omega_{abj}$, proving the
displayed symmetry relation after multiplication by $c^b$.
All terms in this equality belong to $L^3_{\mathrm{loc}}$ here, and to
$L^2_{\mathrm{loc}}$ under the multiplier hypothesis on the metric derivative below. It follows by applying the
weak product rule to the bounded $H^1$ coframe and the smooth
constant Clifford matrices; no derivative of the distributional solution is
involved. We set $S=\rho^{1/2}D\rho^{-1/2}=a^j\del_j+b$, with
\be
 b=P-\tfrac12a^j\del_j\log\rho,\qquad
 b-b^*=\del_j a^j.
\ee
\ese
The second identity uses the bounded Sobolev chain rule. For compact
smooth coefficient fields, integration by parts now proves
$(S\psi,\phi)_{L^2(dx)}=(\psi,S\phi)_{L^2(dx)}$. This also identifies the
conjugated operator obtained from physical $H^1$ tests: multiplication
by $\rho^{\pm1/2}$ is bounded on $H^1_0$, by Sobolev--H\"older
here and by the multiplier estimate for the metric derivative in Proposition~\ref{thm:H1-multiplier-Dirac-entry} below.

\vskip.2cm \noindent
{\bf 2. \sl The graph estimate.} For such a coefficient field $\phi$, decompose
$(a^i)^*a^j=k^{ij}I+C^{ij}$, where $C^{ij}=-C^{ji}$.
Integration by parts in the antisymmetric term gives
\bse
\be
 \int\langle\del_i \phi,C^{ij}\del_j\phi\rangle
       =-\int\langle \phi,(\del_iC^{ij})\del_j\phi\rangle.
\ee
The second derivatives cancel by antisymmetry, and the formula is valid directly for $C\in W^{1,3}_{\mathrm{loc}}\cap L^\infty$ and compactly supported smooth $\phi$.
Its absolute value is bounded by
$C\|dk\|_3\|\phi\|_6\|d\phi\|_2$. Sobolev and
$\|b\phi\|_2\le\|b\|_3\|\phi\|_6$ therefore show, for small
$\epsilon$,
\be
 c\|d\phi\|_2\le\|S\phi\|_2\le C\|d\phi\|_2.
\ee
\ese
Density extends these bounds and symmetry to domain $H^1(\R^3)$.

\vskip.2cm \noindent
{\bf 3. \sl The method of continuity.} To prove surjectivity rather than assume adjoint regularity, join
$k$ to its exterior constant value $k_0$ by
$k_t=(1-t)k_0+tk$. Ellipticity stays uniform and
$\|dk_t\|_3\le\epsilon$. The corresponding volume-conjugated
operators $S_t:H^1\to L^2$ are continuous in operator norm:
their principal coefficients are continuous in $L^\infty$, and
their lower-order coefficients are continuous in $L^3$ by the
chain rule and dominated convergence. Symmetry gives
\be
 \|(S_t\pm i)\phi\|_2^2=\|S_t\phi\|_2^2+\|\phi\|_2^2,
 \qquad
 \|\phi\|_{H^1}\le C\|(S_t\pm i)\phi\|_2.
\ee
At $t=0$ these maps are isomorphisms by Fourier transformation.
The method of continuity now gives isomorphisms for all $t$:
invertibility is open by a Neumann series. If $t_j\to t$ and $(S_{t_j}\pm i)\phi_j=f$, the uniform estimate bounds $\phi_j$ in $H^1$. Applied to $\phi_j-\phi_k$, it bounds their difference by $C\|(S_{t_j}-S_{t_k})\phi_k\|_2$, which tends to zero by operator-norm continuity. The limit solves $(S_t\pm i)\phi=f$ and is unique by the same estimate.
Thus $S$ is self-adjoint on $H^1$, which is therefore also the domain of its adjoint.

\vskip.2cm \noindent
{\bf 4. \sl Localization.} We now construct a comparison metric equal to $g$ on the active ball. On $B_{2r}$ we put $k_0=g_{B_{2r}}$ and choose a smooth cutoff
$\chi$ with $0\le\chi\le1$, equal to one on $B_r$ and zero near $\del B_{2r}$. The eigenvalues of the average $k_0$ lie in the same spectral interval as those of $g$.
We extend
\bse
\be
 k=k_0+\chi(g-k_0)
\ee
to the constant matrix $k_0$ outside $B_{2r}$. Convexity preserves the spectral bounds. Poincar\'e at exponent three gives
\be
 \|dk\|_{L^3(\R^3)}
 \le C\bigl(\|dg\|_{L^3(B_{2r})}
             +r^{-1}\|g-k_0\|_{L^3(B_{2r})}\bigr)
 \le C\|dg\|_{L^3(B_{2r})}.
\ee
\ese
Absolute continuity makes this smaller than the preceding threshold
for small $r$. The extension agrees with $g$ on $B_r$.

Finally, localize the equation in the sense of distributions with a compactly supported smooth
cutoff $\eta$ supported in $B_r$. The volume-conjugated field
$\rho^{1/2}\eta \psi$ lies in the adjoint domain of $S$, with image
$\rho^{1/2}(\eta f+[D,\eta]\psi)\in L^2$. For this assertion the
adjoint identity extends to compactly supported $H^1$ test spinors, since
$D:H^1_0\to L^2$ is bounded; multiplication by the bounded
$W^{1,3}$ volume factor preserves $H^1$. Self-adjointness gives
$\rho^{1/2}\eta \psi\in H^1$, hence $\eta \psi\in H^1$ by the same
multiplier property. Its graph estimate gives the displayed local
bound. Since $\Omega\in L^3$ and $\eta \psi\in L^6$, the actual
connection product is in $L^2$ with the asserted estimate.
\end{proof}


\subsection{Interior regularity from a multiplier estimate on the metric derivative}\label{sec:dirac-form-local}
\paragraph{Control of multiplication by the connection.}
For a metric in $W^{1,3}_{\mathrm{loc}}$, the inclusion $L^3\cdot L^6\subset L^2$ controls multiplication of a spinor by the lower-order coefficient of the Dirac operator. We now impose this multiplication estimate directly on the metric derivative.
\begin{proposition}[Interior regularity from a multiplier estimate on the metric derivative]
\label{thm:H1-multiplier-Dirac-entry}
Let $g\in H^1_{\mathrm{loc}}$ with $g,g^{-1}\in L^\infty_{\mathrm{loc}}$ on a chart.
Suppose that on each relatively compact ball $U$ and for every
$\varepsilon>0$ there is $C_{\varepsilon,U}<+\infty$ such that
\be
\label{eq:full-metric-derivative-form}
 \int_U|dg|^2|v|^2dx
 \le\varepsilon\int_U|dv|^2dx
                +C_{\varepsilon,U}\int_U|v|^2dx
 \qquad(v\in C_c^\infty(U)).
\ee
Then every spinor satisfying $D\psi=f$ in the sense of distributions through the adjoint identity, with $\psi,f\in L^2_{\mathrm{loc}}$,
has $\psi\in H^1_{\mathrm{loc}}$ and actual connection product
$\Omega \psi\in L^2_{\mathrm{loc}}$. On nested relatively compact coordinate neighborhoods, the same estimate for the graph norm as in
Proposition~\ref{thm:critical-Sobolev-Dirac-entry} holds. All metric
derivatives enter the single norm $|dg|$; no direction is singled out.
\end{proposition}
\begin{proof}
{\bf 1. \sl Extension from the active ball.} We construct the comparison metric using the multiplier inequality for $|dg|$. We choose nested balls and a compact
smooth cutoff $\chi$ equal to one on the active ball. We take $0\le\chi\le1$ and choose a constant matrix $k_0$ whose eigenvalues lie in the same spectral interval as those of $g$. We extend
$k=k_0+\chi(g-k_0)$ to $\R^3$. Convexity preserves both spectral bounds. Its derivative is
$dk=\chi\,dg+(g-k_0)d\chi$. Applying
\eqref{eq:full-metric-derivative-form} to $\chi \phi$ and using
boundedness of $g$ shows that for every $a>0$,
\be
 \||dk|\phi\|_2\le a\|d\phi\|_2+B_a\|\phi\|_2
                    \qquad(\phi\in H^1(\R^3)),
\ee
with $B_a$ finite. The cutoff derivative contributes to the zeroth-order coefficient $B_a$.

\vskip.2cm \noindent
{\bf 2. \sl A uniform graph estimate.} We connect $k$ to $k_0$ by $k_t=(1-t)k_0+tk$ and conjugate the physical
Dirac to coordinate $L^2$. We write $S_t=a_t^j\del_j+b_t$.
We use the Clifford relations, the weak chain rule, and the symmetry identity established in the proof of Proposition~\ref{thm:critical-Sobolev-Dirac-entry} above. Uniformly in $t$,
$|da_t|+|b_t|\le C|dk|$. In the antisymmetric integration by parts
of that proof, the principal error is bounded by
\bse
\be
 C\||dk|\phi\|_2\|d\phi\|_2
       \le Ca\|d\phi\|_2^2+CB_a\|\phi\|_2\|d\phi\|_2.
\ee
The lower-order term obeys
$\|b_t\phi\|_2\le C(a\|d\phi\|_2+B_a\|\phi\|_2)$.
To make the absorption explicit, we put $P=\|d\phi\|_2$, $H=\|\phi\|_2$ and $N=\|S_t\phi\|_2$, and denote the uniform principal ellipticity constant by $\lambda>0$. The two estimates give
\be
 \lambda P^2
 \le 2N^2+C(a^2+a)P^2+CB_a^2H^2+CB_aHP.
\ee
We choose $a$ so $C(a^2+a)\le\lambda/4$ and apply Young's inequality
to the last term with derivative coefficient $\lambda/4$.
Then $\lambda P^2/2\le2N^2+C_\lambda B_a^2H^2$.
Thus, uniformly in $t$,
\be
 \|\phi\|_{H^1}\le C_a(\|S_t\phi\|_2+\|\phi\|_2),\qquad
 \|S_t\phi\|_2\le C_a\|\phi\|_{H^1}.
\ee
The integrations initially use smooth $\phi$, bounded $H^1$ principal
coefficients and their actual $L^2$ derivatives. The multiplier estimate for the metric derivative
extends every resulting product continuously to $H^1$.

Symmetry and the last display give
$\|\phi\|_{H^1}\le C_a\|(S_t\pm i)\phi\|_2$.
The maps $t\mapsto S_t$ are continuous from $H^1$ to $L^2$:
principal differences are bounded by $C|t-s|$, and lower-order
differences by $C|t-s|\,|dk|$. To see the latter, we write
$b_t=t\mathcal B(k_t)[dk]$, where $\mathcal B(A)$ is linear in
the derivative argument and smooth in the matrix $A$ on the set determined by its spectral bounds. The matrix $k-k_0$ is bounded, so $\mathcal B$ and its first
derivative are uniformly bounded along the path. The multiplier estimate for the metric derivative then
gives
\be
 \|(S_t-S_s)\phi\|_2\le C|t-s|(\|d\phi\|_2+\|\phi\|_2).
\ee
The method of continuity used
above therefore proves that $S_1\pm i$ are isomorphisms and
that its adjoint domain is $H^1$. For completeness, surjectivity
is open by the inverse mapping theorem and the operator-norm
continuity just proved. It is closed as follows. If $t_j\to t$ and
$(S_{t_j}\pm i)\phi_j=f$, the uniform inverse bound controls $\phi_j$ in
$H^1$. Applying the same bound to $\phi_j-\phi_k$ gives
\be
 \|\phi_j-\phi_k\|_{H^1}
 \le C\|(S_{t_j}-S_{t_k})\phi_k\|_2\longrightarrow0.
\ee
\ese
Its limit solves $(S_t\pm i)\phi=f$. At $t=0$ both inverses are
provided by the constant-coefficient Fourier transform. Finally,
if $w$ belongs to the adjoint domain, we solve
$(S_1-i)\phi=(S_1^*-i)w$ with $\phi\in H^1$. The difference $w-\phi$
is orthogonal to the range of $S_1+i$, hence vanishes. This proves
the asserted equality of domains without assuming regularity of $w$.

\vskip.2cm \noindent
{\bf 3. \sl The localized adjoint equation.} Multiplication by the bounded volume factors $\rho^{\pm1/2}$
preserves $H^1$, since their derivatives are bounded by $C|dg|$
and satisfy the same multiplier estimate. Localizing the original equation in the sense of distributions thus places $\rho^{1/2}\eta \psi$ in this adjoint domain,
with $L^2$ image $\rho^{1/2}(\eta f+[D,\eta]\psi)$.
The adjoint identity defining the equation in the sense of distributions extends from smooth tests to
compactly supported $H^1$ test spinors, because $D:H^1_0\to L^2$ is bounded and
$\psi,f\in L^2$. For a compact smooth coefficient field $\phi$ in the
comparison chart, we use the test
$\tau=\eta\rho^{-1/2}\phi\in H^1_0$. The Sobolev product rule and
the definition of $S$ yield in $L^2$
\bse
\be
 D\tau=\rho^{-1/2}\bigl(\eta S\phi+c(d\eta)\phi\bigr).
\ee
Using the equation in the sense of distributions and the skew-Hermiticity of $c(d\eta)$ gives
\be
 (\rho^{1/2}\eta \psi,S\phi)_{L^2(dx)}
 =\bigl(\rho^{1/2}(\eta f+c(d\eta)\psi),\phi\bigr)_{L^2(dx)}.
\ee
\ese
Since the right-hand side defines a bounded linear functional on $L^2$, the localized spinor belongs to the adjoint domain, with image $\rho^{1/2}(\eta f+c(d\eta)\psi)$. The graph bound gives
$\eta \psi\in H^1$ with its quantitative local
estimate. Finally $|\Omega|\le C|dg|$ and the multiplier estimate for $|dg|$
identify the actual product $\Omega\eta \psi\in L^2$. 
\end{proof}

\paragraph{A fixed small coefficient in the interior estimate.}
Suppose that $\lambda^{-1} g_{\text{Eucl}}\le g\le\lambda g_{\text{Eucl}}$ for some $\lambda\ge1$. There is a universal constant $c>0$ such that the local Dirac regularity proved above holds under the single estimate
\be
 \||dg|v\|_2\le a\|dv\|_2+b\|v\|_2,
 \qquad 0\le a<c\lambda^{-3},\quad 0\le b<+\infty,
 \qquad v\in C_c^\infty.
\ee
\begin{proof}[Proof of the fixed-threshold Dirac estimate]
The cutoff extension $k=k_0+\chi(g-k_0)$ preserves the spectral bounds. Applying the estimate to $\chi \phi$ gives $\||dk|\phi\|_2\le a\|d\phi\|_2+B\|\phi\|_2$ for a finite $B$: the terms containing $d\chi$ enter only in $B$. The same spectral bounds hold along $k_t=(1-t)k_0+tk$.

We now determine how the estimates in the preceding proof depend on $\lambda$. For a metric matrix $q$ whose eigenvalues satisfy the preceding spectral bounds, the Sylvester equation for the derivative of its square root gives $|d(q^{1/2})|\le C\lambda^{1/2}|dq|$. Differentiating its inverse gives $|d(q^{-1/2})|\le C\lambda^{3/2}|dq|$. If $a^j$ are the principal Dirac coefficients, then $|a^j|\le C\lambda^{1/2}$ and the derivative of the antisymmetric part of $(a^i)^*a^j$ is bounded by $C\lambda^2|dq|$. Evaluating $d\theta=d(q^{1/2})\wedge dx$ on two unit frame vectors costs at most a factor $\lambda$. The torsion equation therefore bounds the connection coefficients in that frame by $C\lambda^{3/2}|dq|$. For $\rho=\sqrt{\det q}$ we have $|d\log\rho|\le C\lambda|dq|$, so the volume-conjugation term satisfies the same bound $C\lambda^{3/2}|dq|$. These estimates hold uniformly for $q=k_t$, with $|dq|\le |dk|$.

We write $P=\|d\phi\|_2$, $H=\|\phi\|_2$ and $N=\|S_t\phi\|_2$. The antisymmetric integration by parts and the estimate on the zeroth-order term give
\be
 \lambda^{-1}P^2
 \le 2N^2+C(\lambda^2 a+\lambda^3 a^2)P^2
       +C\lambda^2 BHP+C\lambda^3 B^2H^2 .
\ee
For $a<c\lambda^{-3}$, the second term on the right is bounded by $C(c+c^2)\lambda^{-1}P^2$. We choose $c$ so that $C(c+c^2)\le1/4$. Young's inequality absorbs the mixed term with another $\lambda^{-1}P^2/4$. The resulting graph estimate, the method of continuity and the localization argument above apply with this fixed $a$ and finite $B$.
\end{proof}

A single sufficiently small multiplier coefficient suffices for local $H^1$ regularity of the spinor. The global domain argument below applies this estimate after the local conformal and coordinate changes in Definition~\ref{def:H1-tame}.

\begin{lemma}[Coordinate invariance of the derivative multiplier estimate]
\label{lem:full-derivative-form-intrinsic}
For a metric $g\in H^1_{\mathrm{loc}}$ with uniform spectral bounds relative to $g_{\text{Eucl}}$, condition
\eqref{eq:full-metric-derivative-form} is invariant under smooth
changes of coordinates on compact subcharts. Equivalently, it can
be stated with $|\nabla^h g|_h$, reference volume and reference
energy for any smooth background $h$.
\end{lemma}
\begin{proof}
On a compact overlap, the smooth coordinate change and its inverse
have bounded derivatives and bounded Jacobians. The tensor
transformation rule gives
\be
 |d g_{\rm new}|\le C|d g_{\rm old}|+C',
\ee
where $C'$ is finite because $g$ and the second coordinate
derivatives are bounded. We change variables in the multiplier inequality for scalar functions. The two gradient energies and two volume densities are
uniformly equivalent there. Given a desired derivative coefficient,
we choose the old $\varepsilon$ sufficiently small; the bounded
term contributes only to the zeroth-order constant. Reversing the
coordinate change proves equivalence. Replacing coordinate
derivatives by $\nabla^h$ likewise adds a bounded multiple of $g$,
and changing a smooth reference changes only these bounded terms
and equivalent norms. This proves the intrinsic formulation.
\end{proof}

Every $W^{1,3}_{\mathrm{loc}}$ metric with uniform spectral bounds relative to $g_{\text{Eucl}}$ satisfies
\eqref{eq:full-metric-derivative-form}: we split $|dg|$ at height $M$,
use Sobolev on its $L^3$ tail, and put the bounded part in the
zeroth-order constant. Thus the local regularity result also applies when $g\in W^{1,3}_{\mathrm{loc}}$. At $H^1_{\mathrm{loc}}$ regularity, we use the multiplier estimate as a separate condition. Proposition~\ref{thm:global-critical-Dirac-domain} uses this multiplier estimate to identify the global Dirac domain.

\paragraph{A Lorentz-space sufficient condition.}\label{sec:lorentz-form-bound}
For the potential $|dg|^2$, the inequality~\eqref{eq:full-metric-derivative-form} allows the coefficient of the Dirichlet energy to be arbitrarily small, with a finite zeroth-order term; see Maz'ya--Verbitsky \cite{MazyaVerbitsky2005} for estimates of this kind for general Schr\"odinger potentials. It is implied by a standard Lorentz-space
hypothesis on all components of $dg$. On a relatively compact chart
$U$, suppose that $|dg|$ belongs to the closure of $L^\infty(U)$ in
$L^{3,\infty}(U)$. Then for every $\varepsilon>0$ there is a finite
$C_{\varepsilon,U}$ such that
\be
 \int_U|dg|^2|\varphi|^2dx
 \le\varepsilon\int_U|d\varphi|^2dx
                 +C_{\varepsilon,U}\int_U|\varphi|^2dx.
\ee
\begin{proof}[Proof of the Lorentz-space criterion]
We write $|dg|=f_0+f_1$, with $f_0$ bounded and $\|f_1\|_{3,\infty}$
as small as prescribed. Lorentz--H\"older and the Sobolev embedding
$H^1_0\hookrightarrow L^{6,2}$ give
\be
 \|f_1\varphi\|_2\le C\|f_1\|_{3,\infty}\|\varphi\|_{6,2}
                      \le C\|f_1\|_{3,\infty}\|d\varphi\|_2.
\ee
 The bounded part is
controlled by $\|f_0\|_\infty\|\varphi\|_2$. Squaring and choosing the tail norm small enough gives the multiplier inequality with any prescribed positive coefficient of the Dirichlet energy. On a bounded domain
$L^{3,\infty}\subset L^2$, so this sufficient condition remains within
$H^1_{\mathrm{loc}}$; it does not impose a preferred metric direction.
\end{proof}


\subsection{The scalar measure and the Lichnerowicz identity}\label{sec:weighted-domain}
\paragraph{The integrated identity.}
Following the spinorial method of Witten \cite{Witten1981} and
the presentation given in Lee--LeFloch \cite{LeeLeFloch}, we now derive the local identity that relates spinorial energy to the scalar distribution. Its weighted form also applies when the coframe connection is only locally integrable.
We combine the terms linear in the connection before integrating by parts.

\begin{proposition}[The weighted Lichnerowicz identity]
Let $\theta$ be an oriented coframe with locally bounded coefficients and inverse. Assume $d\theta\in L^1_{\mathrm{loc}}$ and that the contracted quadratic term $T$ of \eqref{eq:coframe-density-test} belongs to $L^1_{\mathrm{loc}}$.
Let $w\in H^1_{\mathrm{loc}}\cap L^\infty_{\mathrm{loc}}$ be real
with $w\omega\in L^2_{\mathrm{loc}}$. For $\psi=w\tau$, with
$\tau$ compact smooth, put $v=w^2|\tau|^2$. Then
\be
\label{eq:weighted-grouped-new}
 \int(|D\psi|^2-|\nabla\psi|^2)dv_g
 =-\frac14\int\rho F\cdot dv-\frac14\int\rho Tv\,dx.
\ee
In the current integral, expand $dv=2w|\tau|^2dw+w^2d|\tau|^2$. Each term is integrable by $wF\in L^2_{\mathrm{loc}}$, the bound on $w$, and compact support of $\tau$.
\end{proposition}
\begin{proof}
The two energies are finite: $d(w\tau)=\tau\,dw+w\,d\tau$ and
$\Omega(w\tau)=(w\Omega)\tau$ both belong to $L^2$. We give the grouping used to integrate their
difference. We choose Hermitian Pauli matrices $\gamma^a$ with
$\gamma^a\gamma^b+\gamma^b\gamma^a=2\delta^{ab}I$ and set
$\gamma^{ab}=[\gamma^a,\gamma^b]/2$. With $c^a=i\gamma^a$ and the
stated torsion convention, the spin connection matrix is
$\Omega_i=\frac14\omega_{ab i}\gamma^a\gamma^b$. We put
$\Gamma^{ij}=e_a^ie_b^j\gamma^{ab}$. Pointwise Clifford multiplication
gives the following three Clifford identities
\be
 \Gamma^{ij}\Omega_j=-\tfrac14 F^iI+\tfrac14B_{ab}^i\gamma^{ab},
 \quad
 -\Omega_i\Gamma^{ij}=-\tfrac14 F^jI-\tfrac14B_{ab}^j\gamma^{ab},
 \quad
 -\Omega_i\Gamma^{ij}\Omega_j=-\tfrac14 TI.
\ee
For these calculations, we evaluate at an orthonormal coframe at the
point. We write $\omega_{ab i}=\varepsilon_{abd}K_{di}$; then
\be
 F^i=2\varepsilon_{ijd}K_{dj},\qquad
 B_{ab}^i\varepsilon_{abe}=-2(K_{ie}-\delta_{ie}\operatorname{tr}K),
 \qquad T=\operatorname{tr}(K^2)-(\operatorname{tr}K)^2.
\ee
Substitution in $\Omega_i=(i/2)K_{di}\gamma^d$, followed by
$\gamma^a\gamma^b=\delta^{ab}I+i\varepsilon_{abc}\gamma^c$,
gives the displayed contractions; we then restore tensor indices.
Thus for coefficient fields $\eta,\xi$, the polarized difference of the two energies has coordinate density
\be
\begin{aligned}
 \mathcal L(\eta,\xi)={}&\frac\rho2 C_{ab}^{ij}
 \langle\del_i\eta,\gamma^{ab}\del_j\xi\rangle
 -\frac14\rho F^i\del_i\langle\eta,\xi\rangle
 -\frac14\rho T\langle\eta,\xi\rangle\\
 &+\frac\rho4 B_{ab}^i\bigl(
 \langle\del_i\eta,\gamma^{ab}\xi\rangle
 -\langle\eta,\gamma^{ab}\del_i\xi\rangle\bigr).
\end{aligned}
\ee
This is an algebraic expansion using the Clifford anticommutator.
Its quadratic scalar coefficient is
$((\operatorname{tr}K)^2-\operatorname{tr}(K^2))/4=-T/4$;
the quadratic bivector terms cancel. Four-form coefficients vanish
in dimension three. Thus the expression groups all rough quadratic
terms into the specified $T$.
We now integrate this expansion for $\eta=\xi=w\tau$, taking
real parts. Antisymmetry in $i,j$ cancels the terms containing two
derivatives of $w$. Skew Hermiticity of $\gamma^{ab}$ combines
the other derivative terms into
\be
 \frac12\int\rho C_{ab}^{ij}
       \del_i\bigl(w^2\operatorname{Re}
                  \langle\tau,\gamma^{ab}\del_j\tau\rangle\bigr)dx.
\ee
Indeed the derivative of $w^2$ accounts for both mixed terms;
the extra second derivative of $\tau$ vanishes after contraction
with $C_{ab}^{ij}$. The divergence identity for the exterior product coefficients therefore makes
this integral equal to
\be
 \frac12\int\rho B_{ab}^j w^2\operatorname{Re}
                  \langle\tau,\gamma^{ab}\del_j\tau\rangle dx.
\ee
This integration by parts is valid at the stated regularity. Its
compactly supported scalar factor is bounded and belongs to $H^1$.
We approximate that factor in $H^1$ and almost everywhere by uniformly
bounded smooth factors. The differentiated term converges against
the bounded coefficients $\rho C_{ab}^{ij}$ by $L^2$ convergence on a compact set; the
undifferentiated term converges against its $L^1$ divergence by
dominated convergence.

The real part of the last line of the algebraic expansion is
\be
 -\frac12\int\rho B_{ab}^j w^2\operatorname{Re}
                  \langle\tau,\gamma^{ab}\del_j\tau\rangle dx.
\ee
The derivatives of $w$ cancel in this expression before integration.
It is integrable because $B\in L^1$ and the other factors are bounded.
Thus these two integrals cancel. The only remaining terms are the
scalar current and $T$ terms. They are integrable since
$wF,dw\in L^2$, $w$ and $\tau$ are bounded, and $T\in L^1$.
This proves the weighted Lichnerowicz identity~\eqref{eq:weighted-grouped-new}.
\end{proof}


\paragraph{The scalar measure as a trace on spinors.}
For $H^1_{\mathrm{loc}}$ metrics, we proved in Section~\ref{sec:weak-coframes} above that a nonnegative scalar curvature measure gives zero mass to sets of zero $H^1$ capacity. To integrate the square of an arbitrary $H^1_0$ spinor, we now prove a trace estimate. For this trace estimate, the coefficient of the Dirichlet norm in the multiplier inequality need only be finite.

\begin{proposition}[The closed scalar identity from a multiplier estimate on the metric derivative]
\label{prop:critical-closed-scalar-identity}
Let $g,g^{-1}\in L^\infty_{\mathrm{loc}}$ and $g\in H^1_{\mathrm{loc}}$. Assume that on every relatively compact ball $V$ there are finite constants $a_V,b_V$ such that
\be
\label{eq:finite-derivative-multiplier}
 \||dg|v\|_{L^2(V)}\le a_V\|dv\|_{L^2(V)}+b_V\|v\|_{L^2(V)}
 \qquad(v\in C_c^\infty(V)).
\ee
In its canonical frame assume that
$\mu=\diver(\rho F)-\rho Tdx$ is nonnegative. On every relatively
compact ball $V$, this measure is diffuse with respect to $H^1$
capacity and has a bounded trace map $H^1_0(V)\to L^2(\mu)$.
For every compactly supported $H^1$ coefficient spinor $\psi$,
where the two spinor norms use $dv_g$,
\be
 \|D\psi\|_2^2=\|\nabla\psi\|_2^2
                  +\tfrac14\int|\widetilde\psi|^2d\mu.
\ee
The trace is its quasi-continuous representative.
\end{proposition}
\begin{proof}
The bounded Sobolev chain rule gives $\omega,F\in L^2_{\mathrm{loc}}$
and $T\in L^1_{\mathrm{loc}}$. For compact smooth $\psi$, we apply the
weighted Lichnerowicz identity with a scalar weight equal to one near
its support. Its right side is precisely
$\tfrac14\langle\mu,|\psi|^2\rangle$, by the definition of $\mu$.
All connection products are integrable at this regularity. Positivity
therefore gives
\be
 \int|\psi|^2d\mu\le4\|D\psi\|_2^2
                         \le C_V\|\psi\|_{H^1(V)}^2.
\ee
For the last inequality use bounded principal coefficients and
$|\Omega|\le C|dg|$, together with the multiplier estimate for the metric derivative.
Thus multiplication by $\Omega$ is bounded from $H^1_0(V)$ to $L^2$.
Taking $\psi=\varphi e$ for a constant unit coefficient spinor $e$
proves the scalar trace bound on smooth functions.

The preceding capacity argument shows that $\mu$ ignores every
capacity-zero Borel set in $V$ independently of the trace estimate.
The trace bound extends by density to $H^1_0(V)$. A strongly convergent
$H^1$ sequence has a subsequence converging quasi-everywhere;
its $L^2(\mu)$ trace limit is therefore the quasi-continuous
representative, since $\mu$ ignores the exceptional set.

Finally, we approximate the given $\psi$ by compact smooth spinors in
$H^1_0(V)$. Their $D$ images and covariant derivatives converge in
$L^2$, because the multiplier estimate for the metric derivative controls multiplication
by $\Omega$ on their $H^1$ differences. Their scalar traces also converge by the established bound, so passage to the limit in the smooth identity proves the assertion.
Polarization gives the corresponding bilinear identity. When the local Dirac regularity estimate holds, this also applies to compact localizations of $L^2$ spinors satisfying the Dirac equation in the sense of distributions.
\end{proof}

Under the hypotheses of the preceding proposition, suppose in addition that solutions of the adjoint equation in the sense of distributions have the local Sobolev regularity stated above, for example by Proposition~\ref{thm:H1-multiplier-Dirac-entry} or its fixed-small-coefficient extension. Let $D\psi=f$ hold in the sense of distributions on a chart, with $\psi,f\in L^2_{\mathrm{loc}}$. For every real smooth
compactly supported cutoff $\eta$, the scalar trace of $\psi$ is defined
locally and
\be
 \int\eta^2|\nabla \psi|^2dv_g
       +\frac14\int\eta^2|\widetilde \psi|^2d\mu
 \le 4\int\eta^2|f|^2dv_g
                  +6\int|d\eta|_g^2|\psi|^2dv_g.
\ee
In particular, a spinor satisfying $D\psi=0$ in the sense of distributions has a locally square-integrable covariant derivative and a trace that is locally square integrable with respect to the scalar curvature measure.
\begin{proof}[Proof of the local energy and trace estimate]
The assumed local regularity first gives $\psi\in H^1_{\mathrm{loc}}$. Thus the
closed identity applies to $\eta \psi$ and all its scalar traces are
defined by the preceding proposition. We expand
$D(\eta \psi)=\eta f+c(d\eta)\psi$ and
$\nabla(\eta \psi)=\eta\nabla \psi+d\eta\otimes \psi$.
The terms $|d\eta|_g^2|\psi|^2$ on the two sides cancel exactly.
Writing $E=\int\eta^2|\nabla \psi|^2dv_g$,
$M=\frac14\int\eta^2|\widetilde \psi|^2d\mu$,
$F=\int\eta^2|f|^2dv_g$ and
$G=\int|d\eta|_g^2|\psi|^2dv_g$, Cauchy--Schwarz gives
\be
 E+M\le F+2\sqrt{FG}+2\sqrt{EG}
       \le 2F+3G+\tfrac12E.
\ee
Hence $E+2M\le4F+6G$, which implies the stated estimate since
$M\ge0$. After local $H^1$ regularity has been established, every product in these identities is integrable on the support of $\eta$.
\end{proof}

\subsection{Dirac domains and local conformal representatives}
\paragraph{The self-adjoint domain from local multiplier estimates.}
Bandara and Saratchandran~\cite{BandaraSaratchandran} studied essential self-adjointness for first-order operators with rough coefficients. Here the local regularity and scalar trace estimates justify the Lichnerowicz identity on the physical Dirac domain.

\begin{proposition}[The global Dirac domain from local multiplier estimates]
\label{thm:global-critical-Dirac-domain}
Let $M$ be an oriented smooth three-manifold without boundary, with a compact interior and one Euclidean exterior chart. One chooses a complete smooth reference metric $h$
equal to the Euclidean metric on that exterior. Let $g$ belong to $H^1_{\mathrm{loc}}$ and have uniform spectral bounds relative to $h$ on $M$. At every point of $M$, including points on the end, suppose there is a coordinate ball $B$ on which the metric derivative satisfies the multiplier estimate below: for some finite $b_B$ and some $a_B$ smaller than the threshold determined by the spectral bounds for $g$ in the local Dirac estimate,
\be
\label{eq:local-small-multiplier}
 \||dg|v\|_{L^2(B)}\le a_B\|dv\|_{L^2(B)}+b_B\|v\|_{L^2(B)}
 \qquad(v\in C_c^\infty(B)).
\ee
Then the following conclusions hold.
\begin{enumerate}
\item The physical Dirac operator on compactly supported smooth reference spinors is essentially self-adjoint. The compactly supported elements of its closed domain are exactly the compactly supported $H^1$ spinors.
\item If $R_g^{\mathrm w}$ is nonnegative, one writes $\mu$ for the associated Radon measure. Every such compactly supported spinor satisfies the closed Lichnerowicz identity with scalar term $\tfrac14\int |\widetilde\psi|^2\,d\mu$.
\end{enumerate}
\end{proposition}
\begin{proof}
{\bf 1. \sl Spin frames and local domains.} An oriented three-manifold admits a spin structure; we fix one for
$h$. We identify its orthonormal frames with those for $g$ by the
positive symmetric metric comparison. In a smooth $h$ coframe
$\eta$, the resulting $g$ coframe is $A\eta$, with
$A$ bounded, uniformly positive and in $H^1_{\mathrm{loc}}$.
Under a smooth reference rotation $S$ this comparison obeys
$A'=SAS^{-1}$ and $\theta'=S\theta$. Thus its spin transitions
are those of the fixed smooth reference bundle. Compact reference
$H^1$ spaces are consequently well defined and the actual local
connection coefficients belong to $L^2$ and are bounded
$H^1$-to-$L^2$ multipliers locally by \eqref{eq:local-small-multiplier} and a smooth partition of unity. A smooth exterior satisfies this condition with $a_B=0$ on each relatively compact ball. On an exterior where the metric is only $H^1_{\mathrm{loc}}$, the support of each compact localization is covered by finitely many of the selected balls. Their zeroth-order constants may grow along the end.

The local Dirac regularity estimate with a fixed small multiplier coefficient applies on the selected balls. To compare them with the coordinate canonical frame on
a contractible chart, we write $\eta=H(x)dx$ and let $G$ be a
positive definite matrix. The rotation relating the two coframes is
\bse
\be
 A(x,G)=(H^{-t}GH^{-1})^{1/2},\qquad
 R(x,G)=A(x,G)HG^{-1/2}\in SO(3).
\ee
The covering $\operatorname{Spin}(3)\to SO(3)$ gives a smooth lift of $R$ on the simply connected product of the chart with the positive definite matrix cone. On a compact subchart, its first derivatives are bounded when the eigenvalues of $G$ range over a fixed interval. Evaluating the lift at $G=g(x)$ and using the Sobolev chain rule yields a transition $U$ with $|dU|\le C(|dg|+1)$. Multiplication by $U$ and its inverse preserves $H^1$:
the derivative product is controlled by the finite multiplier estimate \eqref{eq:local-small-multiplier}.
For a bounded coefficient field $\phi$ in $H^1_0$, the product rule
and the multiplier inequality~\eqref{eq:local-small-multiplier} give
\be
 \|d(U\phi)\|_2\le C\|d\phi\|_2+C\|(|dg|+1)\phi\|_2
                         \le C_U\|\phi\|_{H^1},
\ee
\ese
and density removes the boundedness restriction on $\phi$, with the same conclusion for $U^{-1}$. The actual connection and Dirac transformation
laws identify the local operators as maps $H^1_0\to L^2$.

This also justifies the change of frame in the adjoint equation.
For a solution of the adjoint equation expressed in a reference spin frame, the adjoint identity defining the equation in the sense of distributions initially uses spinors with compactly supported smooth coefficients. Since the physical Dirac is bounded from
$H^1_0$ to $L^2$ on the patch, this identity extends by density
to all compact $H^1$ tests. A smooth test in the coordinate
canonical frame transforms into such an admissible test by the
multiplier bound for $U$ and $U^{-1}$. Hence the transformed field
satisfies the coordinate-frame adjoint equation, without assuming
its $H^1$ regularity. Proposition~\ref{thm:H1-multiplier-Dirac-entry} and its fixed-small-coefficient extension then apply.
The invariance under changes of Sobolev coframe, the agreement with the coordinate formula in Section~\ref{sec:weak-coframes} above, and the closed Lichnerowicz identity show that all local identities use the same nonnegative scalar curvature measure $\mu$.
Under this sign assumption the scalar traces are compatible under changes of spin frame. Indeed, $U\in H^1_{\mathrm{loc}}\cap L^\infty$ is unitary almost everywhere, so its quasi-continuous representative is unitary outside a set of $H^1$ capacity zero. For $\phi\in H^1_{\mathrm{loc}}$, we choose smooth $\phi_j\to \phi$ in $H^1$ on a compact subpatch. The multiplier estimate gives $U\phi_j\to U\phi$ in $H^1$. After extraction, both sequences converge to their quasi-continuous representatives outside a common capacity-zero set. Since $\widetilde{U\phi_j}=\widetilde U\,\phi_j$ there, passing to the limit gives $\widetilde{U\phi}=\widetilde U\,\widetilde \phi$ quasi-everywhere. This also covers unbounded $\phi$. Since $\mu$ does not charge capacity-zero sets, the norm of the trace is unchanged by $U$ almost everywhere with respect to $\mu$.

\vskip.2cm \noindent
{\bf 2. \sl The adjoint equation.} We denote by $D_c$ the operator on compact smooth reference
coefficients. The integrated identity for the frame current makes it symmetric.
When the scalar density is nonnegative, a finite smooth quadratic partition of unity sums the local scalar identities on any compact support. Explicitly, if $\sum_i\chi_i^2=1$ near the support of $\psi$, then
\be
 \sum_i\|D(\chi_i\psi)\|_2^2=\|D\psi\|_2^2+\int\sum_i|d\chi_i|_g^2|\psi|^2\,dv_g,
\ee
and the same identity holds with $D$ replaced by $\nabla$. The cross terms vanish because $\sum_i\chi_i d\chi_i=0$, and the cutoff squares agree by the Clifford relation. The trace terms sum to $\int|\widetilde\psi|^2d\mu$ by the frame compatibility and $\sum_i\chi_i^2=1$. The local Dirac regularity estimate shows that compact localizations of every solution in the adjoint domain belong to $H^1$, and identifies their images under $D$. Conversely compact $H^1$ fields are approximated
by smooth reference coefficients in graph norm, since the
connection multiplication is continuous from $H^1$ to $L^2$
on their support. When $\mu\ge0$, the same approximation and the local scalar trace bounds extend the positive identity to them.

\vskip.2cm \noindent
{\bf 3. \sl Approximation on the complete manifold.} Finally, we choose smooth exhaustion cutoffs $\zeta_R$ equal to one
on the compact interior, supported in the interior together with
$\{r<2R\}$, with $|d\zeta_R|_g\le C/R$ on the exterior.
If $\psi\in\operatorname{Dom}D_c^*$ and $f=D_c^*\psi$, then
\be
 D(\zeta_R\psi)=\zeta_R f+c(d\zeta_R)\psi.
\ee
The local regularity estimate places the compactly supported spinor $\zeta_R\psi$ in the closure of the smooth compact domain for the Dirac graph norm. As $R\to+\infty$, it converges to $\psi$ in $L^2$ and
its displayed image converges to $f$ in $L^2$. Therefore
$\operatorname{Dom}D_c^*\subset\operatorname{Dom}\overline{D_c}$.
The reverse inclusion and equality of the operators follow from
symmetry. This proves self-adjointness and the domain assertions. This completes the proof of Proposition~\ref{thm:global-critical-Dirac-domain}.
\end{proof}


\paragraph{Transport through conformal representatives and coordinates.}
We now transfer the local Dirac estimate to the conformal representatives used in condition~\ref{item:tame-local} of Definition~\ref{def:H1-tame}. The coordinate maps and conformal factors can vary from one neighborhood to the next. This step identifies the physical operator domain without imposing a multiplier estimate on the derivative of the conformal factor.

\begin{proposition}[The Dirac domain in local Sobolev coordinates]\label{prop:local-conformal-charts}
Let $M$ be a connected oriented smooth three-manifold without boundary. Let $h$ be a complete smooth metric that is Euclidean on each of $N\ge1$ exterior ends, with compact complement, and let $g\in H^1_{\mathrm{loc}}(M)$ satisfy uniform spectral bounds relative to $h$. Suppose that at every point $g$ admits the local factorization and small multiplier estimate in condition~\ref{item:tame-local} of Definition~\ref{def:H1-tame}.

Let $\mathcal C$ be the compactly supported smooth sections of the reference spinor bundle. The closure of $D_g$ on $\mathcal C$ is self-adjoint, and every compact localization of an element of its maximal domain admits approximation by elements of $\mathcal C$ in the Dirac graph norm. If $R_g^{\mathrm w}\ge0$, its $L^2$ kernel is zero.
\end{proposition}
\begin{proof}
{\bf 1. \sl Transport of bounded tests.} We fix a spin structure and identify the reference spinor bundle with the physical bundle by the positive symmetric metric comparison. In a selected neighborhood we write
\be
                       g=w^4X^*k .
\ee
We shrink the neighborhood and orient the image coordinates consistently. The inverse of $X$ belongs locally to $W^{2,2}$: its derivative is $(dX)^{-1}\circ X^{-1}$, and differentiation of the matrix inverse gives products of bounded matrices with the weak Hessian of $X$.

We use a smooth reference spin frame on the source and the canonical $k$ frame on the image. The rotation between the source physical coframe and $w^2X^*\theta_k$ belongs to $H^1_{\mathrm{loc}}$, since both coframe matrices and their inverses are bounded and locally $H^1$. On a smaller simply connected neighborhood it has an $H^1$ lift to $\operatorname{Spin}(3)$ by the Sobolev lifting theorem recalled by Van Schaftingen~\cite{VanSchaftingenLifting2025}. We write $P$ for pullback in the resulting spin frames, including this unitary change of coefficients. Both $P$ and $P^{-1}$ preserve bounded compactly supported $H^1$ fields. Change of variables and the bounds on $w,w^{-1}$ give the local $L^2$ bounds needed below.

For bounded compactly supported $H^1$ spinors $\phi$ on the image,
\be
\label{eq:local-rough-conformal-transport}
 D_g(w^{-2}P\phi)=w^{-4}P(D_k\phi),\qquad
 dv_g=w^6X^*dv_k .
\ee
The first identity is the conformal law for the spin Dirac operator in dimension three~\cite{HitchinHarmonic}, together with covariance under $X$. It remains valid at the stated regularity. Indeed $d(X^*\theta_k)=X^*(d\theta_k)$ by the weak chain rule and the symmetry of the weak Hessian of $X$. Uniqueness in the torsion equation identifies the pulled-back connection. The spin lift and the conformal factor obey the bounded $H^1$ chain and product rules. Every differentiated coefficient is locally $L^2$ and multiplies a bounded field.

A bounded compactly supported $H^1$ spinor on either side admits uniformly bounded smooth approximations in $H^1$ converging almost everywhere. Its connection products converge in $L^2$ by dominated convergence against the squared connection coefficients, so these approximations also converge in the Dirac graph norm. Consequently such fields belong to the minimal domain and may be used as tests in an adjoint equation with field and source in $L^2$.

\vskip.2cm \noindent
{\bf 2. \sl Local graph approximation.} We consider $\psi$ in the maximal domain of $D_g$ and put $f=D_g\psi\in L^2(dv_g)$. In the selected chart we set $z=P^{-1}(w^2\psi)$ and test its adjoint equation with $w^{-2}P\phi$, where $\phi$ is compact and smooth in the image. The preceding approximation makes this test admissible. The volume factor $w^6$ and the two factors $w^{-4}$ and $w^{-2}$ in the pairings give
\be
 (z,D_k\phi)_{L^2(dv_k)}
       =(P^{-1}(w^4f),\phi)_{L^2(dv_k)} .
\ee
Both fields are locally $L^2$. The local Dirac estimate for $k$ therefore gives $z\in H^1_{\mathrm{loc}}$ in the image frame.

For a smooth source cutoff $\eta$ compactly supported in this neighborhood, $z_\eta=(\eta\circ X^{-1})z$ is a compactly supported $H^1$ image spinor. We approximate it in $H^1$ by compact smooth $\phi_j$. The derivative multiplier bound for $k$ gives $D_k\phi_j\to D_kz_\eta$ in $L^2$. Each $w^{-2}P\phi_j$ is a bounded compactly supported $H^1$ source spinor, hence belongs to the minimal domain by the preceding step. The fields converge in physical $L^2$ to $\eta\psi$, while~\eqref{eq:local-rough-conformal-transport} gives convergence of their Dirac images to $w^{-4}P(D_kz_\eta)$. Closedness places $\eta\psi$ in the minimal domain; its image is $\eta f+c_g(d\eta)\psi$ by the maximal equation. This argument uses $H^1$ regularity in the image frame and bounded test fields in the source frame.

A smooth partition of unity gives graph approximation for every compact localization. We choose exhaustion cutoffs $\chi_R$ that are equal to one on the compact core and satisfy $|d\chi_R|_g\le C/R$ on every end. Then
\be
 \chi_R\psi\longrightarrow\psi,\qquad
 D_g(\chi_R\psi)=\chi_Rf+c_g(d\chi_R)\psi\longrightarrow f
       \quad\hbox{in }L^2 .
\ee
Thus the maximal domain is contained in the minimal domain. Symmetry gives the reverse inclusion and self-adjointness.

\vskip.2cm \noindent
{\bf 3. \sl The $L^2$ kernel.} We assume $R_g^{\mathrm w}\ge0$ and denote its measure by $\mu$. The integrated Lichnerowicz identity on $\mathcal C$, Kato's inequality, and the scalar Sobolev inequality proved in Lemma~\ref{prop:global-scalar-Sobolev} give
\be
 \|\phi\|_6\le C\|\nabla\phi\|_2
              \le C\|D_g\phi\|_2\qquad(\phi\in\mathcal C).
\ee
Here the scalar inequality depends only on the reference geometry and the spectral bounds; its proof does not use the Dirac domain. If $n\in L^2$ satisfies $D_gn=0$ by adjoint testing, the domain equality provides $\phi_j\in\mathcal C$ with $\phi_j\to n$ in $L^2$ and $D_g\phi_j\to0$ in $L^2$. The displayed inequality gives $\phi_j\to0$ in $L^6$, hence in $L^2$ on every relatively compact set. Therefore $n=0$.
\end{proof}

\section{Bounded Witten spinors and the positive mass identity}\label{sec:witten-end-estimates}\label{sec:mass-identity}

\subsection{Annular estimates and the global Sobolev inequality}
Using the Dirac domain and compact identity established in Section~\ref{sec:isotropic}, we construct bounded solutions with prescribed constant spinors at infinity, called \underline{Witten spinors}. The terminal condition is expressed in $L^6$. The following annular estimates and global Sobolev inequality yield the mass identity in Proposition~\ref{thm:positive-domain-mass} below, including on the two-ended orientation cover.

\begin{lemma}[Three mixed annular errors]
Suppose the coefficient fields $q,\psi$ satisfy $q,\psi\in L^6(\Mext)$ and $dq,d\psi\in L^2(\Mext)$. Then
\be
 \frac1R\int_{A_R}|q||d\psi|\to0,\qquad
 \frac1R\int_{A_R}|dq||\psi|\to0,\qquad
 \frac1R\int_{A_R}|d\psi||\psi|\to0.
\ee
Under the same hypotheses,
$R^{-1}\int_{A_R}|q||d\psi||\psi|\to0$.
\end{lemma}
\begin{proof}
H\"older's inequality and $|A_R|^{1/3}\le CR$ give, for example,
\be
 R^{-1}\int_{A_R}|q||d\psi|
 \le C\|q\|_{6,A_R}\|d\psi\|_{2,A_R}\to0.
\ee
The next two estimates follow by exchanging the factors. Finally
$\|\psi\|_{3,A_R}\le CR^{1/2}\|\psi\|_{6,A_R}$ gives
\be
 R^{-1}\int_{A_R}|q||d\psi||\psi|
 \le C R^{-1/2}\|q\|_{6,A_R}\|d\psi\|_{2,A_R}\|\psi\|_{6,A_R}\to0.
\ee
\end{proof}

\begin{lemma}[Sobolev inequality with an arbitrary compact core]\label{prop:global-scalar-Sobolev}
Let $(M,h)$ be a connected smooth three-manifold with $N\ge1$
Euclidean exterior regions and a compact smooth core. If $g$ has uniform spectral bounds relative to $h$, then
\be
                  \|f\|_{L^6(g)}\le C\|df\|_{L^2(g)}
                  \qquad(f\in C_c^\infty(M)).
\ee
\end{lemma}
\begin{proof}
The spectral bounds for $g$ relative to $h$ reduce the statement to $h$. We fix one exterior
end and an annulus $A=\{R_0<r<2R_0\}$ in it. Since $f$ is compactly supported,
Cauchy--Schwarz along each ray gives
\be
 |f(r\omega)|^2\le r^{-1}
            \int_r^\infty|\del_sf(s\omega)|^2s^2ds.
\ee
Integration on $A$ yields
$\|f\|_{L^2(A)}\le C_{R_0}\|df\|_{L^2(\Mext)}$.
We choose a connected relatively compact domain $K'$ containing the
core and one inner annulus in each end, including $A$. The compact Poincar\'e inequality with an anchored
annular average gives
\be
 \|f\|_{L^2(K')}\le C\bigl(\|df\|_{L^2(K')}
                                      +\|f\|_{L^2(A)}\bigr).
\ee
For completeness, a failure would yield functions normalized in
$L^2(K')$ with gradients and annular norm tending to zero. Rellich
compactness would give a nonzero constant limit on the connected
$K'$, while its annular norm is zero, a contradiction.
Compact Sobolev and the anchored estimate control the $L^6$ norm
on the core by the global gradient norm. On each exterior apply
Euclidean homogeneous Sobolev to a cutoff times $f$ that vanishes
near its inner boundary. The cutoff derivatives are controlled by the
compact estimate just proved. Summing over the finitely many ends
and combining with the core estimate proves the inequality. The spectral bounds for $g$ relative to $h$ then give the estimate with the norms defined by $g$.
\end{proof}


\subsection{Constant terminal spinors and the compact energy identity}
We work on an oriented manifold with a fixed spin structure. Its restriction to each Euclidean exterior is the standard spin structure, so constant spinor coefficients are available in the canonical exterior frame. We write $\mathcal C$ for the compactly supported smooth sections of the reference spin bundle, identified with physical spinors as in Section~\ref{sec:isotropic}. Their connection coefficients are locally square-integrable, and their pointwise Hermitian norm is preserved by this identification.

\begin{lemma}\label{lem:constant-tail-spinor-energy}
Let $M$ be connected, with a compact core and $N\ge1$ Euclidean ends $E_a$. Suppose that $g\in H^1_{\mathrm{loc}}$ satisfies uniform spectral bounds relative to the complete reference metric, has nonnegative scalar density $\mu$, and satisfies the two integrals in~\eqref{eq:end-energy} on each end. Assume $\mu(E_a)<+\infty$, and let $m_a$ be the corresponding dyadic annular ADM mass. One chooses a smooth reference section $s$ equal to a constant spinor $\xi_a$ sufficiently far out on $E_a$. For every $\phi\in\mathcal C$, the field $t=s+\phi$ has finite energy
\be
 \mathcal E(t)\coloneqq\int_M|\nabla t|^2\,dv_g
                +\frac14\int_M|t|^2\,d\mu,
\ee
and satisfies
\be
\label{eq:constant-tail-spinor-energy}
 \mathcal E(s+\phi)
   =4\pi\sum_a m_a|\xi_a|^2+\|D(s+\phi)\|_2^2.
\ee
\end{lemma}
\begin{proof}
Section~\ref{sec:mass-existence} gives finite masses. Local finiteness of $\mu$ on the compact core gives $\mu(M)<+\infty$. The spinor $t$ is bounded, and on each constant tail its covariant derivative and Dirac image are bounded by $C|dg||\xi_a|$. These observations give the stated integrability.

The compact Lichnerowicz identity of Section~\ref{sec:weighted-domain} applies in the reference-induced physical frames. A smooth quadratic partition of unity gives the same identity globally: the terms involving derivatives of the partition cancel after polarization, and all scalar terms use the measure $\mu$ from Section~\ref{sec:weak-coframes}. We choose $\chi_R$ equal to one on the core, with the same smooth radial cutoff on each end. Polarizing with $\eta t$ and $\chi_Rt$, where $\eta=1$ near $\operatorname{supp}\chi_R$, gives
\be
\label{eq:weighted-spinor-current-focused}
 \int_M\chi_R\bigl(|\nabla t|^2-|Dt|^2\bigr)\,dv_g
       +\frac14\int_M\chi_R|t|^2\,d\mu
   =-\sum_a\int_{E_a}J_t\cdot d\chi_R,
\ee
where, in exterior coordinates, $c^i=e_b^ic_b$ and
\be
 J_t^i\coloneqq\rho\operatorname{Re}
        \langle t,(\nabla^i+c^iD)t\rangle,
 \qquad \rho=\sqrt{\det(g_{ij})}.
\ee

We record the algebra that identifies the constant-tail current. We set
$B_g^{ij}=g^{ij}\operatorname{Id}_{\Sigma}+c^ic^j$. In an oriented orthonormal coframe, the contracted spin connection is
\be
 c^b\Omega_b=\sigma_g\operatorname{Id}_{\Sigma}-\frac14c(F^\flat),
 \quad F^b=2\sum_c\omega_{bc c},
 \quad \sigma_g=-\frac14\sum_{b,c,d}\varepsilon_{bcd}\omega_{bc d}\in\R.
\ee
This follows by multiplying the Clifford matrices in the Pauli representation. Skew Hermiticity of $\Omega^i$ and $c^i$, together with the Clifford anticommutator, then gives
\be
\label{eq:spinor-current-algebra-focused}
 J_t^i=\rho\operatorname{Re}\langle t,B_g^{ij}\del_jt\rangle
                       +\frac14\rho F^i|t|^2.
\ee
For completeness, linearizing the canonical coframe $g^{1/2}$ at the Euclidean metric gives the torsion connection
\be
 \omega_{bc i}=\frac12(\del_c h_{bi}-\del_b h_{ci})
 \quad\hbox{to first order},\qquad h=g-g_{\text{Eucl}}.
\ee
Its contraction is $P_i(g)=\del_jh_{ij}-\del_i h_{jj}$. The bounded Sobolev matrix chain rule and smooth dependence within the spectral bounds therefore give
\be
\label{eq:frame-current-ADM-remainder-focused}
 |\rho F-P(g)|\le C|h||dg|,\qquad
 |F|\le C|dg|,\qquad
 |\rho B_g-B_0|\le C|h|,
\ee
where $B_0^{ij}=\delta^{ij}\operatorname{Id}_{\Sigma}+c_ic_j$.
On the support of $d\chi_R$ the coefficients of $t$ equal $\xi_a$, so its derivative term in~\eqref{eq:spinor-current-algebra-focused} vanishes. The error in replacing $\rho F$ by $P(g)$ is bounded by
\be
 C|\xi_a|^2\|h/r\|_{2,A_R}\|dg\|_{2,A_R}\longrightarrow0.
\ee
The smooth annular mass and dyadic annular mass have the same limit by Lemma~\ref{lem:dyadic-annular-mass}. Thus the right side of~\eqref{eq:weighted-spinor-current-focused} tends to $4\pi\sum_a m_a|\xi_a|^2$. Dominated convergence on the left proves~\eqref{eq:constant-tail-spinor-energy}.
\end{proof}


\subsection{The bounded Witten construction}
The domain conclusion of Proposition~\ref{prop:local-conformal-charts} allows us to solve the Dirac equation by completing $\mathcal C$ in its Dirac norm. We first construct this completion; boundedness will then give local $H^1$ regularity of the resulting affine solution in the physical frame.

\begin{proposition}[The homogeneous Witten problem]\label{thm:homogeneous-witten}
Let $M$ be a connected oriented smooth three-manifold with a compact core and $N\ge1$ Euclidean ends. Let $g\in H^1_{\mathrm{loc}}$ satisfy uniform spectral bounds relative to a complete smooth reference metric that is Euclidean on these ends, and let its scalar density be a nonnegative measure $\mu$. Suppose that the compact Lichnerowicz identity holds on $\mathcal C$, that the minimal and maximal $L^2$ domains of $D$ coincide, and that the only $L^2$ spinor satisfying $D\varphi=0$ in the sense of distributions is zero. These properties hold under the local factorizations of Proposition~\ref{prop:local-conformal-charts}. If $V$ is the completion of $\mathcal C$ for $\|D\cdot\|_2$, then
\be
\label{eq:focused-homogeneous-Dirac-isomorphism}
 D:V\longrightarrow L^2(dv_g)
\ee
is an isometric bijection, and $V$ embeds injectively into $L^6$ with $\|\psi\|_6\le C\|D\psi\|_2$.
\end{proposition}
\begin{proof}
Kato's inequality, the preceding global scalar Sobolev inequality, and the positive compact identity give
\be
\label{eq:focused-Dirac-core-coercivity}
 \|\phi\|_6\le C\|\nabla\phi\|_2\le C\|D\phi\|_2,
 \qquad
 \|D\phi\|_2^2=\|\nabla\phi\|_2^2
                   +\frac14\|\phi\|_{L^2(\mu)}^2
 \quad(\phi\in\mathcal C).
\ee
Thus a Dirac-Cauchy sequence has an $L^6$ limit $\psi$, a covariant derivative limit in $L^2$, and a scalar trace limit in $L^2(\mu)$. The latter two limits are, at this stage, the maps defined by completion. If $D\phi_j\to f$ in $L^2$, symmetry against $\tau\in\mathcal C$ gives
\be
 (f,\tau)=\lim_j(\phi_j,D\tau)=(\psi,D\tau),
\ee
because $D\tau$ is compactly supported and belongs to $L^2$, hence to $L^{6/5}$. If $\psi=0$, density of $\mathcal C$ in $L^2$ implies $f=0$. This proves injectivity into $L^6$ and identifies $D\psi=f$ in the sense of distributions.

By construction $D$ is an isometry onto the closure of $D\mathcal C$. If $\varphi\in L^2$ is orthogonal to that range, then $(\varphi,D\tau)=0$ for every $\tau\in\mathcal C$, so $D\varphi=0$ in the sense of distributions. The stated domain result and absence of such nonzero fields imply $\varphi=0$. The closed range is therefore all of $L^2$, proving the proposition.
\end{proof}

\begin{proposition}[A bounded Witten spinor]\label{prop:bounded-Witten-spinor}
Under the hypotheses of Proposition~\ref{thm:homogeneous-witten}, assume the end integrals and finiteness of $\mu(E_a)$ in Lemma~\ref{lem:constant-tail-spinor-energy}. For any constant terminal values $\xi_a$, let $s$ be a smooth reference extension equal to $\xi_a$ on each sufficiently far end. There is a unique field
\be
 \Psi=s+\psi,\qquad \psi\in V,\qquad D\Psi=0.
\ee
Writing $c=\max_a|\xi_a|$, it satisfies
\be
\label{eq:bounded-Witten-regularity}
 \Psi\in H^1_{\mathrm{loc}},\qquad
 |\widetilde\Psi|\le c\ \hbox{quasi-everywhere},\qquad
 \nabla\Psi\in L^2(dv_g),\qquad
 \widetilde\Psi\in L^2(\mu),
\ee
and its energy is
\be
\label{eq:bounded-Witten-mass-energy}
 \int_M|\nabla\Psi|^2\,dv_g
        +\frac14\int_M|\widetilde\Psi|^2\,d\mu
      =4\pi\sum_a m_a|\xi_a|^2.
\ee
On each constant tail, the coefficient derivative $d(\Psi-\xi_a)$ belongs to $L^2$. Here the tilde denotes the actual quasi-continuous representative of the locally $H^1$ spinor.
\end{proposition}

For one end, the energy in~\eqref{eq:bounded-Witten-mass-energy} is $4\pi m_{\rm ADM}|\xi|^2$. On the orientation cover, prescribing $\xi$ on one end and zero on the other selects the mass term of that end.

\begin{proof}
The connection estimate on the ends gives $Ds\in L^2$. The preceding proposition gives a unique $\psi\in V$ with $D\psi=-Ds$. We choose $\phi_j\in\mathcal C$ with $D(s+\phi_j)\to0$ in $L^2$. If all terminal values are zero, uniqueness gives $\Psi=0$, so we assume $c>0$.

\emph{Bounded approximations.} For $\delta>0$, we choose a smooth radial contraction $T_\delta$ of the spinor coefficient space that is the identity on $\{|z|\le c\}$, has differential of norm at most one, and satisfies $|T_\delta(z)|\le\min(|z|,c+\delta)$. It is obtained by smoothing the radial clipping map outside that ball and commutes with unitary changes of frame. Differentiating this equivariance along the skew-Hermitian connection matrices gives
\be
 \nabla T_\delta(t)=dT_\delta(t)\nabla t
\ee
almost everywhere for smooth reference sections $t$. Since $\mu\ge0$, the contraction decreases both terms of $\mathcal E$. The sections
\be
 t_j\coloneqq T_{1/j}(s+\phi_j)
\ee
still belong to $s+\mathcal C$, because they keep every terminal value fixed. We put $M_\xi=4\pi\sum_a m_a|\xi_a|^2$. The constant-tail identity gives
\be
 0\le\|Dt_j\|_2^2=\mathcal E(t_j)-M_\xi
       \le\mathcal E(s+\phi_j)-M_\xi
       =\|D(s+\phi_j)\|_2^2\longrightarrow0.
\ee
Consequently $t_j-s$ converges to $\psi$ in $V$ and in $L^6$. The compact identity applied to $t_j-t_k\in\mathcal C$ makes $\nabla t_j$ and the scalar traces of $t_j$ Cauchy in their global $L^2$ spaces.

\emph{Identification of the physical regularity and trace.} After extraction, $t_j\to\Psi$ almost everywhere for volume. In a relatively compact reference spin frame the connection matrix $\Omega$ belongs to $L^2$. The uniform bound $|t_j|\le c+1/j$ and dominated convergence therefore give
\be
 \Omega t_j\longrightarrow\Omega\Psi\quad\hbox{in }L^2_{\mathrm{loc}}.
\ee
Since $dt_j=\nabla t_j-\Omega t_j$, the coefficient derivatives converge strongly in $L^2$ on every compact subpatch. Their limit is $d\Psi$, so $\Psi\in H^1_{\mathrm{loc}}$, and the covariant derivative and Dirac limits are the actual ones. On a compact exhaustion, a further subsequence converges quasi-everywhere. The measure $\mu$ ignores sets of zero local $H^1$ capacity by Lemma~\ref{lem:scalar-measure-capacity}; hence the $L^2(\mu)$ limit is $\widetilde\Psi$, with $|\widetilde\Psi|\le c$ quasi-everywhere. Passing to the limit in $\mathcal E(t_j)=M_\xi+\|Dt_j\|_2^2$ proves~\eqref{eq:bounded-Witten-mass-energy}.

On a tail where $s=\xi_a$, the field $\psi=\Psi-\xi_a$ is bounded and belongs to $L^6$. The constant coefficient section has $\nabla\xi_a=\Omega\xi_a$, with $\Omega\in L^2$ on the end. Thus $\nabla\psi\in L^2$ there and
\be
 d\psi=\nabla\psi-\Omega\psi\in L^2,
\ee
which completes the proof.
\end{proof}


\subsection{The mass identity and nonnegativity}\label{sec:positive-domain}
\paragraph{The terminal spinorial integral.}
The bounded construction also permits the localized energy identity to be applied directly to the Witten field. We give the boundary estimate needed for that argument.

\begin{lemma}[Integral identification of the spinorial terminal flux]
\label{prop:integral-spinorial-terminal}
On an end with the spectral bounds and integral decay in~\eqref{eq:end-energy}, let $\Psi=\xi+\psi$ be bounded and locally $H^1$, where $\xi$ is constant, $\psi\in L^6$ and $d\psi\in L^2$. For the current $J_\Psi$ in~\eqref{eq:spinor-current-algebra-focused}, one has
\be
 -\int J_\Psi\cdot d\chi_R
       +\frac{|\xi|^2}{4}\int P(g)\cdot d\chi_R\longrightarrow0.
\ee
If the dyadic annular ADM average has a finite limit $m_{\rm ADM}$, then
\be
 -\int J_\Psi\cdot d\chi_R\longrightarrow4\pi m_{\rm ADM}|\xi|^2.
\ee
In particular, this applies to the spinor of Proposition~\ref{prop:bounded-Witten-spinor} on each end.
\end{lemma}
\begin{proof}
The algebra in~\eqref{eq:spinor-current-algebra-focused} and the bounds in~\eqref{eq:frame-current-ADM-remainder-focused} apply almost everywhere. The flat linear derivative term cancels exactly: the Clifford relation gives $B_0^{ij}=-B_0^{ji}$, so integration by parts yields
\be
 \int\del_i\chi_R\,
       \operatorname{Re}\langle\xi,B_0^{ij}\del_j\psi\rangle\,dx=0.
\ee
Indeed the symmetric Hessian of $\chi_R$ contracts to zero with $B_0$; a cutoff extension across the inner region justifies the same calculation on the end.

After subtracting this term and the constant-spinor ADM term, the pointwise remainder is bounded by
\be
 C\bigl(|h||d\psi|+|\psi||d\psi|
       +|dg||\psi|+|dg||\psi|^2+|h||dg|\bigr),
 \qquad h=g-g_{\text{Eucl}}.
\ee
The first three terms have vanishing annular integrals after division by $R$, by the $L^6$--$L^2$ earlier estimates. Boundedness of $\psi$ reduces the fourth term to the third. For the last term, Cauchy--Schwarz gives
\be
 R^{-1}\int_{A_R}|h||dg|\,dx
    \le2\|h/r\|_{2,A_R}\|dg\|_{2,A_R}\longrightarrow0.
\ee
Since $|d\chi_R|\le C/R$, these estimates prove the displayed current comparison.

To identify the finite limit, we write
$F(t)\coloneqq(16\pi)^{-1}\int_{S_t}P(g)\cdot n\,dS$ for almost every $t>R_0$.
Polar integration gives $F\in L^1_{\mathrm{loc}}((R_0,+\infty))$ and
$m_{\text{dya}}(R)=R^{-1}\int_R^{2R}F(t)\,dt$.
We fix $a_0>R_0$ and put $G(t)\coloneqq\int_{a_0}^t F(s)\,ds$.
Thus $G$ is locally absolutely continuous and
$m_{\text{dya}}(R)=[G(2R)-G(R)]/R$ is continuous.
If $m_{\text{dya}}(R)\to m\in\R$, this function is bounded for $R\ge a_0$.
For $t\ge a_0$, we choose $n$ with $a=t2^{-n}\in[a_0,2a_0)$ and obtain
\be
 \frac{G(t)}t
 =\frac{G(a)}t+\sum_{j=1}^n2^{-j}m_{\text{dya}}(t2^{-j})
 \longrightarrow m .
\ee
Indeed $G$ is bounded on $[a_0,2a_0]$, and the boundedness of
$m_{\text{dya}}$ permits dominated convergence in the geometrically weighted sum.

For the radial profile of $\chi$, we put $w(s)\coloneqq-\chi'(s)$.
Then $w(1)=w(2)=0$ and $\int_1^2w(s)\,ds=1$.
Integration by parts gives
\be
 m_\chi(R)=\int_1^2F(Rs)w(s)\,ds
   =-\int_1^2s\,\frac{G(Rs)}{Rs}\,w'(s)\,ds
   \longrightarrow-m\int_1^2s w'(s)\,ds=m .
\ee
The convergence of $G(t)/t$ is uniform on every tail, so the last passage
is valid for this fixed smooth profile.
The preceding current comparison now gives the asserted finite terminal limit.
\end{proof}


\paragraph{The positive mass identity.}
The operator domain and the bounded Witten spinors constructed above give the mass identity directly. We first prove the identity for finitely many ends, so that the same argument applies to an orientation cover, and then state its one-ended consequence.
\begin{proposition}[The positive mass identity]\label{thm:positive-domain-mass}
Assume the hypotheses of Proposition~\ref{prop:bounded-Witten-spinor}, with a finite ADM mass $m_a$ on each Euclidean end and prescribed constant terminal spinors $\xi_a$. The bounded Witten spinor $\Psi$ satisfies
\be
\label{eq:integral-Witten-mass}
 4\pi\sum_a m_a|\xi_a|^2
   =\int_M|\nabla\Psi|^2\,dv_g
      +\frac14\int_M|\widetilde\Psi|^2\,d\mu
   \ge0 .
\ee
If the left side vanishes, $\Psi$ is parallel in the sense of distributions and $|\widetilde\Psi|^2\mu=0$.
\end{proposition}
\begin{proof}
Proposition~\ref{prop:bounded-Witten-spinor} gives a bounded spinor $\Psi\in H^1_{\mathrm{loc}}$ with $D_g\Psi=0$, the prescribed terminal values, $\nabla\Psi\in L^2$, and the scalar trace $\widetilde\Psi\in L^2(\mu)$. On each end $\Psi-\xi_a\in L^6$ and $d(\Psi-\xi_a)\in L^2$ by the same proposition. Thus Lemma~\ref{prop:integral-spinorial-terminal} applies on each end.

For a nonnegative compact smooth real cutoff $\chi$, we choose a compact smooth $\eta$ equal to one near $\operatorname{supp}\chi$. We apply the polarized integrated Lichnerowicz identity to $\eta\Psi$ and $\chi\Psi$. These spinors are bounded and locally $H^1$, so the bounded approximation and scalar-capacity arguments of Sections~\ref{sec:weak-coframes} and~\ref{sec:weighted-domain} justify the pairing with their quasi-continuous traces. Because $D_g\Psi=0$, the Dirac pairing is zero: the two remaining terms $c_g(d\eta)\Psi$ and $c_g(d\chi)\Psi$ have disjoint supports. The terms containing $d\eta$ in the covariant pairing also vanish. Taking real parts gives
\be
 \int_M\chi|\nabla\Psi|^2\,dv_g
  +\frac14\int_M\chi|\widetilde\Psi|^2\,d\mu
 =-\int_M\operatorname{Re}
       \langle\nabla\Psi,\Psi\otimes d\chi\rangle_g\,dv_g .
\ee
The right side is $-\int J_\Psi\cdot d\chi$, since the Dirac term in $J_\Psi$ vanishes.

We choose $\chi=\chi_R$ equal to one on the compact core, with the same nonincreasing smooth radial profile on every end, so that $\chi_R\uparrow1$. The terminal-flux identification gives the limit $4\pi\sum_a m_a|\xi_a|^2$ on the right. Monotone convergence on both nonnegative terms on the left proves~\eqref{eq:integral-Witten-mass}. If the left side is zero, both integrals vanish, which gives the asserted parallel equation and scalar-measure identity.
\end{proof}

\begin{corollary}[Nonnegativity and parallel spinors]\label{thm:reference-one-ended-positivity}
Let $g$ be an $H^1$-tame metric with reference $h$ on an oriented manifold $M$, and suppose that $R_g^{\mathrm w}\ge0$. One denotes the associated scalar curvature measure by $\mu$. Its ADM mass lies in $[0,+\infty]$. When it is finite, each terminal spinor $\xi$ determines the bounded Witten spinor of Proposition~\ref{prop:bounded-Witten-spinor}, and
\be
 4\pi m_{\rm ADM}|\xi|^2
  =\int_M|\nabla\Psi_\xi|^2\,dv_g
       +\frac14\int_M|\widetilde\Psi_\xi|^2\,d\mu .
\ee
If $m_{\rm ADM}=0$ and $\xi\ne0$, then $\Psi_\xi$ is parallel with constant norm $|\xi|$ and $\mu=0$.
\end{corollary}
\begin{proof}
The end conditions in Definition~\ref{def:H1-tame}, together with the assumed nonnegative scalar density, give an ADM limit in $\R\cup\{+\infty\}$, by Proposition~\ref{prop:annular-negative-variation}. If the scalar measure has infinite mass on the end, that limit is $+\infty$. Otherwise it is finite. Proposition~\ref{prop:local-conformal-charts} gives the self-adjoint physical domain and the zero $L^2$ kernel from the local factorizations in condition~\ref{item:tame-local} of the definition. The compact positive identity and the bounded Witten construction therefore apply, and~\eqref{eq:integral-Witten-mass} gives the mass identity and nonnegativity.

At zero mass the spinor is parallel. Metric compatibility makes its norm constant on the connected manifold. Since $\Psi_\xi-\xi\in L^6$ on an end of infinite volume, this constant is $|\xi|$. The scalar term in the identity is therefore $|\xi|^2\mu(M)/4$, so $\mu=0$ when $\xi\ne0$. Section~\ref{sec:rigidity} derives the Euclidean reconstruction from this parallel spinor.
\end{proof}

\section{Coordinate invariance and the Dini regularity criterion}\label{sec:coordinate-invariance}

\subsection{Changes of Euclidean reference}
The comparison of Euclidean references uses the local Sobolev regularity, the spectral bounds and the two quadratic integrals on the end. We then construct the harmonic coordinates needed for the Dini criterion in Proposition~\ref{prop:H1-tame-criteria}.


Following Bartnik's coordinate-invariance question \cite{Bartnik1986}, we compare the annular expressions in two asymptotic charts. The two quadratic integrability conditions determine an admissible class of transition maps. The comparison of the masses at scale $R$ does not require a sign condition on the scalar density and applies to a mass taking values in the extended real line once the existence criterion of Section~\ref{sec:mass-existence} holds.

\begin{proposition}[Integral coordinate invariance]
\label{thm:absolute-integral-coordinate-invariance}
Let $g\in H^1_{\mathrm{loc}}$ on an end, with uniform spectral bounds relative to $g_{\text{Eucl}}$ in the smooth coordinates $x$, and suppose
\be
 \int_{\Mext}\bigl(|dg|^2+|x|^{-2}|g-g_{\text{Eucl}}|^2\bigr)dx<+\infty.
\ee
Let $x=X(y)$ be a $W^{1,1}_{\mathrm{loc}}$ homeomorphism between exterior regions that each contain the complement of a compact set in $\R^3$, and suppose that $X$ takes infinity to infinity. One sets $K\coloneqq DX$ and chooses a Borel representative of $g_x$ satisfying these uniform spectral bounds at every point and set $g_y=K^t(g_x\circ X)K$. Under either condition below this definition is independent almost everywhere of the representative. Then the
following two conditions are equivalent:
\begin{enumerate}
\item $g_y\in H^1_{\mathrm{loc}}$ and satisfies uniform spectral bounds relative to $g_{\text{Eucl}}$ in the $y$ coordinates, and
$\int(|dg_y|^2+|y|^{-2}|g_y-g_{\text{Eucl}}|^2)\,dy<+\infty$;
\item $K$ and $K^{-1}$ are bounded, $dK\in L^2(dy)$, and
$K-Q\in L^6(dy)$ for a constant $Q\in O(3)$.
\end{enumerate}
These conditions are imposed on sufficiently far exterior regions, and $K^{-1}$ denotes the pointwise matrix inverse. On these exterior regions, either condition implies that $X$ and $X^{-1}$ are locally Lipschitz and their weak second derivatives are square integrable on the end. Every fixed smooth radial cutoff function $\chi$ as in \eqref{eq:annular-mass-definition} satisfies
\be
\label{eq:matched-radial-flux}
 m_{\chi,y}(R)-m_{\chi,x}(R)\longrightarrow0.
\ee
Thus a limit for this cutoff function in one chart is also its limit in the other. If the distributional scalar density is a signed Radon measure $\mu$ with finite negative part near infinity, the limit in both charts exists in $\R\cup\{+\infty\}$ and is independent of the cutoff function. It is finite exactly when the positive part of $\mu$ is finite near infinity. For a smooth end these measure conditions are equivalent to integrability of $(R_g)_-$ and $(R_g)_+$, respectively. The comparison uses only the assumptions on the asymptotic end.
\end{proposition}
\begin{proof}
{\bf 1. \sl Bounds for the transition and its inverse.} Under (1), the uniform spectral bounds for $g_x$ and $g_y$ in the identity $g_y=K^t(g_x\circ X)K$ bound $K$ and $K^{-1}$. Under (2) these bounds are given. Since $K$ is the weak derivative, $X$ is locally Lipschitz, and its inverse $Y$ has the same local regularity. To verify this before differentiating the metric, we write $J=|\det K|$, which has a positive lower bound. At differentiability points with invertible derivative, the sign of the determinant agrees with the local topological degree, by comparison with the derivative on a sufficiently small sphere. A homeomorphism has constant orientation on the connected exterior, so this sign is a constant $\sigma$. For each $a$, the vector field with components $J(K^{-1})_{ja}$ has zero divergence in the sense of distributions: up to the constant orientation sign, these components are the coefficients of the closed forms $dX^b\wedge dX^c$. This identity follows by local smooth approximation of the Lipschitz coordinate functions in $W^{1,2}$. For $\phi\in C_c^\infty$ in the image, the area formula and this identity give
\be
 \int Y^i\partial_a\phi\,dx
 =\int y^iJ(K^{-1})_{ja}
               \partial_j(\phi\circ X)\,dy
 =-\int J(K^{-1})_{ia}(\phi\circ X)\,dy .
\ee
Integration by parts allows compact Lipschitz tests since these coefficients are locally bounded. Changing variables back gives $DY=K^{-1}\circ Y$ in the sense of distributions, hence the asserted local Lipschitz property. Both maps now preserve null sets.

\vskip.2cm \noindent
{\bf 2. \sl Second derivatives of the transition.} We assume (1). To obtain the second derivatives of $X$, we put $A_z=\sqrt{\det g_z}\,g_z^{-1}$ for $z=x,y$. These coefficients are bounded, uniformly positive and locally $H^1$. Changing variables in the divergence of $A_x^{ba}$ in the sense of distributions gives, for each component $U=X^a$,
\bse
\be
 \operatorname{div}_y(A_y\,dU)=F,\qquad
 F=J(\partial_b A_x^{ba})\circ X\in L^2_{\mathrm{loc}}.
\ee
Indeed $A_y=JK^{-1}(A_x\circ X)K^{-t}$, and testing with $\phi$ in $y$ amounts to testing the $x$ equation with $\phi\circ Y$. We give the local estimate because the coefficients have only one square-integrable derivative. We choose a directional difference quotient $\delta_s$ and a smooth cutoff $\eta$ on a ball compactly contained in the parameter domain. We apply the difference quotient to the equation and test the resulting equation with $\eta^2\delta_sU$. On a fixed slightly larger ball $B$, the coefficient error is bounded by $\|\delta_s A_y\|_{2,B}\|dU\|_{\infty,B}\|d(\eta^2\delta_sU)\|_2$. Moving the difference quotient off $F$ bounds its contribution by $\|F\|_{2,B}\|d(\eta^2\delta_sU)\|_2$. Ellipticity and Young's inequality therefore give
\be
 \|\eta\,d\delta_s U\|_2
 \le C_\eta\big(\|F\|_{2,B}+\|dA_y\|_{2,B}\|dU\|_{\infty,B}
                           +\|dU\|_{2,B}\big)
\ee
for all sufficiently small $|s|$, enlarging $B$ within the chart for translations. Weak compactness in each coordinate direction gives $X\in W^{2,2}_{\mathrm{loc}}$. We may now differentiate the metric identity. Permuting its three indices and using $\partial_iK^a_j=\partial_jK^a_i$ gives
\be
 \partial_i\partial_jX^a
 =K^a_k\Gamma(g_y)^k_{ij}
 -(\Gamma(g_x)^a_{bc}\circ X)K^b_iK^c_j.
\ee
\ese
All products here are locally square integrable by the Sobolev chain rule for the locally bi-Lipschitz map. The two metric derivative integrals bound the Christoffel symbols in
$L^2$; bounded Jacobians therefore give $dK\in L^2(dy)$. Each
quadratic end condition also gives $g-g_{\text{Eucl}}\in L^6$ in its chart.
We extend $K$ boundedly across a compact collar, with square-integrable
derivative. Sobolev--Poincar\'e on expanding balls and a convergent
subsequence of its bounded averages give a constant $K_\infty$ with
$K-K_\infty\in L^6$. The metric transformation and the two $L^6$
normalizations imply $K^tK-g_{\text{Eucl}}\in L^6$, while boundedness of $K$
gives $K^tK-K_\infty^tK_\infty\in L^6$. Hence the constant
$K_\infty^tK_\infty-g_{\text{Eucl}}$ vanishes, so $Q=K_\infty\in O(3)$.
This proves condition (2). For $Y=X^{-1}$, differentiating
$DY=K^{-1}\circ Y$ also shows $D^2Y\in L^2(dx)$.

\vskip.2cm \noindent
{\bf 3. \sl The converse implication.} Conversely, we assume condition (2). Since $K$ is bounded,
$\det K-\det Q\in L^6$; the fixed sign of $\det K$ and the bound
on $K^{-1}$ imply $\operatorname{sgn}(\det K)=\det Q$. We set
$Z_Q=X-Qy$ and multiply $Z_Q$ by a smooth cutoff vanishing near the
inner boundary and equal to one far out, then extend it by zero.
Its derivative belongs to $L^6(\R^3)$, so the homogeneous Morrey
inequality gives $|Z_Q(y)-Z_Q(0)|\le C|y|^{1/2}$. In particular
$|X(y)|/|y|\to1$; this follows from the integral transition condition. We put $h_x=g_x-g_{\text{Eucl}}$.
Uniform Jacobian bounds and the resulting radius comparability
preserve finiteness of its $L^6$ and weighted $L^2$ norms under composition. The transformation law
$g_y=K^t(g_x\circ X)K$ gives
\be
 |g_y-g_{\text{Eucl}}|\le C(|h_x\circ X|+|K-Q|),\qquad
 |dg_y|\le C(|dg_x\circ X|+|dK|).
\ee
The bounded field $K-Q$ has square-integrable derivative and belongs
to $L^6$. Cutting it off on annuli gives convergence in the
homogeneous gradient norm, since
$R^{-1}\|K-Q\|_{L^2(A_R)}\le C\|K-Q\|_{L^6(A_R)}\to0$.
Hardy's inequality therefore gives its weighted $L^2$ norm.
This proves finiteness of the two quadratic integrals in the $y$ coordinates and also gives $g_y-g_{\text{Eucl}}\in L^6$. In both directions $DY=K^{-1}\circ Y$ gives $D^2Y\in L^2$ by the matrix inverse chain rule and change of variables.

\vskip.2cm \noindent
{\bf 4. \sl The scalar measure and the annular mass.} The scalar measure is also compatible with these coordinates. We pull back the $x$ coframe by $X$. Its coefficients are bounded and locally $H^1$ by the regularity just proved. In its exterior derivative the terms containing second derivatives of $X$ cancel by symmetry. Thus its Cartan connection is $X^*\omega_x$ by torsion uniqueness. The contraction $T$ pulls back as a scalar, and $\rho F$ transforms as a vector density, with factor $JK^{-1}$. Testing the scalar formula in $y$ therefore tests the $x$ formula against $\phi\circ Y$. This compact Lipschitz test admits smooth approximations converging uniformly and in $H^1$, which suffice for both the signed measure and the divergence formula. The pulled-back and canonical coframes in $y$ differ by a bounded $H^1_{\mathrm{loc}}$ orthogonal transformation. A constant reflection accounts for the orientation when needed, and the Sobolev rotation identity in Section~\ref{sec:weak-coframes} gives the same scalar density. In particular its positive and negative parts transform as measures.

For the comparison of the annular masses, we apply Lemma~\ref{lem:integral-end-mass-comparison} below. The preceding steps give its transition assumptions, and both charts satisfy the required weighted metric integral.

If the scalar density is a signed measure with finite negative part near infinity, Proposition~\ref{prop:annular-negative-variation} gives existence, independence of the cutoff function and the finiteness criterion. The comparison just proved identifies the limits in the two charts, including the value $+\infty$. This completes the proof of Proposition~\ref{thm:absolute-integral-coordinate-invariance}.
\end{proof}

\begin{lemma}\label{lem:integral-end-mass-comparison}
Let $x=X(y)$ be a locally bi-Lipschitz $W^{2,2}_{\mathrm{loc}}$ change of exterior coordinates taking infinity to infinity. Suppose that $K=DX$ and $K^{-1}$ are uniformly bounded and, for some $Q\in O(3)$,
\be
 dK\in L^2,\qquad K-Q\in L^6.
\ee
Suppose that $g_y=K^t(g_x\circ X)K$, that both metrics satisfy uniform spectral bounds relative to $g_{\text{Eucl}}$, and that
\be
 g_z\in H^1_{\mathrm{loc}},\qquad
 |z|^{-1}(g_z-g_{\text{Eucl}})\in L^2,\qquad z=x,y.
\ee
Then every fixed smooth radial cutoff function satisfies
$m_{\chi,y}(R)-m_{\chi,x}(R)\to0$.
\end{lemma}
\begin{proof}
An orthogonal change of the $x$ coordinates reduces to $Q=I$. We put $H=K-I$ and $Z=X-y$. On every annulus $A_R^*=\{cR<|y|<CR\}$, with fixed $0<c<C$, H\"older's inequality gives
\be
 \|H\|_{L^2(A_R^*)}=o(R),\qquad
 \|Z\|_{L^2(A_R^*)}=o(R^2).
\ee
For the second estimate, Morrey's inequality on a finite chain of overlapping balls covering a rescaled annulus bounds its oscillation by
$CR^{1/2}\|DZ\|_{L^6(\{c'R<|y|<C'R\})}=o(R^{1/2})$.
The averages of $Z$ on overlapping dyadic annuli differ by the same bound. Summing from one fixed annulus gives an average of size $o(R^{1/2})$, since a sum of increments $\varepsilon_j2^{j/2}$ with $\varepsilon_j\to0$ is $o(2^{j/2})$. This proves the asserted $L^2$ estimate and also gives $|X(y)|/|y|\to1$ uniformly. Thus the supports of the radial tests in either chart lie in a fixed $A_R^*$.

We write $h_z=g_z-g_{\text{Eucl}}$ and $h=h_x\circ X$. The weighted integrals and the bounded Jacobians give
$\|h\|_{L^2(A_R^*)}+\|h_y\|_{L^2(A_R^*)}=o(R)$.
For $T_R=D^2\chi_R-g_{\text{Eucl}}\Delta\chi_R$, the definition of the smooth annular expression is
$16\pi m_{\chi,z}(R)=\int h_z:T_R\,dz$.
Expanding the pullback gives
\be
 h_y=h+H+H^t+E_R,\qquad
 |E_R|\le C\bigl(|H|^2+|h||H|\bigr).
\ee
Its error paired against $T_R$ tends to zero, because
$\|T_R\|_\infty\le CR^{-2}$.
Each row of $T_R$ has zero divergence and $H=DZ$, so
$\int(H+H^t):T_R=0$ by integration by parts; one may first cut $Z$ off outside a neighborhood of the support of $T_R$.
Finally, we put $J=|\det K|$. Since $|J-1|\le C|H|$ and
$\|DT_R\|_\infty\le CR^{-3}$, the difference between
$\int h:T_R(y)\,dy$ and $\int h:J T_R(X)\,dy$ is bounded by
\be
 C\|h\|_{L^2(A_R^*)}
 \left(R^{-2}\|H\|_{L^2(A_R^*)}
             +R^{-3}\|Z\|_{L^2(A_R^*)}\right)=o(1).
\ee
The latter integral is the annular expression in $x$ by change of variables.
\end{proof}


\subsection{Derivative control from nonnegative scalar curvature in harmonic coordinates}\label{sec:harmonic-scalar-control}
In harmonic coordinates, the nonnegativity of the scalar density yields the bound on the metric derivatives used in the local Dirac estimate. We first prove this implication for coefficients close to a constant matrix, then construct harmonic charts from a condition on mean oscillation in the original atlas.

\begin{proposition}[Scalar measure and metric derivative in a harmonic chart]\label{prop:harmonic-scalar-control}
Let $k\in H^1_{\mathrm{loc}}(B)$ on a ball in $\R^3$, and suppose its scalar density $\mu_k$ is nonnegative. One sets $A=\sqrt{\det k}\, k^{-1}$. Suppose $\partial_i A^{ij}=0$ in the sense of distributions and $\|k-g_{\text{Eucl}}\|_{L^\infty(B)}\le\delta$. For universal constants $\delta_0>0$ and $C$, if $0\le\delta\le\delta_0$, then
\be
\label{eq:harmonic-scalar-control}
 \int_B z^2\,d\mu_k+\int_B z^2|dk|^2\,dx
      \le C\delta\int_B|dz|^2\,dx
       \qquad(z\in C_c^\infty(B;\R)).
\ee
In particular, arbitrarily large derivatives are allowed in this estimate. A constant positive definite comparison matrix can be reduced to $g_{\text{Eucl}}$ by a linear change of harmonic coordinates.
\end{proposition}
\begin{proof}
We write $h=k-g_{\text{Eucl}}$, $\rho=\sqrt{\det k}$ and $\ell=\log\rho$. In the coordinate scalar identity of Section~\ref{sec:weak-coframes}, metric compatibility gives
\be
 \mu_k=-\partial_i(A^{ij}\partial_j\ell)+Q,\qquad
 Q=A^{ij}\Gamma^a_{ib}\Gamma^b_{ja}.
\ee
Indeed the divergence coefficient in general coordinates is $A^{ij}\Gamma^a_{ij}-A^{aj}\partial_j\ell$, and the integrable remainder is $A^{ij}(\Gamma^a_{ib}\Gamma^b_{ja}-\Gamma^a_{ij}\partial_a\ell)$. These expressions follow by inserting $\partial_i A^{jk}=A^{jk}\partial_i\ell-\Gamma^j_{ib}A^{bk}-\Gamma^k_{ib}A^{jb}$ in the integrated coordinate formula. Here $A^{ij}\Gamma^a_{ij}=-\partial_i A^{ia}=0$. In particular $|Q|\le C|dh|^2$.

In the quadratic expression for $Q$, replacing the undifferentiated metric coefficients by their Euclidean values costs at most $C\delta|dh|^2$. In the resulting contraction, the two factors, before division by four, are $\partial_i h_{ab}\pm(\partial_bh_{ia}-\partial_a h_{ib})$. Summation gives
\be
 Q=-\tfrac14|dh|^2+\tfrac12\partial_bh_{ia}\,\partial_a h_{ib}
                     +\Ocal(\delta|dh|^2).
\ee
We put $X=\|z\,dh\|_2$, $Y=\|dz\|_2$, $L=\|z\,d\ell\|_2$ and $M_z=\int z^2\,d\mu_k$. The integral of the crossed product against $z^2$ differs from $\int z^2|\operatorname{div}h|^2$ by at most $C\delta XY$. Indeed, we pair $\partial_a\partial_b h_{ib}=\partial_b\partial_a h_{ib}$ in $H^{-1}$ with $z^2h_{ia}$ and sum. The two errors contain $2zh_{ia}\,dz$ and one derivative of $h$. Also $\operatorname{div} h=d\ell+\Ocal(\delta|dh|)$, by differentiating $A=\rho k^{-1}$ in the harmonic condition. Testing the scalar identity with $z^2$ now gives
\be
 M_z+\tfrac14 X^2\le C\big(LY+L^2+\delta X^2+\delta XY\big).
\ee

To control the volume derivative, we use $\ell z^2$ in the same scalar identity. This bounded compact $H^1$ test is admissible by the scalar measure and capacity result of Section~\ref{sec:weak-coframes}. The uniform ellipticity of $A$ and the bound $\|\ell\|_\infty\le C\delta$ give
\be
 L^2\le C\delta(M_z+X^2)+C\delta LY
      \le C'\delta(M_z+X^2)+\tfrac12 L^2+C'\delta^2Y^2 .
\ee
Here the sign of the scalar measure gives $|\int\ell z^2\,d\mu_k|\le C\delta M_z$. Substitute the resulting bound for $L$ in the preceding inequality. Young's inequality absorbs the term $C\sqrt{\delta(M_z+X^2)}Y$ at a cost $C'\delta Y^2$. Taking $\delta_0$ sufficiently small absorbs the remaining multiple of $\delta(M_z+X^2)$ and proves~\eqref{eq:harmonic-scalar-control}.
\end{proof}

We now construct the harmonic charts for the Dini criterion in Proposition~\ref{prop:H1-tame-criteria}. A pointwise Dini modulus of the metric coefficients implies the local Dini mean oscillation condition.

\begin{proposition}[Harmonic coordinates and the local Dirac estimate]\label{thm:dini-harmonic-mass}
Let $g$ be a measurable Riemannian metric with $H^1_{\mathrm{loc}}$ regularity and locally uniform spectral bounds relative to a smooth reference metric. If its coefficients have local Dini mean oscillation in a smooth chart about a point, there is a bi-Lipschitz $W^{2,2}_{\mathrm{loc}}$ harmonic chart about that point in which the metric coefficients are continuous. Whenever a harmonic chart of this regularity with continuous metric coefficients is available and $R_g^{\mathrm w}\ge0$ in the neighborhood, the chart verifies condition~\ref{item:tame-local} of Definition~\ref{def:H1-tame}, with $w=1$. If $R_g^{\mathrm w}\ge0$ and these charts cover a manifold satisfying conditions~\ref{item:tame-reference}--\ref{item:tame-decay} of Definition~\ref{def:H1-tame}, they establish its remaining local condition, so the metric is $H^1$-tame.
\end{proposition}
\begin{proof}
{\bf 1. \sl Construction of harmonic coordinates.} We work first in a smooth chart in which the coefficients have local Dini mean oscillation. The matrix $A=\sqrt{\det g}\,g^{-1}$ also has local Dini mean oscillation. Indeed the map from $g$ to $A$ is Lipschitz on the set of symmetric matrices satisfying the prescribed uniform spectral bounds, and we can subtract its value at the mean of $g$ on each ball. The metric has a continuous representative: the means on concentric balls of radii $s, s/2$ differ by at most $C\omega_{g,K}(s)$, and the dyadic tail is bounded by $C\int_0^{2s}\omega_{g,K}(t)\,dt/t$ since $\omega_{g,K}(s)\le C\omega_{g,K}(t)$ for $s\le t\le2s$. The continuous ball means converge locally uniformly and agree almost everywhere with the metric.

We center coordinates at a point and on $B_4(0)$ set $A_r(y)=A(ry)$. We solve $\operatorname{div}(A_r dU_r^j)=0$ with boundary values $y^j$ by the Dirichlet principle. We subtract the affine function and test its difference against itself. Since $A(0)$ is constant, ellipticity yields
$\|dU_r-I\|_{2,B_4}\le C\|A_r-A(0)\|_{2,B_4}\to0$.
For a uniform gradient estimate, we use the standard Euclidean basis $(\mathbf e_j)$ and apply Dong and Kim's interior bound to $V_r^j=U_r^j-y^j$, which solves the divergence equation with data $(A(0)-A_r)\mathbf e_j$ \cite[Theorem~1.5 and estimate (2.17)]{DongKimDini}. We put $e_r=\int_0^{r}\omega_{A,K}(t)\,dt/t+\omega_{A,K}(r)$ on a fixed larger compact subball. Then
\be
 \|dV_r\|_{L^\infty(B_1)}\le C\big(\|dV_r\|_{L^1(B_4)}+e_r\big)\longrightarrow0 .
\ee
Indeed the mean oscillations of the rescaled coefficient and of the data at scale $s\le1$ are bounded by $\omega_{A,K}(rs)$. The modified modulus in that estimate sums dilates at scales $\min(\kappa^{-j}s,1)$ with geometric weights $2^{-j}$ for a fixed $0<\kappa<1$. Its Dini integral is at most $C e_r$, by a change of variable on $s<\kappa^j$ and $\sum_{j\ge1} j2^{-j}<+\infty$ on the remaining intervals. Thus the coefficient integral used for absorption in (2.17) is small uniformly for small $r$. Since the same theorem also gives $U_r\in C^1(B_1)$, the inverse function theorem yields, for some $r>0$, a local $C^1$ bi-Lipschitz harmonic coordinate map $X(x)=rU_r(x/r)$.

\vskip.2cm \noindent
{\bf 2. \sl Sobolev regularity of the charts.} These coordinates also belong locally to $W^{2,2}$. For a coordinate difference quotient $\Delta_s$, the equation for $\Delta_s X$ has leading matrix $A(\cdot+s\mathbf e_i)$ and divergence data $-(\Delta_s A)dX$. The data have uniformly bounded local $L^2$ norm since $A\in H^1_{\mathrm{loc}}$ and $dX$ is locally bounded. Testing with $\eta^2\Delta_s X$ yields a uniform $L^2$ bound for $\eta d\Delta_s X$ for each compact cutoff. Weak compactness gives the Hessian, and the formula for $dX^{-1}$ together with the bounded inverse Jacobian gives the same regularity for $X^{-1}$. Consequently $k=(X^{-1})^*g$ is continuous and locally $H^1$. Changing variables in the harmonic coordinate equations in the sense of distributions gives $\partial_i(\sqrt{\det k}\,k^{ij})=0$.

\vskip.2cm \noindent
{\bf 3. \sl The estimate in harmonic coordinates.} We now use a harmonic chart constructed above or given directly. The scalar density on the image is $X_*\mu_g$. To check this at the stated regularity, we use the same physical orthonormal coframe on both sides by pullback. Its exterior derivative commutes with pullback by the Sobolev chain rule and symmetry of the weak Hessian of $X$. The torsion equation identifies the connection one-forms under pullback, as in~\eqref{eq:local-rough-conformal-transport}. In the integrated coframe curvature formula the derivative of the connection is paired by parts with the complementary coframe times the scalar test. These factors are bounded and locally $H^1$ and the connection is locally $L^2$; the remaining quadratic products are locally $L^1$. Change of variables therefore gives the same pairing, using bounded $H^1$ tests in the source chart. Thus $\mu_k\ge0$. A constant linear normalization and continuity make $\|k-g_{\text{Eucl}}\|_\infty$ arbitrarily small on a neighborhood of the point. Estimate~\eqref{eq:harmonic-scalar-control} gives~\eqref{eq:local-small-multiplier} for $k$ with arbitrarily small leading coefficient. This proves the local claim with $w=1$. Under conditions~\ref{item:tame-reference}--\ref{item:tame-decay} of Definition~\ref{def:H1-tame} and the assumed nonnegative scalar density, the local condition just proved makes the metric $H^1$-tame. All coordinate changes used here are local, so the exterior chart and the integral decay hypotheses for computing the mass are those originally prescribed. This completes the proof of Proposition~\ref{thm:dini-harmonic-mass}.
\end{proof}

\section{Rigidity and global Euclidean coordinates}\label{sec:rigidity}

\subsection{From a parallel spinor to global Euclidean coordinates}
\paragraph{The closed coframe.}
The zero-mass identity gives a nonzero parallel spinor whose quaternionic companion determines a closed orthonormal coframe. We integrate this coframe locally and globally, then identify the resulting Euclidean and integral distances.

\begin{lemma}\label{lem:parallel-spinor-closed-coframe}
Let $M$ be a connected oriented smooth three-manifold with a spin structure, and let $g\in H^1_{\mathrm{loc}}$ satisfy locally uniform spectral bounds relative to a smooth reference metric. A nonzero spinor $\Psi\in L^2_{\mathrm{loc}}$ that is parallel in the sense of distributions has constant positive norm. It determines a global oriented orthonormal coframe
$\alpha=(\alpha^1,\alpha^2,\alpha^3)$ whose coefficients are locally in $H^1\cap L^\infty$ and satisfy $d\alpha^a=0$ in the sense of distributions.
\end{lemma}
\begin{proof}
In smooth reference spin frames the equation is $d\Psi=-\Omega\Psi$, with $\Omega\in L^2_{\mathrm{loc}}$. Its right-hand side is locally integrable, hence $\Psi\in W^{1,1}_{\mathrm{loc}}$. For $\epsilon>0$, the Sobolev chain rule and skew-Hermiticity give
\be
 d\sqrt{|\Psi|^2+\epsilon^2}
 =-\frac{\operatorname{Re}\langle\Psi,\Omega\Psi\rangle}
          {\sqrt{|\Psi|^2+\epsilon^2}}=0.
\ee
Connectedness makes $|\Psi|$ a constant $c>0$. The parallel equation now gives $\Psi\in H^1_{\mathrm{loc}}\cap L^\infty$.

The quaternionic conjugation
$\mathsf J(z_1,z_2)=(-\overline z_2,\overline z_1)$
is parallel and commutes with spin transitions. Thus $c^{-1}\Psi$ and its conjugate form a global orthonormal spin frame. We define $\alpha^a$ by requiring their Clifford actions in this frame to be the three fixed skew-Hermitian Pauli matrices. Both the spin frame and Clifford multiplication transform by the same spin transition, so the one-forms agree on overlaps and
$g=\sum_a\alpha^a\otimes\alpha^a$.

Locally we write $U=c^{-1}(\Psi,\mathsf J\Psi)\in SU(2)$, let $\theta$ be a physical orthonormal coframe, and put
$S=\operatorname{rot}(U^{-1})$, so $\alpha=S\theta$.
The bounded matrix chain rule gives $S\in H^1_{\mathrm{loc}}\cap L^\infty$ and the parallel equation gives $dU=-\Omega U$.
The exterior product rule
\be
 d(S\theta)=dS\wedge\theta+S\,d\theta
\ee
holds by approximating $S$ in $H^1$ and almost everywhere by uniformly bounded smooth matrices. The first product converges in $L^2$ because $\theta$ is bounded, and the second in $L^1$ by dominated convergence. The equation for $U$ and the torsion equation cancel these terms. Hence $d\alpha=0$; the bounded $H^1$ product rule also gives the asserted coefficient regularity.
\end{proof}

\paragraph{Local coordinates and global reconstruction.}
In an oriented coordinate ball $B$ the coframe just obtained satisfies $\alpha\in H^1_{\mathrm{loc}}\cap L^\infty$, $d\alpha=0$, and $\det(\alpha_i^a)=\sqrt{\det g}\ge c_B>0$ on each smaller ball, by ellipticity. These are precisely the hypotheses needed here for a Cartan--Whitney presentation on a Euclidean domain. Its component one-forms are flat forms because their coefficients and distributional exterior derivatives are bounded. The orientation bound is the lower bound on their wedge product.

We apply Heinonen--Keith \cite{HeinonenKeith} on a relatively compact subball. Their theorem gives a bi-Lipschitz map $F$ near each point and bounds $\|dF-\alpha\|_{L^\infty(B_r)}$ by $r\|d\alpha\|_{L^\infty(B_r)}$ on smaller concentric balls. Here the right side is zero. Hence $dF=\alpha$ almost everywhere. Every local primitive agrees with $F$ up to translation on connected overlaps, since their difference has zero weak derivative. Local invertibility in this exact case also follows from Heinonen--Kilpel\"ainen \cite{HeinonenKilpelainen}. Thus the closed coframe already has the regularity needed for reconstruction once the mass identity has furnished a parallel spinor. 

\begin{proposition}[Global rigidity from a spinor parallel in the sense of distributions]
\label{prop:global-form-rigidity}
Let $M$ be a connected oriented smooth three-manifold without boundary with a complete smooth reference metric $h$. Let $g\in H^1_{\mathrm{loc}}$ satisfy uniform spectral bounds relative to $h$ and admit a nonzero spinor in $L^2_{\mathrm{loc}}$ that is parallel in the sense of distributions. Then $M$ carries local flat coordinates for $g$ whose transitions are translations, and their complete flat manifold is a quotient of $\R^3$ by a discrete group of translations. If, for one $o\in M$,
\be
 \limsup_{R\to\infty} R^{-2}\operatorname{vol}_h(B_h(o,R))=+\infty ,
\ee
then there is an
orientation-preserving homeomorphism $\Phi:M\to\R^3$, globally
bi-Lipschitz with respect to the smooth reference metric, such that
$g=\Phi^*g_{\text{Eucl}}$ almost everywhere. In smooth local charts both $\Phi$ and its inverse belong to $W^{2,2}_{\mathrm{loc}}$. Moreover $M$ is diffeomorphic to $\R^3$.
\end{proposition}
\begin{proof}
{\bf 1. \sl The local Euclidean metric.} The preceding construction gives a global closed orthonormal coframe
$\alpha$, whose local primitives are bi-Lipschitz coordinates.
Their transitions on connected overlaps are translations. These
charts define a new smooth structure and a smooth flat metric $q$
on the same underlying topological manifold. The new structure is
locally bi-Lipschitz equivalent to the original one; smoothness of
$q$ here refers to the new atlas. The local area formula and the Sobolev chain rule give
\be
 dv_q=dv_g,\qquad C^{-1}d_h\le d_q\le C d_h.
\ee
The derivative and inverse-derivative bounds relative to $h$ give uniform local Lipschitz comparisons. Subdividing curves among these neighborhoods gives the global comparison between the distances of the two smooth metrics in their respective atlases. Completeness of $h$ then implies completeness of $q$.

\vskip.2cm \noindent
{\bf 2. \sl Global reconstruction.} The universal Riemannian cover of the complete flat manifold
$(M,q)$ is Euclidean space. The globally defined parallel coframe
$\alpha$ shows that its deck transformations preserve the coordinate
coframe, and hence are translations. If a nonzero translation of
length $a$ occurred, the cyclic intermediate quotient would have a
fundamental slab of width $a$ and volume growth at most $CR^2$.
The covering from the cyclic quotient onto $M$ maps each ball onto the ball of the same radius about its projected center, by lifting minimizing geodesics. As it is a local isometry, the area formula bounds the volume of the image by the volume of the ball upstairs. On the other hand $B_h(o,R/C)\subset B_q(o,R)$ and $dv_q=dv_g\ge c\,dv_h$. The volume hypothesis therefore excludes every nonzero translation. A Euclidean reference end satisfies the hypothesis because, for all large $R$, the ball $B_h(o,R)$ contains an exterior Euclidean ball of radius proportional to $R$. We may choose its center along a fixed ray at radius $R/2$ and its radius $R/4$.
Thus the deck group is trivial. Its Euclidean developing map gives
the asserted global $\Phi$, and the distance comparison gives its
global bi-Lipschitz bounds. The local differential is $\alpha$,
so $g=\Phi^*g_{\text{Eucl}}$ almost everywhere. On each relatively compact chart $D\Phi=\alpha\in H^1$ and its inverse matrix belongs to $H^1$ by the bounded matrix inverse rule. The derivative of $\Phi^{-1}$ equals $(D\Phi)^{-1}\circ\Phi^{-1}$ almost everywhere. Composition with this bi-Lipschitz map preserves local $H^1$, by smooth approximation and change of variables. This proves the asserted second derivative regularity in both directions.

Finally the resulting homeomorphism to $\R^3$ implies diffeomorphism of the original smooth three-manifold with $\R^3$, by uniqueness of smooth structures in dimension three~\cite{MoiseTopology}.
\end{proof}

\paragraph{The end coordinates and the integral distance.}
We suppose that the reference metric has a Euclidean end with the two quadratic metric integrals of Section~\ref{sec:end-integrals}. The global map $z=\Phi(x)$ sends the end onto the complement of the compact image of the core. Both coordinate systems satisfy condition (1) of Proposition~\ref{thm:absolute-integral-coordinate-invariance}, the second with metric $g_{\text{Eucl}}$. The proposition yields a unique $Q\in O(3)$ such that
\be
 D\Phi-Q\in L^6(\Mext),\qquad D^2\Phi\in L^2(\Mext),\qquad
 (D\Phi-Q)/r\in L^2(\Mext).
\ee
The last assertion follows by the cutoff and Hardy estimate in that proof. In particular, Lemma~\ref{lem:dyadic-annular-mass} identifies the annular limit in the original chart with zero. This converse uses only the coordinate and integrability hypotheses.

For a measurable metric, the length of a particular curve can depend on the values chosen for the metric on a set of volume zero. De Cecco and Palmieri~\cite{DeCeccoPalmieri} introduced an integral distance determined by the energies of test functions, which depend only on the metric almost everywhere. We use the formulation recalled by Lee, Litzinger, and Simon~\cite{LeeLitzingerSimon}. For $p>3$, with $\operatorname{Lip}(M)$ defined using the smooth reference distance $d_h$, we set
\be
 d_p(x,y)\coloneqq\sup\left\{|f(x)-f(y)|:
 f\in\operatorname{Lip}(M),\ \int_M|df|_g^p\,dv_g\le1\right\}.
\ee
The integral distance is the limit $d_g(x,y)\coloneqq\lim_{p\to\infty}d_p(x,y)$. For the metrics in the following lemma, its proof also establishes existence of this limit.

\begin{lemma}[Identification of the integral distance]\label{prop:integral-distance-identification}
Whenever $\Phi:(M,d_h)\to\R^3$ is a global bi-Lipschitz homeomorphism and $g=\Phi^*g_{\text{Eucl}}$ almost everywhere, the integral distance equals
\be
                         d_g(x,y)=|\Phi(x)-\Phi(y)|.
\ee
\end{lemma}
\begin{proof}
The global bi-Lipschitz change of variables and $g=\Phi^*g_{\text{Eucl}}$
identify this supremum with the corresponding Euclidean supremum:
for $f=F\circ\Phi$, the weak chain and area formulas give
\be
                   \int_M|df|_g^p\,dv_g
                              =\int_{\R^3}|dF|^pdx.
\ee
Both directions preserve the Lipschitz test class.

We also verify that the Euclidean suprema tend to the Euclidean distance.
For two points at distance $d>0$, the function
$F(z)=\min\{|z-x|,d\}$ has gradient norm
$|B_d|^{1/p}$ and point difference $d$. It proves the lower limit
is at least $d$. For the upper limit we take $p_j\to\infty$ and
nearly maximizing tests $F_j$, normalized by $F_j(x)=0$ and
$\|dF_j\|_{p_j}\le1$. For every fixed $q>3$ and bounded ball $B$,
\be
 \|dF_j\|_{L^q(B)}\le |B|^{1/q-1/p_j}\qquad(p_j\ge q).
\ee
Local Morrey compactness and the normalization yield a subsequence
converging locally uniformly. A diagonal extraction over balls
centered at $x$ and integer exponents $q>3$, followed by weak lower
semicontinuity, gives a limit $F$ satisfying
$\|dF\|_{L^q(B)}\le |B|^{1/q}$ for every finite $q>3$.
Letting $q\to\infty$ shows $|dF|\le1$ almost everywhere.
Thus $F$ is one-Lipschitz, and the limiting point difference is at
most $d$. Applying this to a subsequence attaining the upper limit
proves the upper bound and the claimed Euclidean limit. Pullback
by $\Phi$ proves the assertion on $M$.
\end{proof}

\subsection{Intrinsic local multiplier conditions}\label{sec:sobolev-multiplier-concentration}
We justify the reference comparisons used in the two sufficient criteria of the introduction. On every relatively compact domain, the condition $\kappa_h(g;x)=0$ at each point is equivalent to the multiplier estimate with an arbitrarily small coefficient of the Dirichlet norm. One implication follows by testing on shrinking balls and using Poincar\'e's inequality. Conversely, we cover a compact set by finitely many balls with multiplier coefficient at most $\eta>0$ and choose a smooth quadratic partition $\sum_i\chi_i^2=1$ there. Summing the squared estimates for $\chi_i v$ gives
\be
 \||\nabla^h g|_h v\|_{L^2(h)}^2
 \le\eta^2\bigl(\|dv\|_{L^2(h)}^2+C_\chi\|v\|_{L^2(h)}^2\bigr).
\ee
For another smooth reference $\bar h$, norms and volumes are comparable on compact neighborhoods, and
$|\nabla^{\bar h}g|_{\bar h}\le C(|\nabla^h g|_h+1)$.
Poincar\'e's inequality on shrinking balls therefore gives
$\kappa_{\bar h}(g;x)\le C\kappa_h(g;x)$; interchanging the references gives the reverse comparison. Finiteness and vanishing are consequently intrinsic. The numerical values can differ; for $\bar h=c^2h$ they satisfy $\kappa_{\bar h}(g;x)=c^{-2}\kappa_h(g;x)$.

The quantity $\mathcal K(g;x)$ in~\eqref{eq:metric-concentration} is independent of the smooth connection used in its definition. Changing this connection adds a locally bounded tensor, whose multiplier contribution on a reference ball of radius $r$ is $O(r)$ by Poincar\'e's inequality and the spectral bounds. With $\nabla^0=\nabla^h$, formula~\eqref{eq:trace-free-shape-derivative} and the spectral bounds compare $\mathcal K(g;x)$ in both directions with $\kappa_h(b;x)$, where $b$ is the volume-normalized metric. The comparison constants depend only on the spectral bounds. These facts prove both the asserted reference independence and the intrinsic interpretation of the first sufficient criterion.

\subsection{Proof of the main theorem}\label{sec:main-proofs}
\paragraph{Existence and invariance of the mass.}
We assume that $g$ is an $H^1$-tame metric with $R_g^{\mathrm w}\ge0$, as assumed in Theorem~\ref{thm:small-concentration-mass}. Conditions~\ref{item:tame-reference}--\ref{item:tame-decay} of Definition~\ref{def:H1-tame}, Proposition~\ref{prop:automatic-end-reference}, Proposition~\ref{prop:annular-negative-variation} and Lemma~\ref{lem:dyadic-annular-mass} give the existence and finiteness criterion for $m_{\rm ADM}=\lim m_{\text{dya}}(R)$. Proposition~\ref{thm:absolute-integral-coordinate-invariance} gives independence of the admissible Euclidean reference. Conversely, if $\Phi$ has the properties stated in Theorem~\ref{thm:small-concentration-mass}, restricting $\Phi^{-1}$ to the exterior of a sufficiently large ball gives an admissible reference in which the metric is $g_{\text{Eucl}}$. Invariance identifies the original mass with zero. If the scalar measure has infinite mass on the end, the mass is $+\infty$. We may therefore assume below that it is finite.

\paragraph{The oriented case.}
We choose a spin structure, with its standard restriction to the Euclidean end. In each factorization in condition~\ref{item:tame-local} of Definition~\ref{def:H1-tame}, the strict multiplier bound for $k$ gives the local Dirac estimate of Section~\ref{sec:dirac-form-local}. The coordinate and conformal transport and the graph approximation of Proposition~\ref{prop:local-conformal-charts} give the physical domain required by the homogeneous Witten construction. For a unit terminal spinor, Proposition~\ref{prop:bounded-Witten-spinor} and the mass identity of Section~\ref{sec:mass-identity} therefore give
\be
 4\pi m_{\rm ADM}
 =\int_M|\nabla^g\Psi|_g^2\,dv_g
      +\tfrac14\int_M|\widetilde\Psi|^2\,d\mu\ge0.
\ee
At zero mass, $\Psi$ is parallel and has norm one, since its difference from the terminal spinor belongs to $L^6$ on the end. Lemma~\ref{lem:parallel-spinor-closed-coframe} gives the closed $H^1_{\mathrm{loc}}$ coframe. Proposition~\ref{prop:global-form-rigidity} gives the global bi-Lipschitz Euclidean coordinates and the $W^{2,2}_{\mathrm{loc}}$ regularity of both maps. Lemma~\ref{prop:integral-distance-identification} identifies the integral distance, and the coordinate comparison above verifies admissibility at infinity.

\paragraph{The orientation cover and its two ends.}
If $M$ is nonorientable, we denote its connected orientation double cover by $\pi:\widehat M\to M$ and put $\widehat h=\pi^*h$, $\widetilde g=\pi^*g$. The reference $\widehat h$ is complete and has two Euclidean ends $\Mext^+$ and $\Mext^-$, each carrying the lifted metric of the original end. The analytic conditions~\ref{item:tame-regularity}--\ref{item:tame-local} of Definition~\ref{def:H1-tame} and the assumed nonnegativity of the scalar curvature density lift unchanged: the spectral bounds and the condition that the scalar curvature density is nonnegative hold locally, the tensor decay holds on each end, and every local factorization has the same representative $k$ and multiplier bound. We use these lifted conditions directly on the two-ended cover. In particular, we do not impose the single-end requirement of condition~\ref{item:tame-reference} on $\widehat M$.

Lemma~\ref{prop:global-scalar-Sobolev} applies to the connected two-ended cover and gives the global Sobolev inequality. Proposition~\ref{prop:local-conformal-charts} gives the physical self-adjoint Dirac domain and zero $L^2$ kernel there. Thus the hypotheses of the finite-ended Witten construction hold on $\widehat M$.

We prescribe a nonzero terminal spinor $\xi$ on $\Mext^+$ and zero on $\Mext^-$. A smooth reference spinor with these values has finite covariant energy by the derivative integral, and finite scalar energy by the finiteness criterion. The finite-end version of Proposition~\ref{prop:bounded-Witten-spinor} applies with these two terminal values. The exhaustion gives the terminal contribution $4\pi m_{\rm ADM}|\xi|^2$ on the first end and zero on the second, so
\be
 4\pi m_{\rm ADM}|\xi|^2
 =\int_{\widehat M}|\widehat\nabla\Psi|^2\,dv_{\widetilde g}
       +\tfrac14\int_{\widehat M}|\widetilde\Psi|^2\,d\widehat\mu
 \ge0.
\ee
At zero mass the spinor would be parallel and have constant norm on the connected cover. Since both ends have infinite volume,
$\Psi-\xi\in L^6(\Mext^+)$ forces that norm to be $|\xi|$, whereas
$\Psi\in L^6(\Mext^-)$ forces it to vanish. Thus the mass is strictly positive in the nonorientable case. This completes the proof of Theorem~\ref{thm:small-concentration-mass}.

\vfill

\paragraph{Acknowledgments.}
I am deeply grateful to Bruno Le Floch (LPTHE, Sorbonne Univ.) for many discussions that have contributed to the development of this work.

\paragraph{Funding.}
This work was supported in part by the Agence Nationale de la Recherche through Einstein-PPF (ANR-23-CE40-0010-02) and by the European Union through the Marie Sk{\l}odowska-Curie Staff Exchanges project Einstein-Waves (grant 101131233).

\paragraph{Use of AI.}
 ChatGPT~6 was used to reorganize the original manuscript, improve the flow of its presentation. The author takes responsibility for the final text.

\paragraph{Competing interests and data availability.}

The author declares no competing interests. This work is theoretical and has no associated datasets.


\end{document}